# HYPERBOLIC CLASSIFICATION OF NATURAL NUMBERS AND GOLDBACH CONJECTURE

BY FERNANDO REVILLA
Department of Mathematics. I.E.S. Sta Teresa de Jesús.
Comunidad de Madrid. Spain
frej0002@ficus.pntic.mec.es

ABSTRACT In a letter to Euler in 1742, Goldbach stated that "Every even integer greater than 2 is the sum of two prime numbers". This statement, known as the Goldbach Conjecture, was published in Edward Waring's *Meditationes Algebraicae* in England in 1770, though without proof. To this date, all efforts to prove this statement have been unsuccessful. In this paper we expose a theory which provides a new perspective on natural numbers ( Hyperbolic Classification of Natural Numbers ). This theory leads to a theorem that characterizes the conjecture in terms of a structure isomorphic to $(\mathbb{N},+,\cdot)$. This characterization depends on the way we arrange natural numbers on the real line.

## Introduction

For a natural number $n>1$ the fact of being a prime is equivalent to stating that the hyperbola $xy=n$ does not contain non-trivial natural number coordinate points that is, the only natural coordinate points in the hyperbola are $(1,n)$ and $(n,1)$. The first part of this paper establishes a family of bijective functions between non-negative real numbers and a half-open interval of real numbers. Bijectivity allows us to transport usual real number operations, sum and product, to the interval. It also allows us to deform the $xy=k$ hyperbolas with $k$ as a real number in such a way that we can distinguish whether a natural number $n$ is a prime or not by its behaviour in terms of gradients of the deformed hyperbolas near the deformed $xy=n$ (Hyperbolic Classification of Natural Numbers). In the second part of this paper, using the concept of essential regions associated to a hyperbola, we construct the area function for the transformed region of the union of regions $R \equiv x \geq 2$, $y \geq x$, $xy \leq k$, $S \equiv x \geq 2$, $y \geq x$, $\alpha - k \leq xy \leq \alpha - 4$ for a given even number $\alpha$ $(\alpha \geq 16)$ when $k$ covers the interval $[4,\frac{\alpha}{2}]$. This is a class 2 function except for a finite set of points. Using the concept of essential point, linked with the second derivative, we obtain a theorem which characterizes that $\alpha$ is the sum of two prime numbers, in terms of the multiplicity of the essential points. The third part proves that for some elements of the aforementioned family of bijective functions, the second derivative is continuous in the interval $[4,\frac{\alpha}{2}]$ (Goldbach Conjecture function). Using this function we prove a statement equivalent to that of the Goldbach Conjecture in the infinite set

$$\mathfrak{N} = \left\{ \alpha \in \mathbb{N} : (\alpha \ even) \wedge (\alpha \geq 16) \wedge (\frac{\alpha}{2} \ non\text{-}prime) \ \wedge (\alpha - 3 \ non\text{-}prime) \right\}$$

that depends on the way we arrange natural numbers on the real line.

*Contents*

# 1. HYPERBOLIC CLASSIFICATION OF NATURAL NUMBERS

## 1.1 $\mathbb{R}^+$ coding function

**1.1.0 Introduction** In this chapter, we define a function $\psi$ which ranges from non-negative real numbers to a half-open interval, strictly increasing, continuous and class 1 in each interval $[m, m+1]$ $(m \in \mathbb{N} = \{0,1,2,3,...\})$. The bijectivity of $\psi$ allows to transport the usual sum and product of $\mathbb{R}^+$ to the set $\widehat{\mathbb{R}^+} := \psi(\mathbb{R}^+)$ in the usual manner. That is, calling $\hat{x} = \psi(x)$, we define $\hat{s} \oplus \hat{t} = \widehat{s+t}$, $\hat{s} \otimes \hat{t} = \widehat{st}$. Therefore, $(\widehat{\mathbb{R}^+}, \oplus, \otimes)$ is an algebraic structure isomorphic to the usual one $(\mathbb{R}^+, +, \cdot)$ and, as a result, we obtain an algebraic structure $(\widehat{\mathbb{N}^+} := \psi(\mathbb{N}), \oplus, \otimes)$ isomorphic to the usual one, $(\mathbb{N}, +, \cdot)$. Thus we transport the notation from $\mathbb{R}^+$ to $\widehat{\mathbb{R}^+}$, that is $\hat{n}$ is natural iff $n$ is natural, $\hat{p}$ is prime iff $p$ is prime, $\hat{x}$ is rational iff $x$ is rational, etc. The function $\psi$ also preserves the usual orders. Assume that, for example $\hat{0} = 0$, $\hat{1} = 0'72$, $\hat{2} = 1'3$, $\hat{3} = 3'0001$, $\hat{4} = \pi$, $\hat{5} = 6'3$, $\hat{6} = 7'2$, $\hat{7} = 7'21, \ldots$ ,$\widehat{12} = 9'0\hat{3}$, … the following situation would arise:

*The even number* $9'0\widehat{3}$ *is the "sum" of the " prime numbers"* $6'3$ *and* $7'21$, *and the number* $\pi$ *is the "product" of the numbers 0'72 and* $\pi$.

Obviously, until now, we have added nothing relevant to the Goldbach Conjecture, we have only actually changed the symbolism by means of the function $\psi$. If we call $\hat{x}\hat{y}$ plane, the set $\left(\psi(\mathbb{R}^+)\right)^2$, the hyperbolas $xy = k(k>0)$ of the $xy$ plane are transformed, by means of the function $\psi\times\psi$, at the $\hat{x}\otimes\hat{y}=\hat{k}$ "hyperbolas" of the $\hat{x}\hat{y}$ plane. We will restrict our attention to the points in the $\hat{x}\hat{y}$ plane that satisfy $\hat{x}\geq\hat{1}$ and $\hat{y}\geq\hat{x}$. Then, with these restrictions for the $\psi$ function, it is possible to choose right-hand and left-hand derivatives of $\psi$ where $m\in\mathbb{N}^*=\{1,2,3,...\}$ such that we can characterize the natural number coordinate points in the $\hat{x}\hat{y}$ plane in terms of differentiability of the functions which determine the transformed hyperbolas. As a result, we can distinguish prime numbers from natural non-prime numbers in the aforementioned terms (Hyperbolic Classification of Natural Numbers).

**1.1.1 Definition** Let $\psi:\mathbb{R}^+\to\mathbb{R}$ be a map and let $\psi_m$ be the restriction of $\psi$ to each closed interval $[m,m+1](m\in\mathbb{N})$. We say that $\psi$ is an $\mathbb{R}^+$ coding function iff:

i) $\psi(0)=0$ ii) $\psi\in\mathcal{C}(\mathbb{R}^+)$ iii) $\forall m\in\mathbb{N}, \psi_m\in\mathcal{C}^1([m,m+1])$ with positive derivative in $[m,m+1]$.

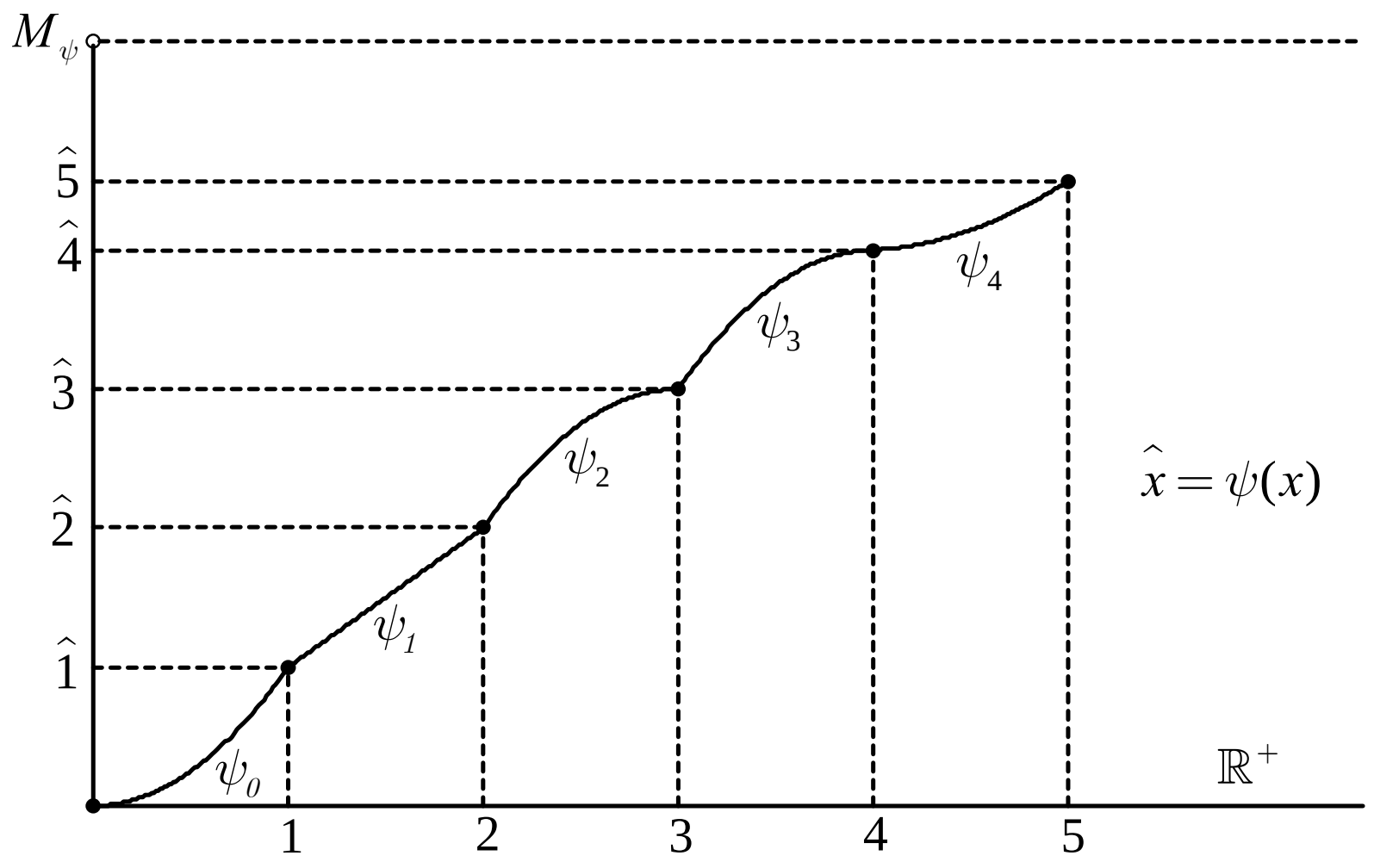

(Fig. 1)

Easily proved, if $\psi$ is an $\mathbb{R}^+$ coding function then it is strictly increasing and consequently, injective. If $M_\psi := \sup\{f(x): x\in\mathbb{R}^+\}$ then, $M_\psi\in(0,+\infty]$, being $M_\psi=+\infty$ iff $\psi$ is not bound, and so $\psi(\mathbb{R}^+)=[0,M_\psi)$. Therefore,

$$\psi:\mathbb{R}^+\to\psi(\mathbb{R}^+)=[0,M_\psi)$$

is bijective, and here onwards we will refer to the $\psi$ function as a bijective function. We will frequently use the notation $\hat{x}=\psi(x)$. Due to the $\psi$ bijection, we transport the sum and the product from $\mathbb{R}^+$ to $\left[0,M_\psi\right)$ in the usual manner, that is we define in $\left[0,M_\psi\right)$ the operations:

a) Sum: $\hat{x}\oplus\hat{y}=\psi(x+y)=\widehat{x+y}$

b) Product: $\hat{x}\otimes\hat{y}=\psi(x\cdot y)=\widehat{x\cdot y}$

c) Subtraction: $\hat{x}\sim\hat{y}=\psi(x-y)=\widehat{x-y}\ (x\geq y)$

d) Quotient: $\hat{x}\div\hat{y}=\psi\left(\frac{x}{y}\right)=\widehat{\left(\frac{x}{y}\right)}\ (y\neq 0)$

Thus, $\psi:\left(\mathbb{R}^+,+,\cdot\right)\to\left(\left[0,M_\psi\right),\oplus,\otimes\right)$ is an isomorphism. The $\psi$ function preserves the usual orders, that is, $\hat{s}\leq\hat{t}\Leftrightarrow s\leq t,\ \hat{s}=\hat{t}\Leftrightarrow s=t$. For $\hat{x}\in\left[0,M_\psi\right)$ we say that $\hat{x}$ is a natural number iff $x$ is a natural number, $\hat{x}$ is prime iff $x$ is prime, $\hat{x}$ is rational iff $x$ is rational, etc. When we work on the set $\left[0,M_\psi\right)^2$, we say that we are on the $\hat{x}\hat{y}$ plane.

## 1.2 Natural coordinate points in the $\hat{x}\hat{y}$ plane

Let $\alpha\in\mathbb{N}^*$. We want to characterise the $(u,v)$ natural numberr coordinate points of the $\hat{x}\hat{y}$ plane whose coordinates sum is $\hat{\alpha}$ ( $u\oplus v=\hat{\alpha}$ ). For this, we begin with the function

$$f_\alpha:[0,\alpha]\to[0,\alpha],\ f_\alpha(x)=\alpha-x.$$

Let $\psi$ be an $\mathbb{R}^+$ coding function and let us apply the function $\psi\times\psi$ to

$$\Gamma(f_\alpha)=\left\{(x,y)\in(\mathbb{R}^+)^2:x\in[0,\alpha]\wedge y=f(x)\right\}$$

We then obtain the transformed curve:

$$\Gamma\left(\widehat{f_\alpha}\right)=(\psi\times\psi)(\Gamma(f_\alpha))=\left\{\left(\hat{x},\hat{y}\right)\in\left[0,M_\psi\right)^2:\hat{x}\in\left[0,\hat{\alpha}\right]\wedge\hat{y}=\widehat{f_\alpha}\left(\hat{x}\right)\right\}$$

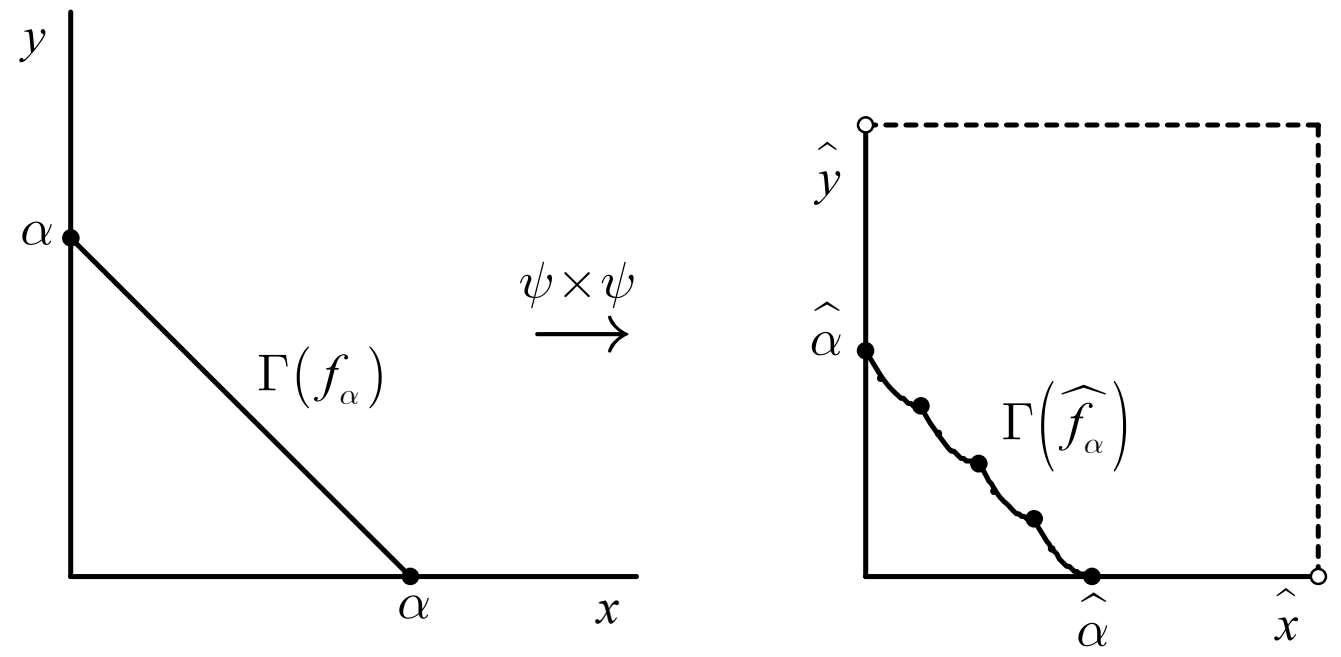

(Fig. 2)

The $\widehat{f_\alpha}:[0,\hat{\alpha}]\to[0,\hat{\alpha}]$ function which determines the graph of the transformed curve is

$$\widehat{f_\alpha}(u)=\psi(\alpha-\psi^{-1}(u))$$

Of course, $(u,v)\in\Gamma(\widehat{f_\alpha})$ has natural number coordinates iff $u$ is natural. The following theorem will allow a characterisation of the natural coordinate points whose sum is $\hat{\alpha}$.

**1.2.1 Theorem** Let $\alpha \in \mathbb{N}^*$, $\psi:\mathbb{R}^+ \to [0, M_\psi)$ an $\mathbb{R}^+$ coding function and

$$\widehat{f_\alpha}:[0,\hat{\alpha}]\to[0,\hat{\alpha}], \widehat{f_\alpha}(u)=\psi(\alpha-\psi^{-1}(u))$$. Then,

a) $\widehat{f_\alpha}$ is continuous and strictly decreasing.

b) Let $m \in \mathbb{N}:[\hat{m},\hat{m}\oplus\hat{1}]\subset[0,\hat{\alpha}]$. Then $\widehat{f_\alpha}$ is differentiable in the interval $(\hat{m},\hat{m}\oplus\hat{1})$ with derivative

$$\left(\widehat{f_\alpha}\right)'(u)=-\frac{(\psi_{\alpha-1-m})'\left(\alpha-\psi^{-1}(u)\right)}{(\psi_m)'\left(\psi^{-1}(u)\right)}$$

c) $\forall m=0,1,2,...,\alpha-1$, we verify

$$\left(\widehat{f_\alpha}\right)'_+\left(\widehat{m}\right)=-\frac{(\psi_{\alpha-m-1})'_-(\alpha-m)}{(\psi_m)'_+(m)}$$

d) $\forall m=1,2,3,...,\alpha$, we verify

$$\left(\widehat{f_\alpha}\right)'_-\left(\widehat{m}\right)=-\frac{(\psi_{\alpha-m})'_+(\alpha-m)}{(\psi_{m-1})'_-(m)}$$

<u>*Proof*</u> a) We have $\psi^{-1}:[0,\widehat{\alpha}]\to[0,\alpha], f_\alpha:[0,\alpha]\to[0,\alpha], \psi:[0,\alpha]\to[0,\widehat{\alpha}]$. Therefore, $\widehat{f_\alpha}=\psi\circ f_\alpha\circ\psi^{-1}$ is a composition of continuous functions, as a result it is continuous. In addition:

$$0\le s<t\le\hat{\alpha}\Rightarrow\psi^{-1}(s)<\psi^{-1}(t)\Rightarrow\alpha-\psi^{-1}(s)>\alpha-\psi^{-1}(t)\Rightarrow$$
$$\psi\left(\alpha-\psi^{-1}(s)\right)>\psi\left(\alpha-\psi^{-1}(t)\right)\Rightarrow\widehat{f_\alpha}(s)>\widehat{f_\alpha}(t)\Rightarrow$$

$\widehat{f_\alpha}$ is strictly decreasing.

b) $(\hat{m},\hat{m}\oplus\hat{1})\xrightarrow{\psi^{-1}}(m,m+1)\xrightarrow{f_\alpha}(\alpha-m-1,\alpha-m)\xrightarrow{\psi}(\hat{\alpha}\sim\hat{m}\sim\hat{1},\hat{\alpha}\sim\hat{m})$ In other words, $\widehat{f_\alpha}$ maps $\widehat{f_\alpha}:(\hat{m},\hat{m}\oplus\hat{1})\to(\hat{\alpha}\sim\hat{m}\sim\hat{1},\hat{\alpha}\sim\hat{m})$. For the $\widehat{f_\alpha}$ function, all the hypotheses of the compound and inverse derivative theorems are fulfilled, and $\forall u\in\left(\widehat{m},\widehat{m}\oplus\hat{1}\right)$:

$$\left(\widehat{f_\alpha}\right)'(u)=\psi'\left(\alpha-\psi^{-1}(u)\right)\cdot\frac{-1}{\psi'\left(\psi^{-1}(u)\right)}=(\psi_{\alpha-m-1})'\left(\alpha-\psi^{-1}(u)\right)\cdot\frac{-1}{(\psi_m)'\left(\psi^{-1}(u)\right)}$$

c) Let $\varepsilon>0:\hat{m}<\hat{m}\oplus\varepsilon<\hat{m}\oplus\hat{1}$ and let $\delta>0:\psi^{-1}(\hat{m}\oplus\varepsilon)=m+\delta$ We obtain

$$[\hat{m},\hat{m}\oplus\varepsilon)\xrightarrow{\psi^{-1}}[m,m+\delta)\xrightarrow{f_\alpha}(\alpha-m-\delta,\alpha-m]\xrightarrow{\psi}\left(\hat{\alpha}\sim\hat{m}\sim\hat{\delta},\hat{\alpha}\sim\hat{m}\right]$$

$\widehat{f_\alpha}$ therefore maps $\widehat{f_\alpha}:[\hat{m},\hat{m}\oplus\varepsilon)\to\left(\hat{\alpha}\sim\hat{m}\sim\hat{\delta},\hat{\alpha}\sim\hat{m}\right]$

Consequently $\forall u\in[\hat{m},\hat{m}\oplus\varepsilon)$ ,$\left(\widehat{f_\alpha}\right)(u)=\psi\left(\alpha-\psi^{-1}(u)\right)=\psi_{\alpha-m-1}\left(\alpha-\psi^{-1}(u)\right)$

.Therefore: $\left(\widehat{f_\alpha}\right)'_+(\hat{m}) = (\psi_{\alpha-m-1})'_-(\alpha-m)\cdot\frac{-1}{\psi'\left(\psi^{-1}(\hat{m})\right)} = -\frac{(\psi_{\alpha-m-1})'_-(\alpha-m)}{(\psi_m)'_+(m)}$

d) We can similarly reason. Let $\varepsilon > 0 : \hat{m} \sim \hat{1} < \hat{m} \sim \varepsilon$ (or $\varepsilon < \hat{1}$). We obtain

$$(\hat{m}\sim\varepsilon,\hat{m}] \xrightarrow{\psi^{-1}} (m-\delta,m] \xrightarrow{f_\alpha} [\alpha-m,\alpha-m+\delta) \xrightarrow{\psi} \left[\hat{\alpha}\sim\hat{m},\hat{\alpha}\sim\hat{m}\oplus\hat{\delta}\right)$$

( Note that $0<\delta<1$), $\widehat{f_\alpha}$ therefore maps

$$\widehat{f_\alpha} : (\hat{m}\sim\varepsilon,\hat{m}] \to \left[\hat{\alpha}\sim\hat{m},\hat{\alpha}\sim\hat{m}\oplus\hat{\delta}\right)$$

As a result, $\widehat{f_\alpha}(u) = \psi\left(\alpha-\psi^{-1}(u)\right) = \psi_{\alpha-m}\left(\alpha-\psi^{-1}(u)\right)$

Therefore, $\left(\widehat{f_\alpha}\right)'_-(\hat{m}) = (\psi_{\alpha-m})'_+(\alpha-m)\cdot\frac{-1}{\psi'\left(\psi^{-1}(\hat{m})\right)} = -\frac{(\psi_{\alpha-m})'_+(\alpha-m)}{(\psi_{m-1})'_-(m)}$

The theorem is thus proven.

Let us now call $\begin{cases} a_k = (\psi_{k-1})'_-(k) \\ b_k = (\psi_k)'_+(k) \end{cases} \quad (k=1,2,3...)$

From the previous theorem , $\left(\widehat{f_\alpha}\right)'_+(\hat{m}) = -\frac{a_{\alpha-m}}{b_m}$ , $\left(\widehat{f_\alpha}\right)'_-(\hat{m}) = -\frac{b_{\alpha-m}}{a_m}$

Therefore, $\widehat{f_\alpha}$ is differentiable at $\widehat{m}$ $\left(\hat{m} = \hat{1},\hat{2},...,\hat{\alpha}\sim\hat{1}\right)$ iff

$$a_m \cdot a_{\alpha-m} = b_m \cdot b_{\alpha-m}$$

The $\widehat{f_\alpha}$ function is not differentiable for any $\widehat{m} \in \left\{\hat{1},\hat{2},...,\widehat{\alpha}\sim\hat{1}\right\}$ iff

$$a_m \cdot a_{\alpha-m} \neq b_m \cdot b_{\alpha-m} \quad (\forall m = 1,2,...,\alpha-1)$$

If $\widehat{f_\alpha}$ is not differentiable for any $\widehat{m} \in \left\{\hat{1},\hat{2},...,\widehat{\alpha}\sim\hat{1}\right\}$, this circumstance allows us to immediately visualise the $\Gamma_{\widehat{f_\alpha}}$ points with natural number coordinates, and from this we may see something deeper.

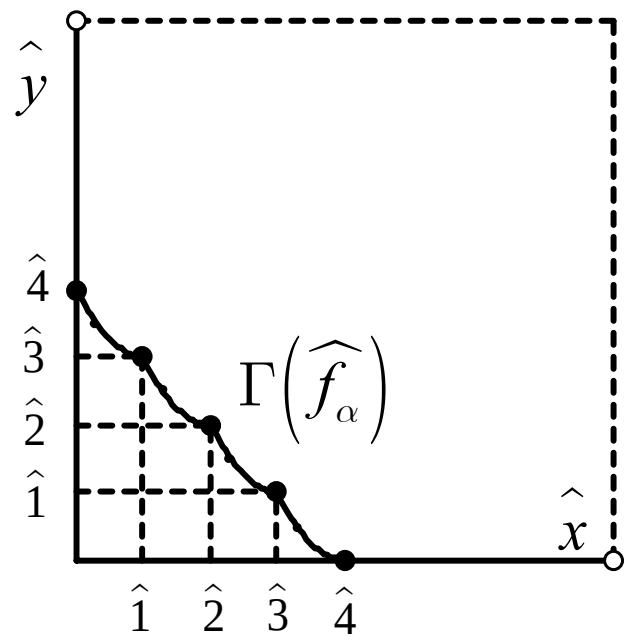

(Fig. 3)

**1.2.2 Definition** Let $\alpha \in \mathbb{N}^*$ where $\alpha \geq 2$, and $\psi : \mathbb{R}^+ \to \left[0, M_\psi\right)$ an $\mathbb{R}^+$ coding function. It is said that the $\psi$ function identifies natural numbers in $\left[\hat{0},\hat{\alpha}\right]$ iff $\forall u \in \left(\hat{0},\hat{\alpha}\right)$ "$u$ is a natural number $\Leftrightarrow \widehat{f_\alpha}$ is not differentiable at $u$"

**1.2.3 Corollary** $\alpha \in \mathbb{N}^* \ (\alpha \geq 2)$, $\psi : \mathbb{R}^+ \to [0, M_\psi)$ an $\mathbb{R}^+$ coding function. Then,

$\psi$ identifies natural numbers in $[0, \hat{\alpha}] \Leftrightarrow a_m \cdot a_{\alpha-m} \neq b_m \cdot b_{\alpha-m} \quad (\forall m = 1, 2, ..., \alpha - 1)$

## 1.3 Hyperbolas in the $\hat{x}\hat{y}$ plane

The aim here is to study the transformed curves of the $y = \frac{k}{x}\left(k \in \mathbb{R}^+ - \{0\}\right)$ hyperbolas by means of an $\mathbb{R}^+$ coding function in terms of differentiability. Consider the function

$$h_k : (0, +\infty) \to (0, +\infty), \quad h_k(x) = \frac{k}{x}$$

We call $\widehat{h_k}$ the function that determines the transformed curve graph of $\Gamma(h_k)$ by means of $\psi \times \psi$. We obtain $\widehat{h_k} : (0, M_\psi) \to (0, M_\psi)$ , $\widehat{h_k}(u) = \psi\left(\frac{k}{\psi^{-1}(u)}\right)$

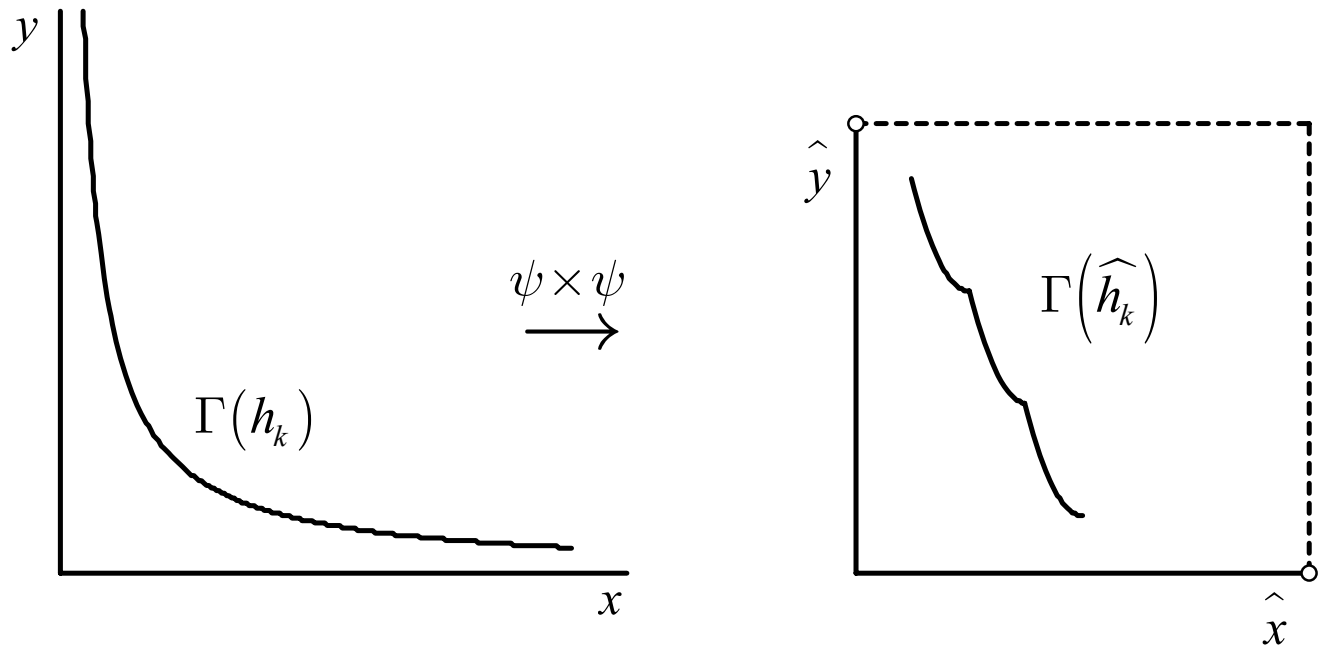

(Fig. 4)

**1.3.1 Proposition** Let $\psi : \mathbb{R}^+ \to [0, M_\psi)$ be an $\mathbb{R}^+$ prime coding. Then,

$$\widehat{h_k} : (0, M_\psi) \to (0, M_\psi) \text{ is continuous and strictly decreasing.}$$

*Proof* $(0, M_\psi) \xrightarrow{\psi^{-1}} (0, +\infty) \xrightarrow{h_k} (0, +\infty) \xrightarrow{\psi} (0, M_\psi)$ , thus $\widehat{h_k} = \psi \circ h_k \circ \psi^{-1}$ is a composition of continuous functions, and is consequently continuous. In addition,

$$0 < s < t < M_\psi \Rightarrow \psi^{-1}(s) < \psi^{-1}(t) \Rightarrow \frac{k}{\psi^{-1}(s)} > \frac{k}{\psi^{-1}(t)} \Rightarrow \psi\left(\frac{k}{\psi^{-1}(s)}\right) >$$

$$\psi\left(\frac{k}{\psi^{-1}(t)}\right) \Rightarrow \widehat{h_k}(s) > \widehat{h_k}(t) \Rightarrow \widehat{h_k} \text{ is strictly decreasing .}$$

We will now analyse the differentiability of $\widehat{h_k}$ distinguishing, for this, the cases in which the dependent and/or independent variable takes natural or non-natural number values.

*Case 1* $u \in (\hat{n}, \hat{n} \oplus \hat{1}) \ (n \in \mathbb{N})$.

We obtain $(\hat{n}, \hat{n} \oplus \hat{1}) \xrightarrow{\psi^{-1}} (n, n+1) \xrightarrow{h_k} \left(\frac{k}{n+1}, \frac{k}{n}\right) \xrightarrow{\psi} \left(\hat{k} \div (\hat{n} \oplus \hat{1}), \hat{k} \div \hat{n}\right)$, that is $\widehat{h_k}$ maps

$$\widehat{h_k} : (\hat{n}, \hat{n} \oplus \hat{1}) \to \left(\hat{k} \div (\hat{n} \oplus \hat{1}), \hat{k} \div \hat{n}\right)$$

1.a) $\widehat{h_k}(u)$ *is not a natural number.*

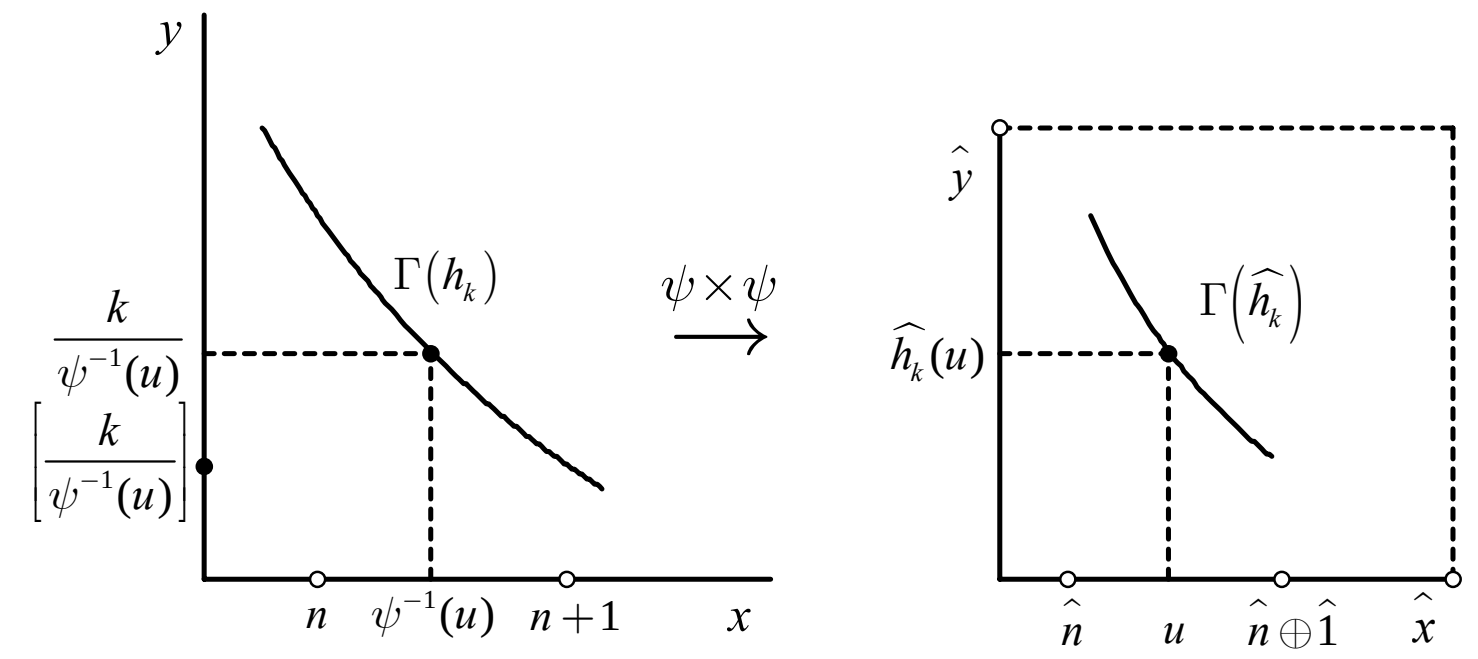

(Fig. 5)

Since $\frac{k}{\psi^{-1}(u)}$ is not a natural number, in a neighbourhood of $u$, the expression of the $\widehat{h_k}$ function is:

$$\widehat{h_k}(t) = \psi_{\left[\frac{k}{\psi^{-1}(u)}\right]}\left(\frac{k}{\psi^{-1}(t)}\right)$$

Therefore, $\left(\widehat{h_k}\right)'(u) = \left(\psi_{\left[\frac{k}{\psi^{-1}(u)}\right]}\right)'\left(\frac{k}{\psi^{-1}(u)}\right) \cdot \frac{-k}{\left(\psi^{-1}(u)\right)^2} \cdot \frac{1}{(\psi_n)'\left(\psi^{-1}(u)\right)}$

and consequently, $\widehat{h_k}$ is differentiable at $u$.

1.b) $\widehat{h_k}(u)$ *is a natural number.*

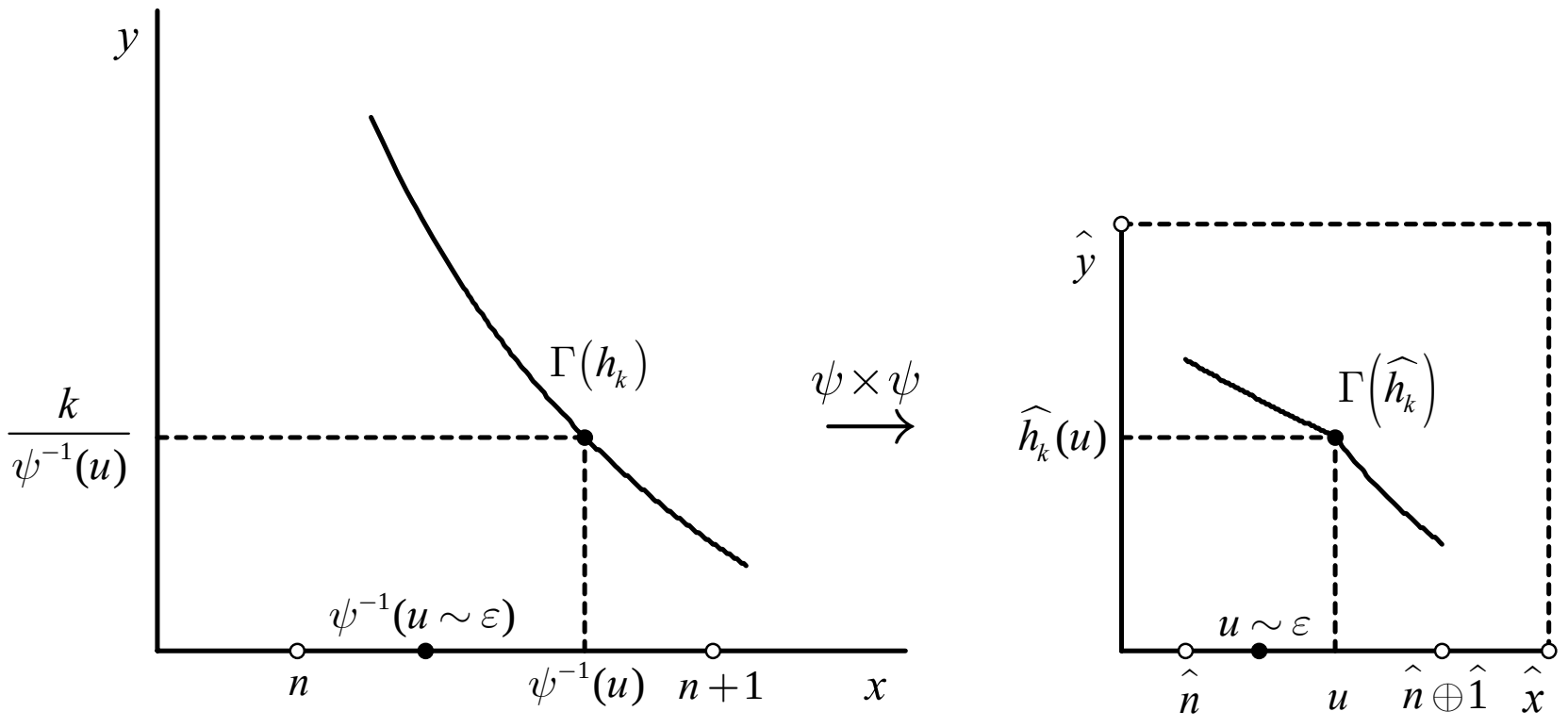

(Fig. 6)

This is equivalent to say that $\frac{k}{\psi^{-1}(u)}$ is a natural number. For a sufficiently small $\varepsilon > 0$,

$$(u \sim \varepsilon, u] \xrightarrow{\psi^{-1}} \left(\psi^{-1}(u \sim \varepsilon), \psi^{-1}(u)\right] \xrightarrow{h_k} \left[\frac{k}{\psi^{-1}(u)}, \frac{k}{\psi^{-1}(u \sim \varepsilon)}\right) \xrightarrow{\Re} \left[\hat{k} \div \widehat{\psi^{-1}(u)}, \hat{k} \div \widehat{\psi^{-1}(u \sim \varepsilon)}\right)$$

Then, $\forall t \in (u \sim \varepsilon, u]$ we verify $\widehat{h_k}(t) = \psi_{\frac{k}{\psi^{-1}(u)}}\left(\frac{k}{\psi^{-1}(t)}\right)$ and:

$$\left(\widehat{h_k}\right)'_-(u)=\left(\psi_{\frac{k}{\psi^{-1}(u)}}\right)'_+\left(\frac{k}{\psi^{-1}(u)}\right)\cdot\frac{-k}{\left(\psi^{-1}(u)\right)^2}\cdot\frac{1}{(\psi_n)'\left(\psi^{-1}(u)\right)}$$

Let us now examine the value of $\left(\widehat{h_k}\right)'_+(u)$. For a sufficiently small $\varepsilon>0$ we obtain

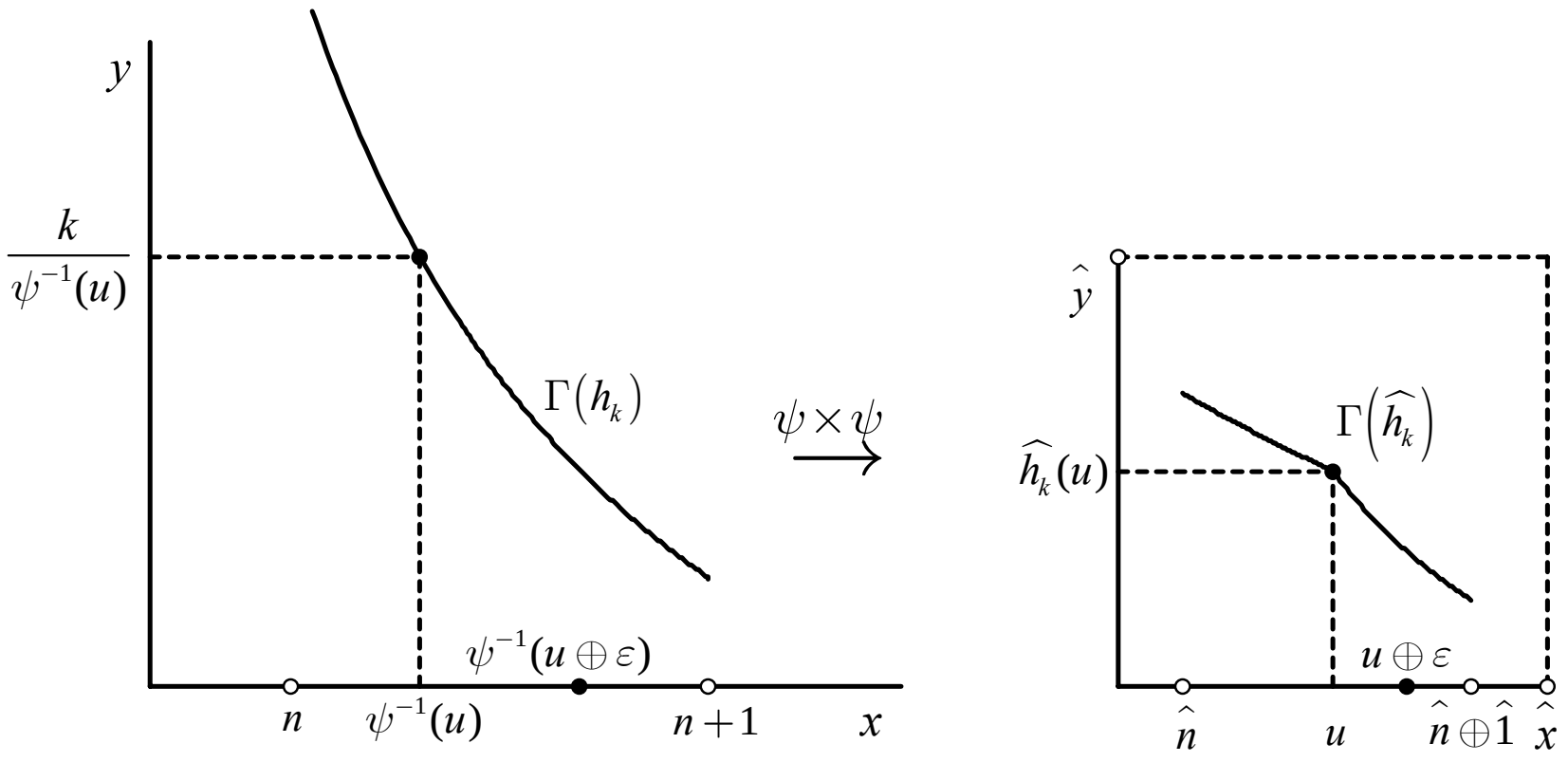

(Fig. 7)

$$[u,u\oplus\varepsilon)\xrightarrow{\psi^{-1}}\left[\psi^{-1}(u),\psi^{-1}(u\oplus\varepsilon)\right)\xrightarrow{h_k}\left(\frac{k}{\psi^{-1}(u\oplus\varepsilon)},\frac{k}{\psi^{-1}(u)}\right]\xrightarrow{\psi}\left(\hat{k}\div\widehat{\psi^{-1}(u\oplus\varepsilon)},\hat{k}\div\widehat{\psi^{-1}(u)}\right]$$

Therefore $\forall t\in[u,u\oplus\varepsilon)$ we verify $\widehat{h_k}(t)=\left(\psi_{\frac{k}{\psi^{-1}(u)}-1}\right)\left(\frac{k}{\psi^{-1}(t)}\right)$. Would result

$$\left(\widehat{h_k}\right)'_+(u)=\left(\psi_{\frac{k}{\psi^{-1}(u)}-1}\right)'_-\left(\frac{k}{\psi^{-1}(u)}\right)\cdot\frac{-k}{\left(\psi^{-1}(u)\right)^2}\cdot\frac{1}{(\psi_n)'\left(\psi^{-1}(u)\right)}$$

*<u>Case 2</u>* $u=\hat{n}$ $\left(n\in\mathbb{N}^*\right)$

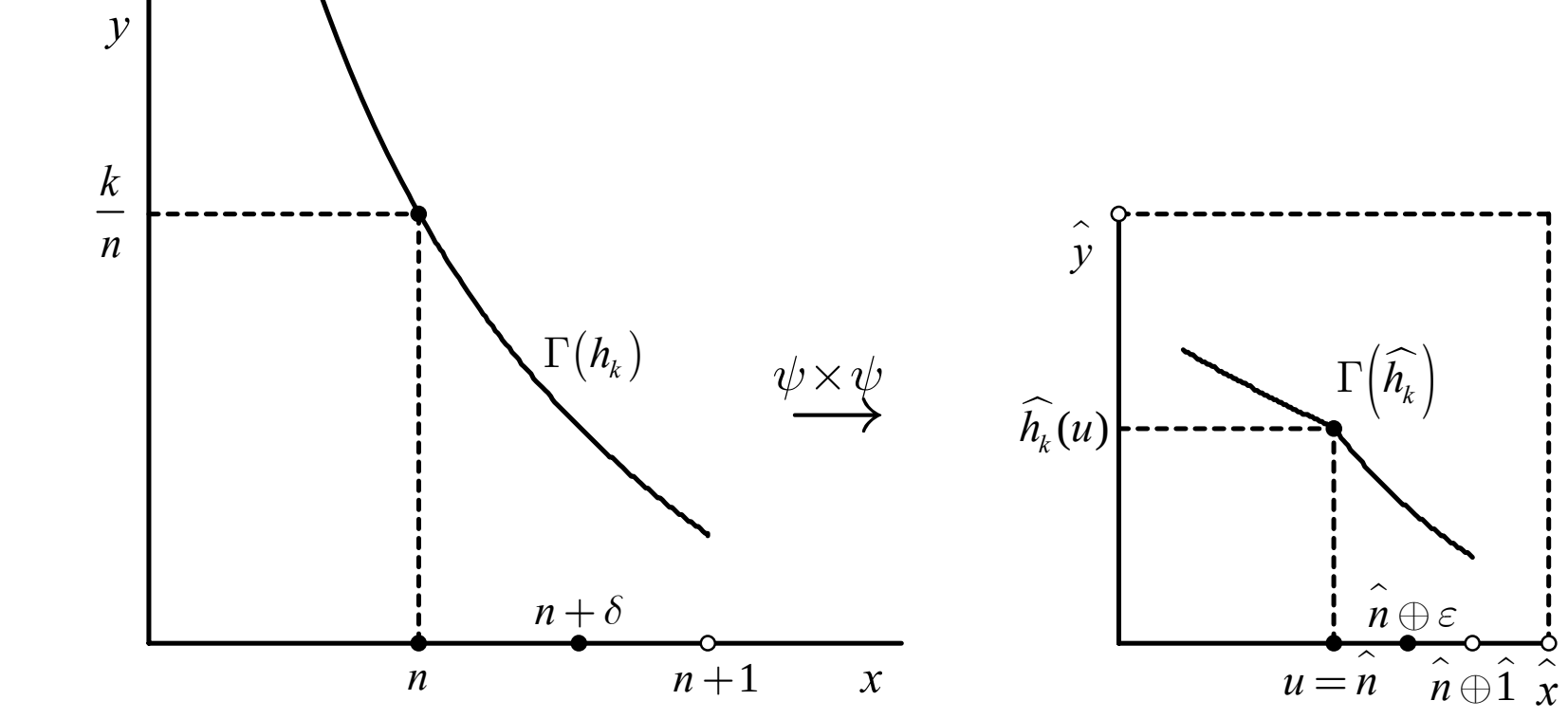

(Fig. 8)

For a sufficiently small $\varepsilon>0$, we obtain:

$$[\hat{n},\hat{n}\oplus\varepsilon)\xrightarrow{\psi^{-1}}[n,n+\delta)\xrightarrow{h_k}\left(\frac{k}{n+\delta},\frac{k}{n}\right]\xrightarrow{\psi}\left(\hat{k}\div(\hat{n}\oplus\hat{\delta}),\hat{k}\div\hat{n}\right]$$

Then , $\forall t \in [\hat{n}, \hat{n} \oplus \varepsilon)$ we verify : $\widehat{h_k}(t) = \psi_{\left[\frac{k}{n}\right]}\left(\frac{k}{\psi^{-1}(t)}\right)$(if $\frac{k}{n} \notin \mathbb{N}^*$) or

$$\widehat{h_k}(t) = \psi_{\frac{k}{n}-1}\left(\frac{k}{\psi^{-1}(t)}\right) \quad (\text{ if } \frac{k}{n} \in \mathbb{N}^*)$$

Then, $\left(\widehat{h_k}\right)'_+(\hat{n}) = \left(\psi_{\left[\frac{k}{n}\right]}\right)'\left(\frac{k}{n}\right) \cdot \frac{-k}{n^2} \cdot \frac{1}{(\psi_n)'_+(n)}$ ( if $\frac{k}{n} \notin \mathbb{N}^*$) or

$$\left(\widehat{h_k}\right)'_+(\hat{n}) = \left(\psi_{\frac{k}{n}-1}\right)'_-\left(\frac{k}{n}\right) \cdot \frac{-k}{n^2} \cdot \frac{1}{(\psi_n)'_+(n)} \quad (\text{ if } \frac{k}{n} \in \mathbb{N}^*).$$

Finally we have to study the differentiability of $\widehat{h_k}$ at $u = \hat{n}$ from the left side.

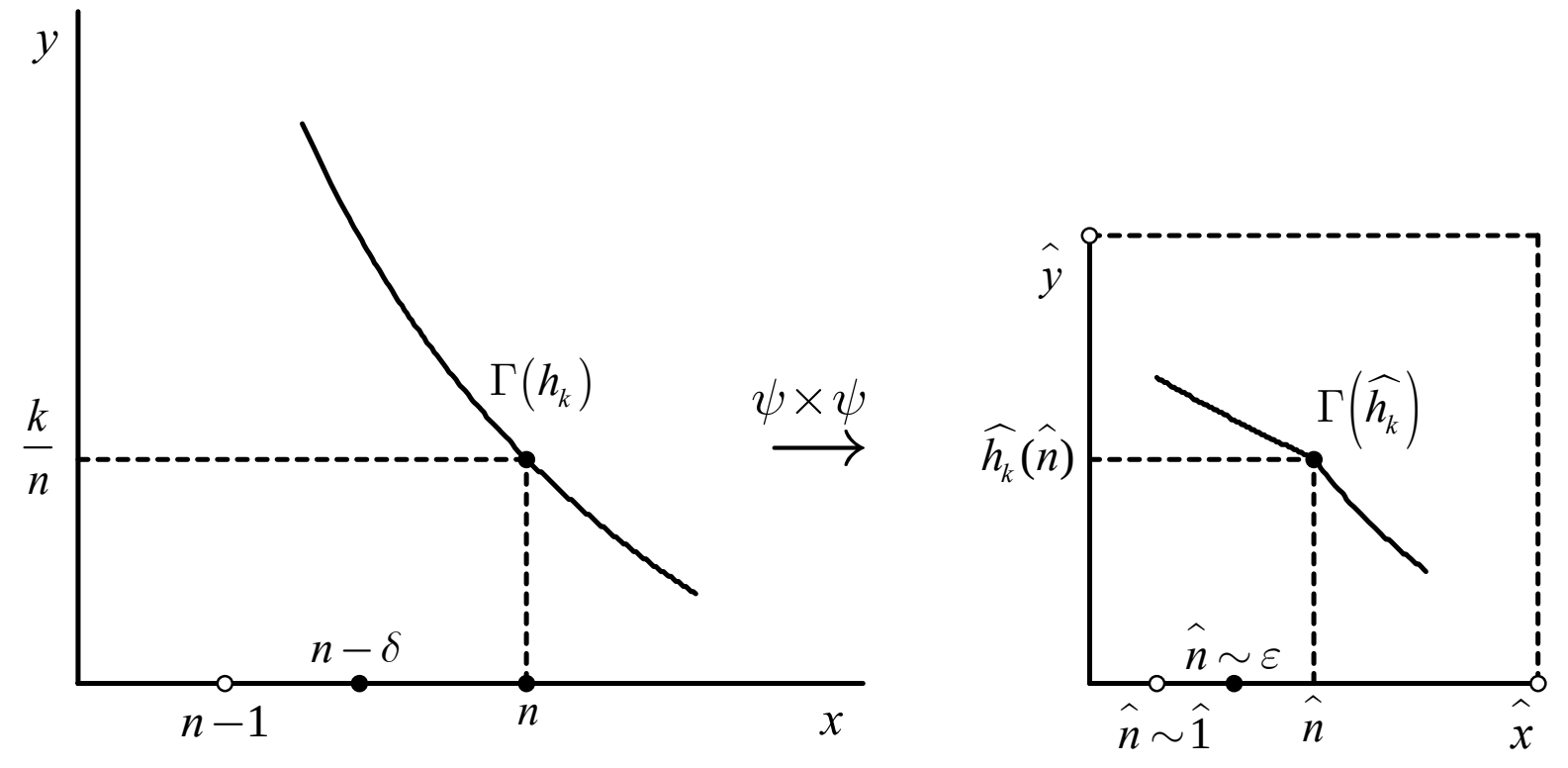

(Fig. 9)

For a sufficiently small $\varepsilon > 0$ we have that

$$(\hat{n} \sim \varepsilon, \hat{n}] \xrightarrow{\psi^{-1}} (n-\delta, n] \xrightarrow{h_k} \left[\frac{k}{n}, \frac{k}{n-\delta}\right) \xrightarrow{\psi} \left[\hat{k} \div \hat{n}, \hat{k} \div (\hat{n} \sim \hat{\delta})\right)$$

then $\forall t \in (\hat{n} - \varepsilon, \hat{n}]$ we verify $\widehat{h_k}(t) = \psi_{\left[\frac{k}{n}\right]}\left(\frac{k}{\psi^{-1}(t)}\right)$ regardless of whether $\frac{k}{n}$ is a natural number or not. This therefore would result

$$\left(\widehat{h_k}\right)'_-(\hat{n}) = \left(\psi_{\left[\frac{k}{n}\right]}\right)'_+\left(\frac{k}{n}\right) \cdot \frac{-k}{n^2} \cdot \frac{1}{(\psi_{n-1})'_-(n)}.$$

We have completed our examination of the differentiability of $\widehat{h_k}$ when dependent and/or independent variables take natural or non-natural number values. We can conclude all of the above in the following outline, bearing in mind that:

$$(\psi_{m-1})'_-(m) = a_m, \; (\psi_m)'_+(m) = b_m$$

**1.3.2 Conclusion** Where $\hat{x}, \hat{y} \notin \mathbb{N}$ , $\hat{n}, \hat{n}' \in \mathbb{N}^*$, $[x] = n, [y] = n'$ we obtain

1.-$(\hat{x}, \hat{y}) \in \Gamma\left(\widehat{h_k}\right)$, then $\left(\widehat{h_k}\right)'(\hat{x}) = \frac{-k}{x^2} \cdot \frac{\left(\psi_{[y]}\right)'(y)}{(\psi_n)'(x)}$

2.- $(\hat{x},\hat{n}') \in \Gamma\left(\widehat{h_k}\right)$, then $\left(\widehat{h_k}\right)'_+(\hat{x}) = \dfrac{-k}{x^2}\cdot\dfrac{a_{n'}}{(\psi_n)'(x)}$, $\left(\widehat{h_k}\right)'_-(\hat{x}) = \dfrac{-k}{x^2}\cdot\dfrac{b_{n'}}{(\psi_n)'(x)}$

3.- $(\hat{n},\hat{y}) \in \Gamma\left(\widehat{h_k}\right)$, then $\left(\widehat{h_k}\right)'_+(\hat{n}) = \dfrac{-k}{n^2}\cdot\dfrac{(\psi_{[y]})'(y)}{b_n}$, $\left(\widehat{h_k}\right)'_-(\hat{n}) = \dfrac{-k}{n^2}\cdot\dfrac{(\psi_{[y]})'(y)}{a_n}$

4.- $(\hat{n},\hat{n}') \in \Gamma\left(\widehat{h_k}\right)$, then $\left(\widehat{h_k}\right)'_+(\hat{n}) = \dfrac{-k}{n^2}\cdot\dfrac{a_{n'}}{b_n}$, $\left(\widehat{h_k}\right)'_-(\hat{n}) = \dfrac{-k}{n^2}\cdot\dfrac{b_{n'}}{a_n}$

Using the above diagram/summary, we obtain the following theorem:

**1.3.3 Theorem** Let $\psi$ be an $\mathbb{R}^+$ coding function, assume $[x]=n,[y]=n'$ and

$\widehat{h_k}:(0,M_\psi)\to(0,M_\psi)$, $\widehat{h_k}(u)=\psi\left(\dfrac{k}{\psi^{-1}(u)}\right)$. Then:

i) $(\hat{x},\hat{y}) \in \Gamma\left(\widehat{h_k}\right)$ where $x\notin\mathbb{N}\wedge y\notin\mathbb{N}$ then : $\widehat{h_k}$ is differentiable at $\hat{x}$

ii) $(\hat{x},\hat{n}') \in \Gamma\left(\widehat{h_k}\right)$ where $x\notin\mathbb{N}\ \wedge n'\in\mathbb{N}^*$ then: "$\widehat{h_k}$ is differentiable at $\hat{x} \Leftrightarrow a_{n'}=b_{n'}$"

iii) $(\hat{n},\hat{y}) \in \Gamma\left(\widehat{h_k}\right)$ where $n\in\mathbb{N}^*\wedge y\notin\mathbb{N}$ then: "$\widehat{h_k}$ is differentiable at $\hat{n} \Leftrightarrow a_n=b_n$"

iv) $(\hat{n},\hat{n}') \in \Gamma\left(\widehat{h_k}\right)$ where $n\in\mathbb{N}^*\wedge n'\in\mathbb{N}^*$ then:

"$\widehat{h_k}$ is differentiable at $\hat{n} \Leftrightarrow a_n\cdot a_{n'}=b_n\cdot b_{n'}$"

If we want the $\widehat{h_k}$ functions to be only differentiable at the points where both the ordinate and the abscissa are not natural numbers, we must select $\psi$ in such a way that the following is true

$$(*)\begin{cases} a_n\neq b_n \\ a_{n'}\neq b_{n'} \\ a_n\cdot a_{n'}\neq b_n\cdot b_{n'}\ (n\neq n') \end{cases} \qquad \text{or equivalently, } a_n\cdot a_{n'}\neq b_n\cdot b_{n'}\left(\forall n\in\mathbb{N}^*,\forall n'\in\mathbb{N}^*\right).$$

Selecting $\psi$ in this way we say that $\psi$ is an $\mathbb{R}^+$ coding function which identifies primes.

**1.3.4 Definition** We say that an $\mathbb{R}^+$ coding function identifies primes iff the $\widehat{h_k}$ functions are only differentiable at the non-natural number abscissa and ordinate points.

## 1.4 Classification of points in the $\hat{x}\hat{y}$ plane

Let $\psi:\mathbb{R}^+\to[0,M_\psi)$ be an $\mathbb{R}^+$ coding function that identifies primes, we have the following cases:

1.- $(\hat{x},\hat{y})\in[0,M_\psi)^2$ ($x\notin\mathbb{N}\wedge y\notin\mathbb{N}$). Then in a neighbourhood $V$ of $(\hat{x},\hat{y})$ we verify

*"$\forall(\hat{s},\hat{t})\in V$, the hyperbola which contains $(\hat{s},\hat{t})$ is differentiable at $\hat{s}$"*

(Of course, we mean to say the function which represents the graph of the transformed hyperbola $\Gamma\left(\widehat{h_{st}}\right)$).

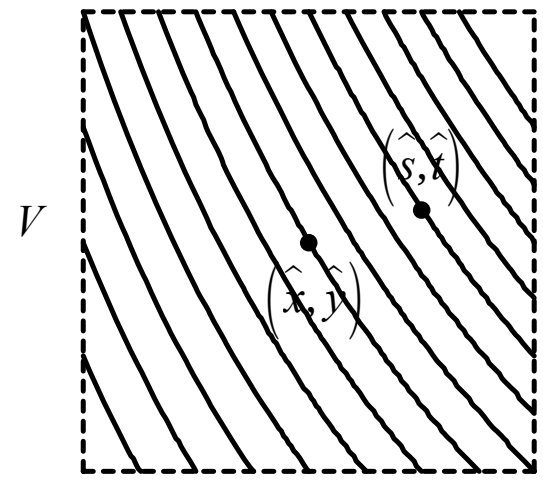

(Fig. 10)

2.- $(\hat{x},\hat{n}')\in\left[0,M_\psi\right)^2$ ( $x\notin\mathbb{N}\wedge n'\in\mathbb{N}^*$ ). Then, in a neighbourhood $V$ of $(\hat{x},\hat{n}')$ we verify

"$\forall(\hat{s},\hat{t})\in V$ *, the hyperbola which contains* $(\hat{s},\hat{t})$ *is differentiable at* $\hat{s}$ $\Leftrightarrow$ $\hat{t}\neq\hat{n}'$ "

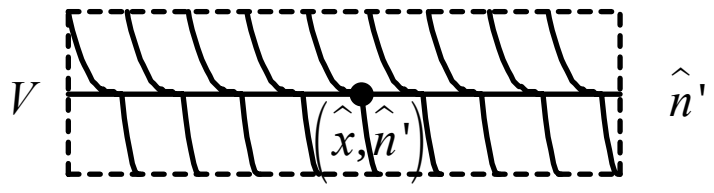

(Fig. 11)

3.- $(\hat{n},\hat{y})\in\left[0,M_\psi\right)^2$ ( $n\in\mathbb{N}^*\wedge y\notin\mathbb{N}$ ). Then in a neighbourhood $V$ of $(\hat{n},\hat{y})$ we verify

"$\forall(\hat{s},\hat{t})\in V$ *, the hyperbola which contains* $(\hat{s},\hat{t})$ *is differentiable at* $\hat{s}$ $\Leftrightarrow$ $\hat{s}\neq\hat{n}$ "

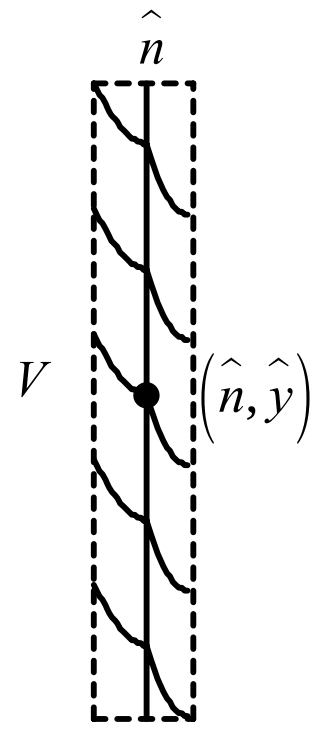

(Fig. 12)

4.- $(\hat{n},\hat{n}')\in\left[0,M_\psi\right)^2$ ( $n\in\mathbb{N}^*\wedge n'\in\mathbb{N}^*$ ). Then in a neighbourhood $V$ of $(\hat{n},\hat{n}')$ we verify

"$\forall(\hat{s},\hat{t})\in V$ *, the hyperbola which contains* $(\hat{s},\hat{t})$ *is differentiable at* $\hat{s}$ $\Leftrightarrow$ $\hat{s}\neq\hat{n}$ *and* $\hat{t}\neq\hat{n}'$"

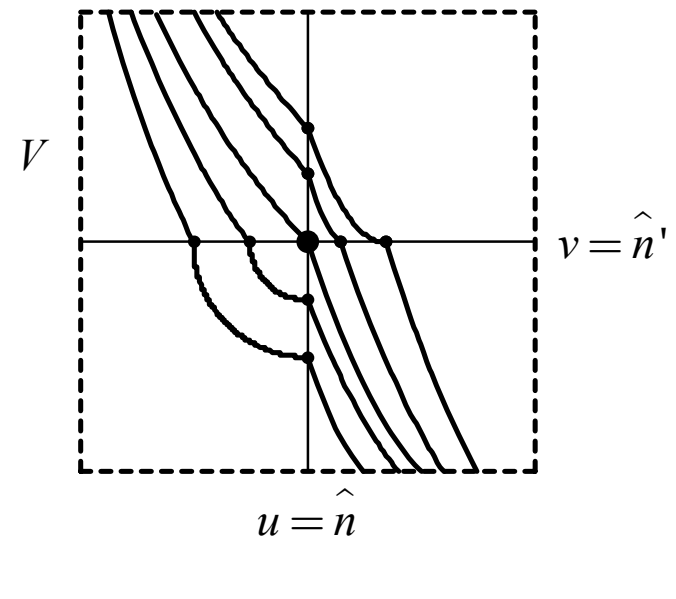

(Fig. 13)

Given the symmetry with respect to the bisector line of the first quadrant, we will only work where $\hat{y} \geq \hat{x}$. We will additionally use $\hat{x} \geq \hat{1}$. Clearly $k \in \mathbb{N}^*$ is not prime, iff $\hat{x} > \hat{1}$ and the "hyperbola" $\hat{x} \otimes \hat{y} = \hat{k}$ contains a point of type 4 or a point with the following shape:

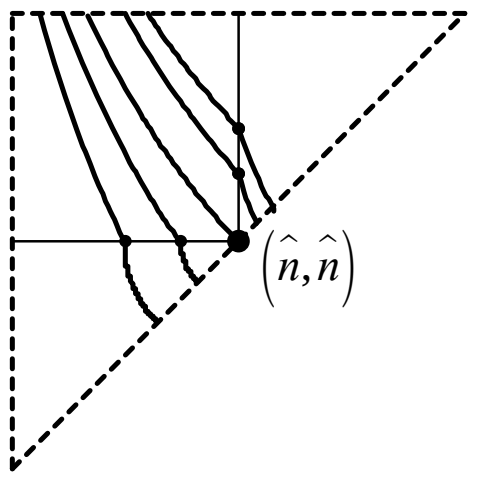

(Fig. 14)
($k$ could be a perfect square)

Graphically, we obtain:

(i) (ii) (iii)

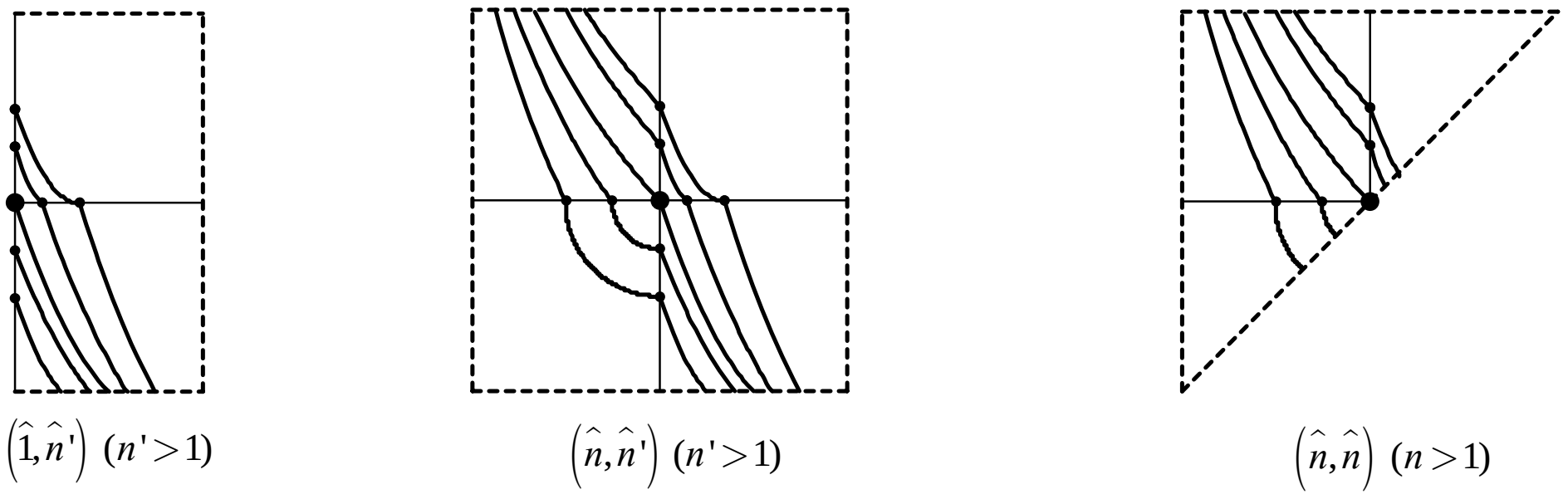

(Fig. 15)

**1.4.1 Definition** Let $\psi$ be an $\mathbb{R}^+$ coding function that identifies primes and assume that $(\hat{x},\hat{y}) \in \left[\hat{1}, M_\psi\right)^2$ and $\hat{y} \geq \hat{x}$. If $(\hat{x},\hat{y})$ is of type (i) we say that it is a semi-vortex point with respect to $\psi$. If $(\hat{x},\hat{y})$ is of type (ii) or (iii) we say that it is a vortex point with respect to $\psi$. We should insist that we are considering the graphs of $\hat{x} \otimes \hat{y} = \hat{k}$ where $\hat{x} \geq \hat{1}$ and. $\hat{y} \geq \hat{x}$.

## 1.5 Hyperbolic classification of natural numbers

Let $\hat{k} \in \left(\hat{1}, M_\psi\right)$ .According to the statements made above, we may classify $\hat{k}$ in terms of the behaviour of "hyperbolas" that are near the $\hat{x} \otimes \hat{y} = \hat{k}$ "hyperbola". We obtain the following classification:

1) $\hat{k}$ is natural iff the $\hat{x} \otimes \hat{y} = \hat{k}$ "hyperbola" contains a semi-vortex point.

2) $\hat{k}$ is prime iff the $\hat{x}\otimes\hat{y}=\hat{k}$ "hyperbola" contains a semi-vortex point and no vortex points.
3) $\hat{k}$ is a non-prime natural number iff the $\hat{x}\otimes\hat{y}=\hat{k}$ "hyperbola" contains a semi-vortex point and at least one vortex point.
4) $\hat{k}$ is not a natural number iff the $\hat{x}\otimes\hat{y}=\hat{k}$ "hyperbola" does not contain semi-vortex (which would imply that it has not vortex points either).

**1.5.1 Remark** Note that there is no difference between a non-natural number and a prime number, except that a semi-vortex point appears in the latter case.

**1.5.2 Question** In which sense is the difference between a non-natural number and a prime relevant according to the hyperbolic classification of natural numbers? .Or, according to hyperbolic classification, given a prime number $\hat{k}$, is it more like a non-natural number or a non-prime natural number? .This question will be answered in the third chapter of this paper and we will see that it has great conceptual depth.

## 1.6 $\mathbb{R}^+$ prime coding

Amongst the $\mathbb{R}^+$coding functions that identifies primes, it will be interesting to select those given by $\psi_m:[m,m+1]\to\mathbb{R}^+\ (m=0,1,2,...)$, that is, functions that are affine

$$\psi_m(x)=\xi_m(x-m)+B_m\ \ (B_0=0,\ B_m=\sum_{j=0}^{m-1}\xi_j\ \text{ if } m\geq 1)$$

$$(\xi_m>0\ \ \forall m\in\mathbb{N})$$

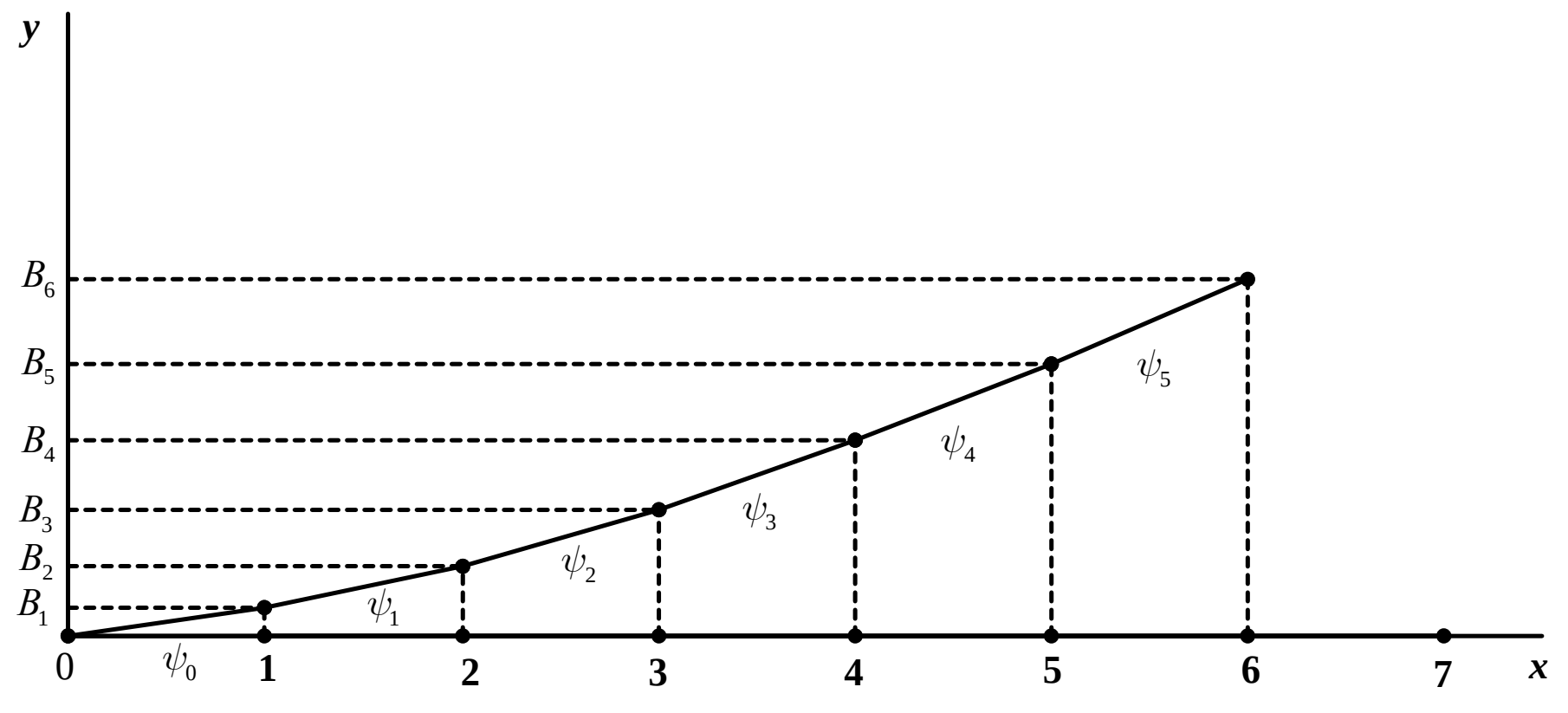

(Fig. 16)

We can easily prove that the $\psi$ functions defined by means of the sequence $(\psi_m)_{m\geq 0}$ are $\mathbb{R}^+$coding functions. The (*) conditions in section 1.3.3 for $\psi$ to identify primes can now thus be expressed:

$$\psi \text{ identifies primes} \Leftrightarrow \begin{cases}\xi_i\neq\xi_{i+1}\\ \xi_i\cdot\xi_j\neq\xi_{i+1}\cdot\xi_{j+1}\end{cases}\quad (**)$$

Or equivalently, $\psi$ identifies primes $\Leftrightarrow \xi_i\cdot\xi_j\neq\xi_{i+1}\cdot\xi_{j+1}\ \ (\forall i\forall j\in\mathbb{N})$

The fulfilment of the inequality above is guaranteed by choosing $\xi_i$ such that $0 < \xi_0 < \xi_1 < \xi_2 < \ldots$ though this is not the only way of choosing it.

**1.6.1 Definition** Any $\mathbb{R}^+$ coding function $\psi$ that is defined by means of $\psi_m$ affine functions that also satisfies $0 < \xi_0 < \xi_1 < \xi_2 < \ldots$ it is said to be $\mathbb{R}^+$ *prime coding.* We call the numbers $\xi_0, \xi_1, \xi_2, \xi_3, \ldots$ coefficients of the $\mathbb{R}^+$ *prime coding.* .

# 2. ESSENTIAL REGIONS

## 2.1 Essential regions associated with a hyperbola

**2.1.0 Introduction** In the $\hat{x}\hat{y}$ plane and for any given even number $\hat{\alpha} \geq \widehat{16}$ we will consider the function in which any number $\hat{k}$ of the closed interval $\left[\hat{4}, \hat{\alpha} \div \hat{2}\right]$ corresponds to the area of the region of $\hat{x}\hat{y}$ limited by the curves $\hat{x} = \hat{2}, \hat{y} = \hat{x}$, $\hat{x} \otimes \hat{y} = \hat{k}$ (called lower area) and also the function that associates each $\hat{k}$ to the area of the region of $\hat{x}\hat{y}$ limited by the curves $\hat{x} = \hat{2}, \hat{y} = \hat{x}$ , $\hat{x} \otimes \hat{y} = \hat{\alpha} \sim \hat{k}$ (called upper area). The $\hat{x}\hat{y}$ plane is considered imbedded in the $xy$ plane with the usual measure. This means that for any given even number $\hat{\alpha} \geq \widehat{16}$ we have $\hat{\alpha} = \hat{k} \oplus (\hat{\alpha} \sim \hat{k})$ and, associated to this decomposition, two data pieces, lower and upper areas. We will study if $\hat{\alpha}$ is the sum of the two prime numbers $\hat{k}_0$ and $\hat{\alpha} \sim \widehat{k_0}$ taking into account the restrictions $\hat{\alpha} \sim \hat{3}$ and $\hat{\alpha} \div \hat{2}$ non-prime. The upper and lower area functions will not yet yield any characterizations to the Goldbach Conjecture. We will need the second derivative of the total area function ( the sum of the lower and upper areas ).

To this end, we define the concept of essential regions associated to a hyperbola which, simply put, is any region in the $xy$ plane with the shape $[n, n+1] \times [n', n'+1]$ where $n$ and $n'$ are natural numbers and $n' > n > 1$ and the hyperbola intersects it in more than one point or else the shape $[n, n+1]^2$ where $n > 1$ and $x \leq y$ and the hyperbola intersects in more than one point.

These essential regions are then transported to the $\hat{x}\hat{y}$ plane by means of the $\psi \times \psi$ function, and we will find the total area function adding the areas determined by each hyperbola in the respective essential regions, and the second derivative of this area function in each essential region. After this process we obtain the formula which determines the second derivative function of the total area $\widehat{A_T}$ in each sub-interval $\left[\widehat{k_0}, \widehat{k_0} \oplus \hat{1}\right]$, a derivative which is continuous

$$\left(\widehat{A_T}\right)''\left(\hat{k}\right) = \frac{x_{k_0}}{\xi_{k_0}^2} \cdot \frac{1}{k} + \frac{y_{k_0}}{\xi_{\alpha-k_0-1}^2} \cdot \frac{1}{\alpha - k}$$

$$(\hat{k} \in \left[\widehat{k_0}, \widehat{k_0 \oplus 1}\right], k_0 = 4,5,..., \frac{\alpha}{2} - 1)$$

Both $x_{k_0}$ and $y_{k_0}$ are numeric values in homogeneous polynomial of degree two obtained from substituting in their variables the $\xi_i$ coefficients of the $\psi$*-affine*. We call $P_{k_0} = \left(x_{k_0}, y_{k_0}\right)$ an essential point. The study of the behaviour of the second derivative in these intervals allows the following characterization of the Goldbach Conjecture for any even number $\alpha \geq 16$ with the restrictions $\alpha - 3$ and $\frac{\alpha}{2}$ non-prime.

*Given* $\alpha \geq 16$*, an even number, then,* $\alpha$ *is the sum of two prime numbers* $k_0$ *and* $\alpha - k_0$ *(*$5 \leq k_0 < \frac{\alpha}{2}$*) iff the consecutive essential points* $P_{k_0-1}$ *y* $P_{k_0}$ *are repeated, that is, if* $P_{k_0-1} = P_{k_0}$

**2.1.1 Definition** Consider the family of functions

$$\mathcal{H} = \left\{ h_k : \left[2, \sqrt{k}\,\right] \to \mathbb{R},\ h_k(x) = \frac{k}{x},\ k \geq 4 \right\}$$

whose graphs represent the pieces of the hyperbolas $xy = k (k \geq 4)$ included in the subset of $\mathbb{R}^2$, $S \equiv (x \geq 2) \wedge (x \leq y)$. Consider the subsets of $\mathbb{R}^2$:

a) $R_{(n,n')} = \left[n, n+1\right] \times \left[n', n'+1\right] (2 \leq n < n', n \text{ and } n' \text{ natural numbers})$

b) $R_{(n,n)} = \left(\left[n, n+1\right] \times \left[n, n+1\right]\right) \cap \left\{\left(x, y\right) \in \mathbb{R}^2 : y \geq x\right\} (2 \leq n,\ n \text{ natural number})$

Let $h_k$ be an element of $\mathcal{H}$. We say that $R_{(n,n')}$ is a square essential region of $h_k$ iff $R_{(n,n')} \cap \Gamma\left(h_k\right)$ contains more than one point. We say that $R_{(n,n)}$ is a triangular essential region of $h_k$ iff $R_{(n,n)} \cap \Gamma\left(h_k\right)$ contains more than one point.

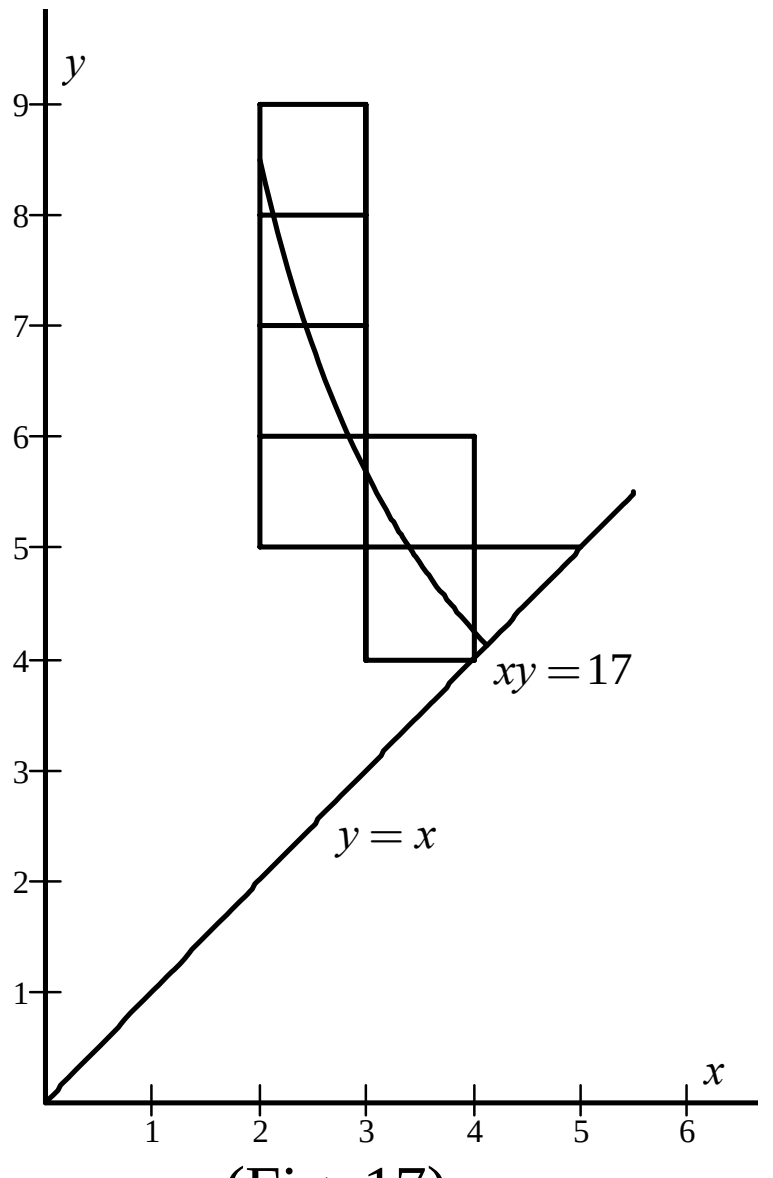

(Fig. 17)

The essential regions of the $xy = 17$ hyperbola are:

$$R_{(2,8)}, R_{(2,7)}, R_{(2,6)}, R_{(2,5)}, R_{(3,5)}, R_{(3,4)}, R_{(4,4)}$$

Analyse the different types of essential regions depending on the way the hyperbola $xy = k$ intersects with $R_{(n,n')}$ $(n' > n)$. If the hyperbola passes through point $P(n, n'+1)$, then the equation for the hyperbola is $xy = n(n'+1)$.

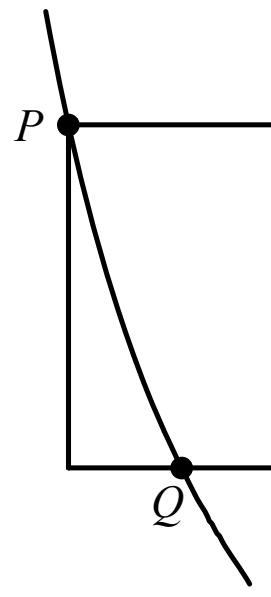

(Fig. 18)

The abscissa of the $Q$ point is $x = \frac{n(n'+1)}{n'}$. We verify that $n < \frac{n(n'+1)}{n'} < n+1$.

This is equivalent to say $\begin{cases} nn' < nn' + n \\ nn' + n < n'n + n' \end{cases}$ or equivalently $\begin{cases} 0 < n & (1) \\ n < n' & (2) \end{cases}$

(1) is trivial, and (2), by hypothesis. The remaining types are reasoned in an almost trivial way.

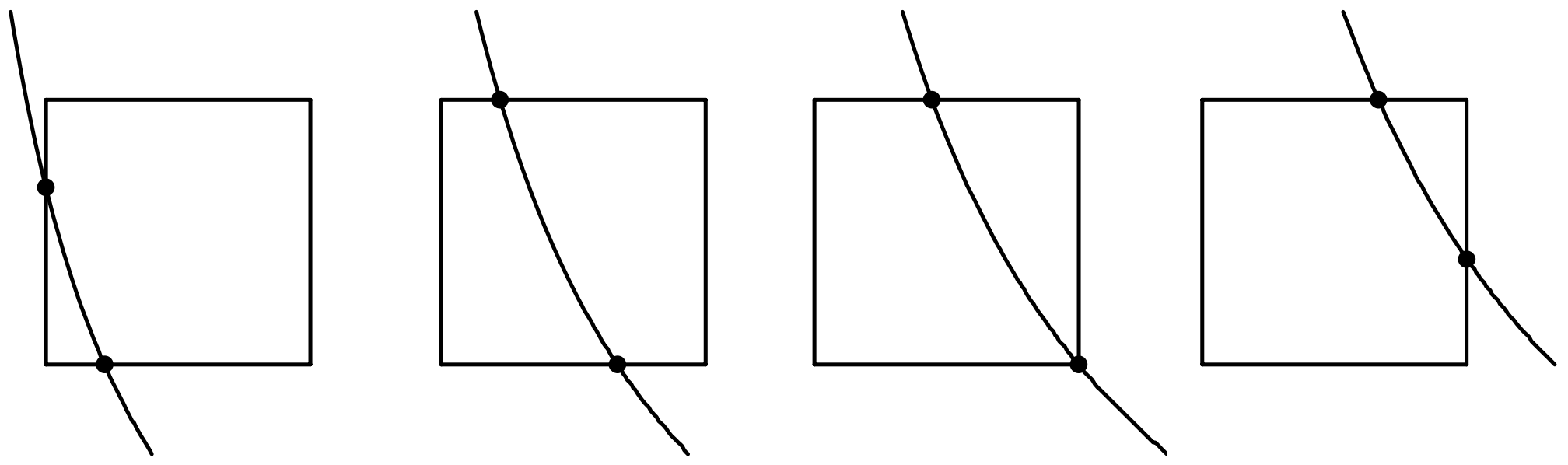

(Fig. 19)

We use the same considerations for the triangular essential regions $R_{(n,n)}$.

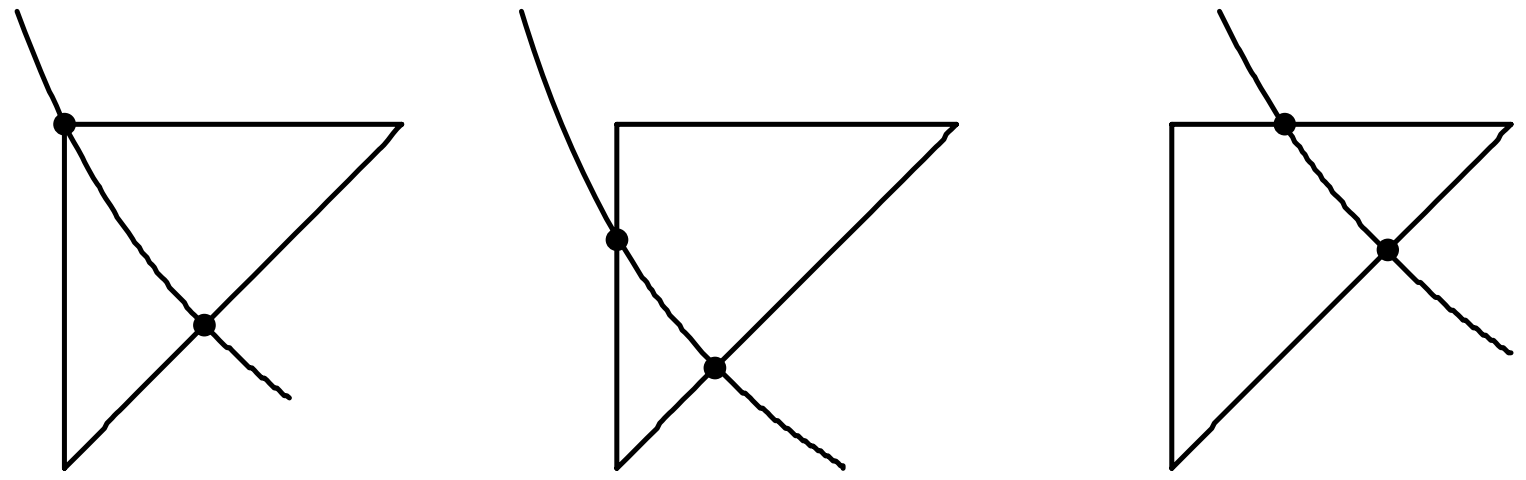

(Fig. 20)

Let $k_0 \in \mathbb{N}, k_0 \geq 4$. We will examine which are the types of essential regions for the hyperbolas $xy = k (y \geq x)$ where $k_0 < k < k_0 + 1$. The passage through essential regions of points $P_0, Q_0$ of the $xy = k_0$ hyperbola with relation to $P, Q$ points of the $xy = k$ hyperbola corresponds to the following diagrams:

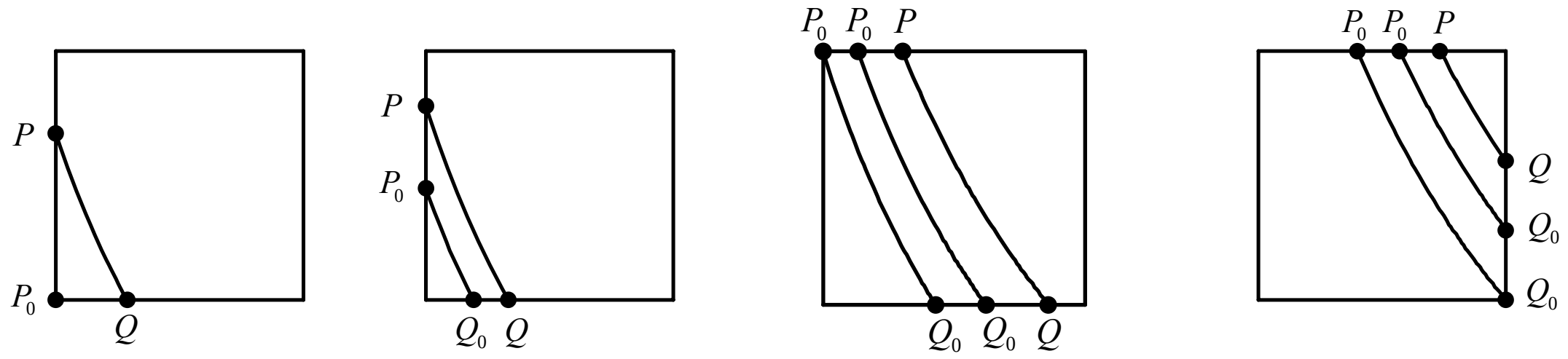

(Fig. 21)

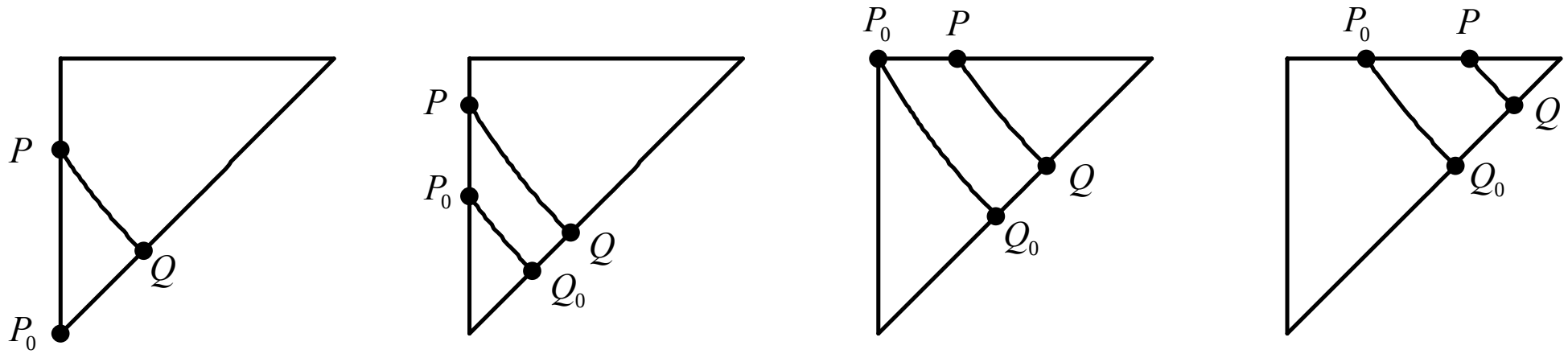

(Fig. 22)

Note that in any other case there would be a $xy = s$ $(k_0 < s < k)$ hyperbola included between $xy = k_0$ and $xy = k$ which would pass through a natural coordinate point, which is absurd.

**2.1.2 Consequence** The essential regions for the hyperbola $xy = k$ are of the following types:

a) Square essential regions $R_{(n,n')}$

Type 1 Type 2 Type 3 Type 4 Type 5

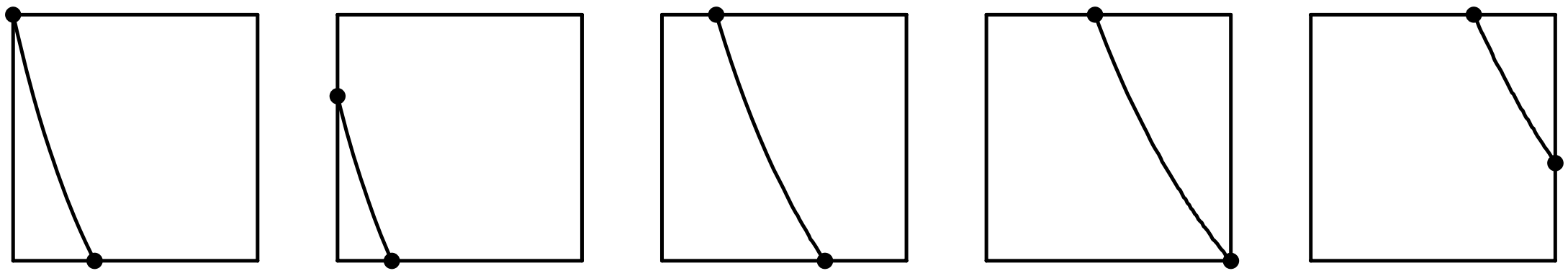

(Fig. 23)

b) Triangular essential regions $R_{(n,n)}$

Type 6 Type 7 Type 8

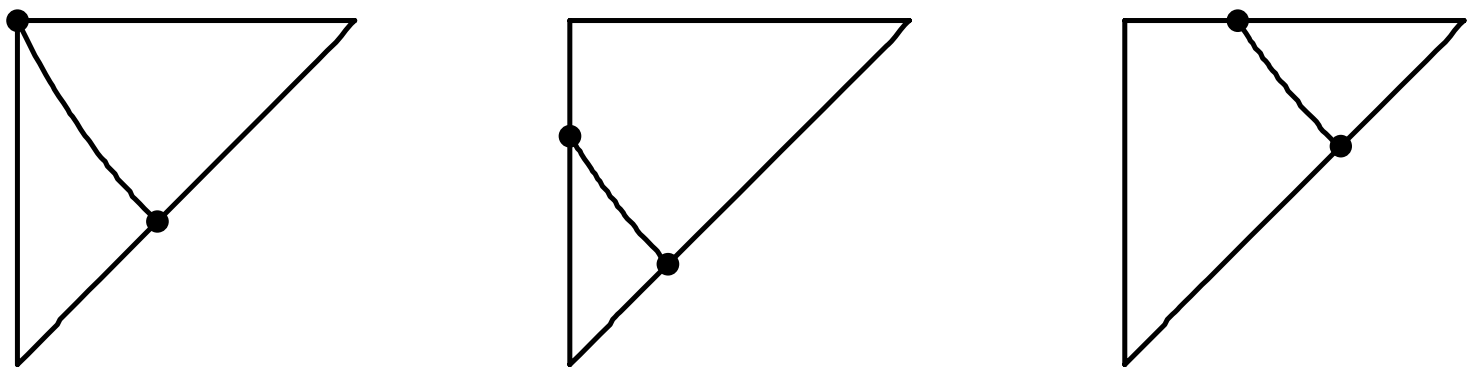

(Fig. 24)

We will find the essential regions of the $xy = k$ hyperbolas with the conditions:

$$k_0 \in \mathbb{N},\ k_0 \geq 4,\ k_0 < k < k_0 + 1\ .$$

The abscissa of $xy = k_0$ varies in the interval $\left[2, \sqrt{k_0}\right]$.

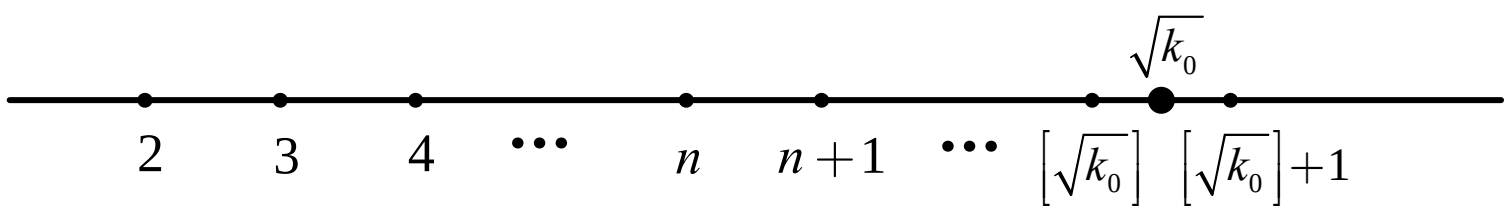

(Fig. 25)

a) For $n \in \{2,3,...,\left[\sqrt{k_0}\right]-1\}$ the $R_{(n,n')}$ essential regions of the hyperbola are obtained when $n'$ varies in the set

$$\{\left[\frac{k_0}{n+1}\right], \left[\frac{k_0}{n+1}\right]+1, ..., \left[\frac{k_0}{n}\right]\}$$

We can easily verify that if $n' = \left[\frac{k_0}{n}\right]$ then $R_{(n,n')}$ is a square essential region of Type 2, if $n' = \left[\frac{k_0}{n+1}\right]$, $R_{(n,n')}$ is a square essential region of Type 5 and the remaining $R_{(n,n')}$ are of Type 3 (Fig. 26).

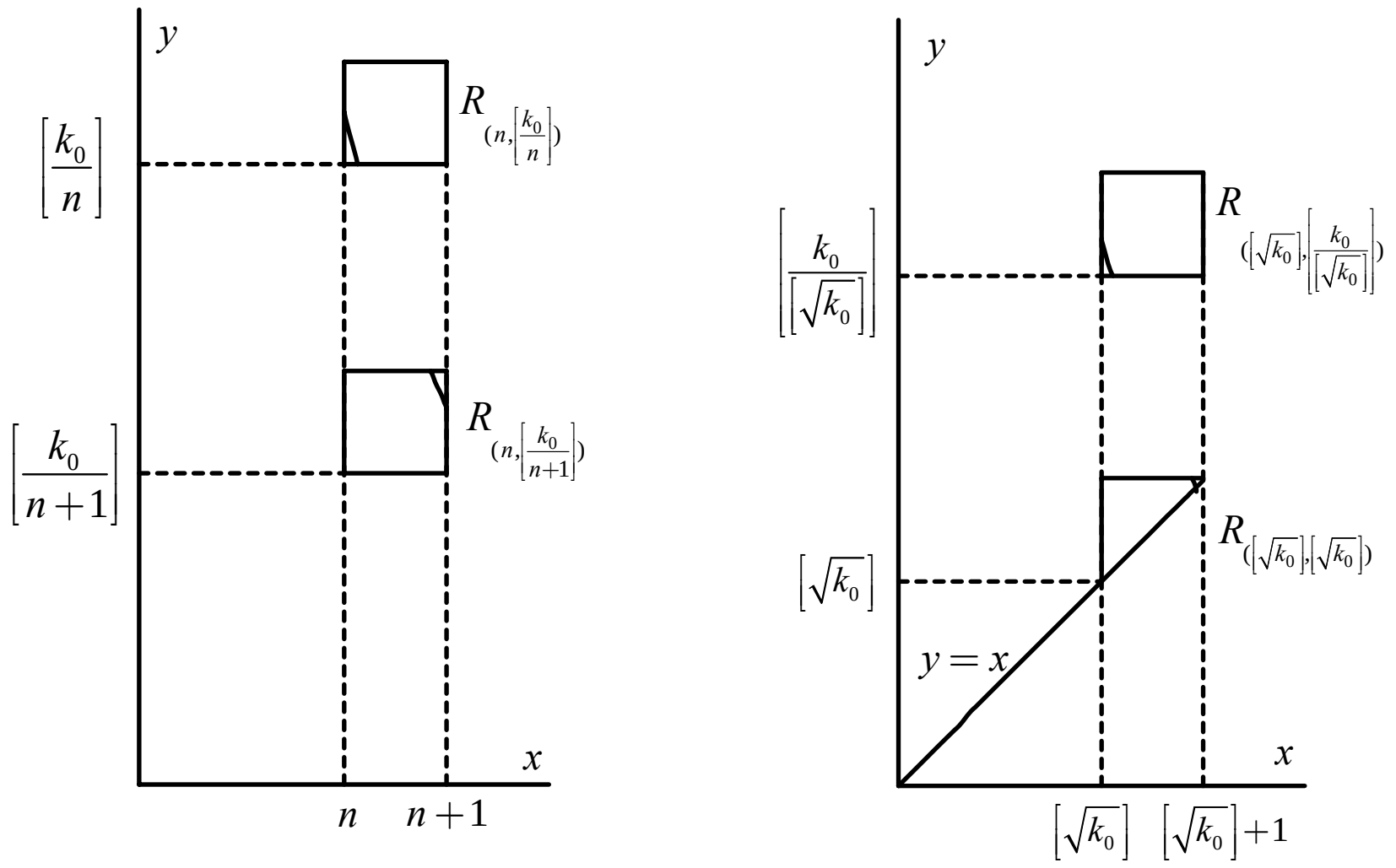

(a) (b)

(Fig. 26)

b) For $n = \left[\sqrt{k_0}\right]$, the $R_{\left(\sqrt{k_0},m\right)}$ essential regions are obtained when $m$ varies in the set

$$\{\left[\sqrt{k_0}\right], \left[\sqrt{k_0}\right]+1, ..., \left[\frac{k_0}{\left[\sqrt{k_0}\right]}\right]\}$$

If $m = \left[\sqrt{k_0}\right]$ we obtain a triangular essential region and could eventually exist a square essential region $R_{\left(\sqrt{k_0},m\right)}$ (Fig. 26 (b)).

Consider the set of indexes $\{(n, i(n))\}$ such that

$$(1) \quad n = 2,3,...,\left[\sqrt{k_0}\right]-1, \quad i(n) = \left[\frac{k_0}{n+1}\right], \left[\frac{k_0}{n+1}\right]+1, ..., \left[\frac{k_0}{n}\right]$$

(2) $n=\left[\sqrt{k_0}\right],\ i(n)=\left[\sqrt{k_0}\right],\left[\sqrt{k_0}\right]+1,\ldots,\left[\frac{k_0}{\left[\sqrt{k_0}\right]}\right]$

Let $E_s(k_0)$ be the set $\{n,i(n)\}$, where $(n,i(n))$ are pairs of type (1) or of type (2) . We obtain the following proposition:

**2.1.3 Proposition** Let $k_0 \in \mathbb{N}^*$ ($k_0 \geq 4$). Then**,**
i) All the $xy=k\ (k_0<k<k_0+1)$ hyperbolas have the same essential regions and the same type. ii) The $xy=k$ essential regions are elements of the set

$$\left\{R_{n,i(n)}:\left(n,i(n)\right)\in E_s\left(k_0\right)\right\}$$

**2.1.4 Example** For $k_0=18$ the essential regions of the $xy=k\ (18<k<19)$ hyperbolas are (Fig. 27):

<u>n=2</u> $\begin{cases}\text{Type 2}: R_{(2,9)}\\ \text{Type 3}: R_{(2,8)}, R_{(2,7)}\\ \text{Type 5}: R_{(2,6)}\end{cases}$, <u>n=3</u> $\begin{cases}\text{Type 2}: R_{(3,6)}\\ \text{Type 3}: R_{(3,5)}\\ \text{Type 5}: R_{(3,4)}\end{cases}$, <u>n=4</u> Type 7 : $R_{(4,4)}$

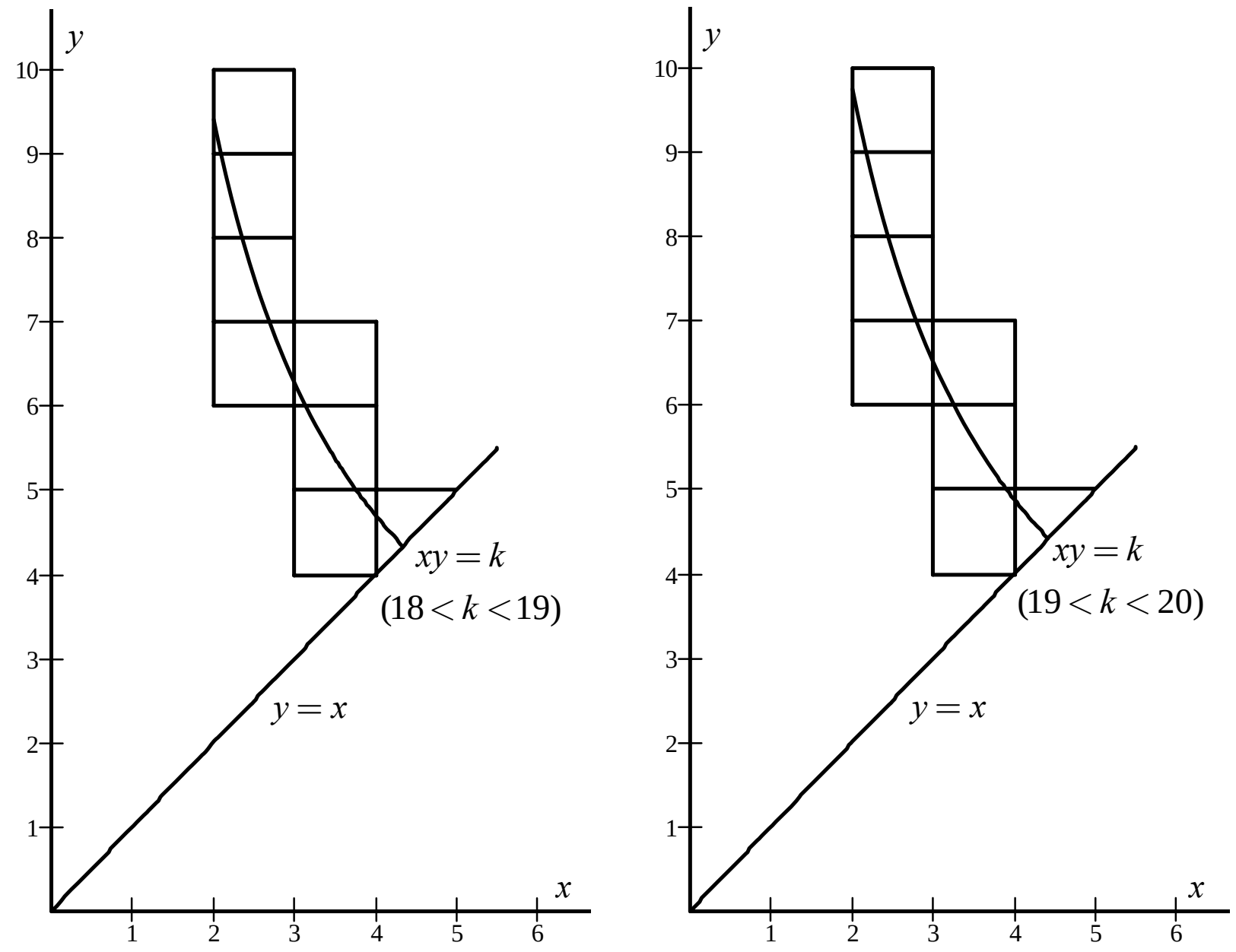

(Fig. 27)

The essential regions of the $xy=k$ (19<*k*<20) hyperbolas are exactly the same, due to the fact that 19 is a prime number.

## 2.2 Areas in essential regions associated with a hyperbola

To every $R_{(n,n')}(n\leq n')$ essential region of the $xy=k\ \left(k\notin\mathbb{N}^*, k>4\right)$ hyperbola, we will associate the region of the *xy* plane below the hyperbola (we call it $D_{(n,n')}(k)$ ). We call $A_{(n,n')}(k)$, the area of $D_{(n,n')}(k)$, considering the $\hat{x}\hat{y}$ plane embedded in the $xy$ plane. We have the following cases:

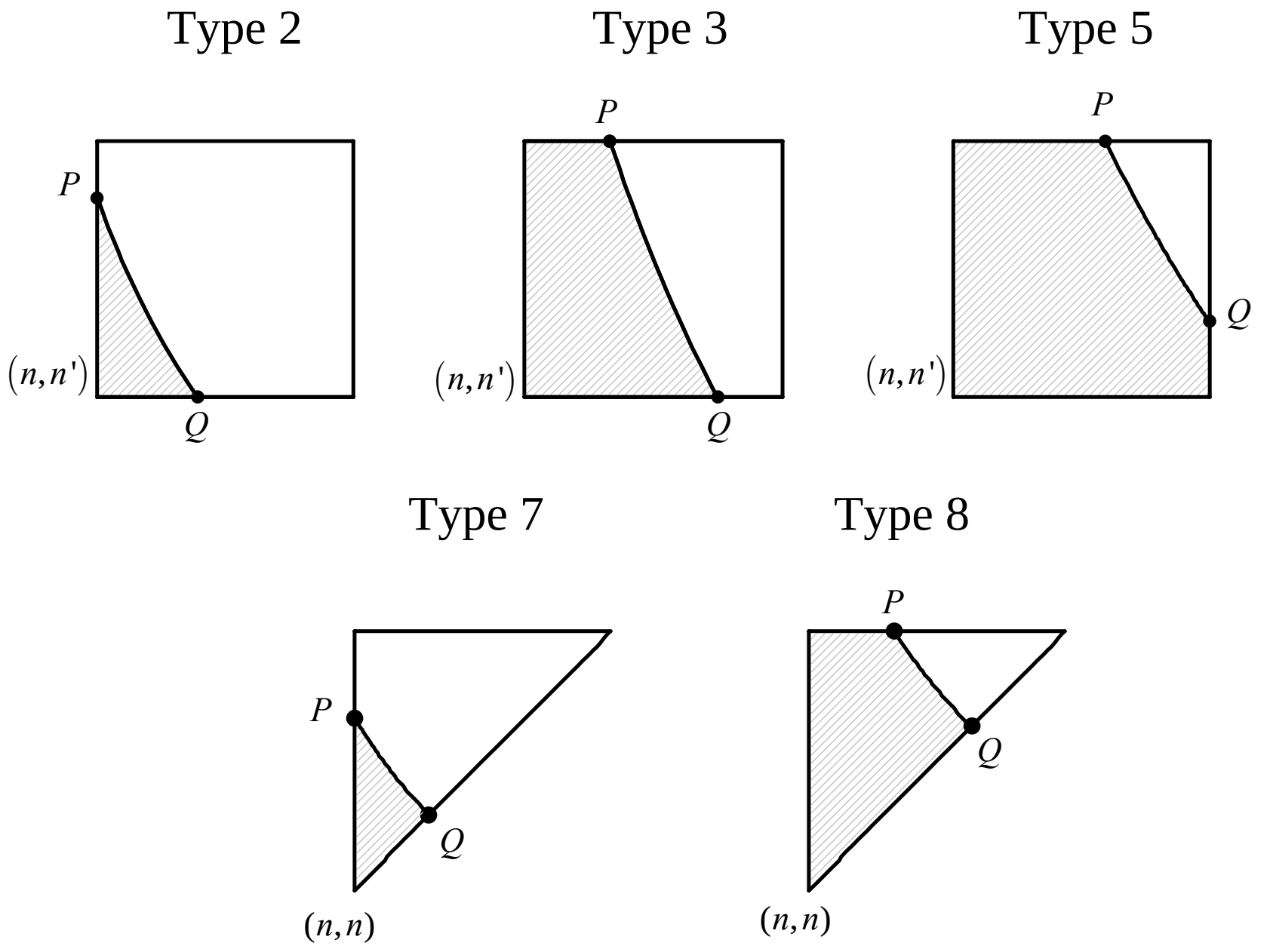

(Fig. 28)

1.- *Type* 2 *essential region.* $A_{(n,n')}(k)=\iint_{D_{(n,n')}} dxdy$, $\left(D_{(n,n')}\equiv n\le x\le \frac{k}{n'},\ n'\le y\le \frac{k}{x}\right)$

$$A_{(n,n')}(k)=\int_{n}^{k/n'}dx\int_{n'}^{k/x}dy=\int_{n}^{k/n'}\left(\frac{k}{x}-n'\right)dx=k\log\frac{k}{nn'}+nn'-k$$

This is the expression of the area. If $k\in[k_0,k_0+1]$ ( $k_0\ge 4$, natural number), then $A'_{(n,n')}(k)=\log\frac{k}{nn'}$ and the second derivative is $A''_{(n,n')}(k)=\frac{1}{k}$. Note that we have used the closed interval $[k_0,k_0+1]$ so we may extend the definition of the essential region for $k\in\mathbb{N}^*$ in a natural manner. In some cases, the "essential region" would consist of a single point (null area).

2.- *Type* 3 *essential region.* $A_{(n,n')}(k)=\iint_{D_{(n,n')}} dxdy$. In this case $D_{(n,n')}=D'\cup D''$, where

$$D'=\left[n,\frac{k}{n'+1}\right]\times[n',n'+1],\ D''\equiv\frac{k}{n'+1}<x\le\frac{k}{n'},\ n'\le y\le\frac{k}{x}\quad (D'\cap D''=\varnothing)$$

$$A_{(n,n')}(k)=\iint_{D_{(n,n')}}dxdy=\frac{k}{n'+1}-n+\iint_{D''}dxdy=\frac{k}{n'+1}-n+k\log\frac{n'+1}{n'}-n'\left(\frac{1}{n'}-\frac{1}{n'+1}\right)k$$

If $k_0\le k\le k_0+1$, then $A''_{(n,n')}(k)=0$

3.-*Type* 5 *essential region.* $A_{(n,n')}(k)=\iint_{D_{(n,n')}} dxdy$. In this case $D_{(n,n')}=D'\cup D''$ , where

$$D' = \left[n, \frac{k}{n'+1}\right]\times\left[n'.n'+1\right], D'' \equiv \frac{k}{n'+1} < x \le n+1,\ n' \le y \le \frac{k}{x}\ \ \left(D'\cap D'' = \varnothing\right)$$

$$A_{(n,n')}(k) = \iint_{D_{(n,n')}} dxdy = \frac{k}{n'+1} - n + \iint_{D''} dxdy = \frac{k}{n'+1} - n + k\log\frac{(n+1)(n'+1)}{k} -$$

$-n'\left(n+1-\frac{k}{n'+1}\right)$. In the interval $\left[k_0, k_0+1\right]$ we obtain $A'_{(n,n')}(k) = \log\frac{(n+1)(n'+1)}{k}$ and

$$A''_{(n,n')}(k) = -\frac{1}{k}$$

4.-*Type* 7 *essential region.* $A_{(n,n)}(k) = \iint_{D_{(n,n)}} dxdy\ \left(D_{(n,n)} \equiv n \le x \le \sqrt{k},\ x \le y \le \frac{k}{x}\right)$

$$A_{(n,n)}(k) = \int_n^{\sqrt{k}} dx \int_x^{k/x} dy = \int_n^{\sqrt{k}} \left(\frac{k}{x} - x\right) dx = \left[k\log x - \frac{x^2}{2}\right]_n^{\sqrt{k}} = \frac{k}{2}\log k - \frac{k}{2} - k\log n + \frac{n^2}{2}$$

If $k_0 \le k \le k_0+1$ we obtain $A'_{(n,n)}(k) = \frac{1}{2}\log k - \log n$ and $A''_{(n,n)}(k) = \frac{1}{2k}$.

5.-*Type* 8 *essential region.* $A_{(n,n)}(k) = \iint_{D_{(n,n)}} dxdy$. In this case $D_{(n,n)} = D' \cup D''$, where

$$D' \equiv n \le x < \frac{k}{n+1},\ x \le y \le n+1,\ D'' \equiv \frac{k}{n+1} < x \le \sqrt{k},\ x \le y \le \frac{k}{x}\ \ \left(D'\cap D'' = \varnothing\right)$$

$$A_{(n,n)}(k) = \iint_{D'} dxdy + \iint_{D''} dxdy = \int_n^{k/(n+1)} dx \int_x^{n+1} dy + \int_{k/(n+1)}^{\sqrt{k}} dx \int_x^{k/x} dy$$

$$A_{(n,n)}(k) = \int_n^{\frac{k}{n+1}} (n+1-x)\,dx + \int_{\frac{k}{n+1}}^{\sqrt{k}} \left(\frac{k}{x} - x\right) dx = \frac{k}{2} - n(n+1) + \frac{n^2}{2} + k\log\frac{n+1}{\sqrt{k}}.$$

If $k_0 \le k \le k_0+1$ we obtain $A'_{(n,n)}(k) = \log\frac{n+1}{\sqrt{k}}$ and $A''_{(n,n)}(k) = -\frac{1}{2k}$.

## 2.3 Areas of essential regions in the $\hat{x}\hat{y}$ plane

Consider in the $xy$ plane, an essential region $R_{(n,n')}\left(n \le n'\right)$ of the $xy = k\left(k \ge 4\right)$ hyperbola and $\psi$ an $\mathbb{R}^+$ prime coding function. Let $\widehat{R}_{(n,m)}$ be the corresponding region in the $\hat{x}\hat{y}$ plane that is, $\widehat{R}_{(n,n')} = \left(\psi\times\psi\right)\left(R_{(n,n')}\right)$. We call $\widehat{A}_{(n,n')}$ the area of $\widehat{D}_{(n,n')} = \left(\psi\times\psi\right)\left(D_{(n,n')}\right)$ considering the $\hat{x}\hat{y}$ plane embedded in the $xy$ plane. Then

$$\hat{A}_{(n,n')} = \iint_{\hat{D}_{(n,n')}} d\hat{x}d\hat{y}$$

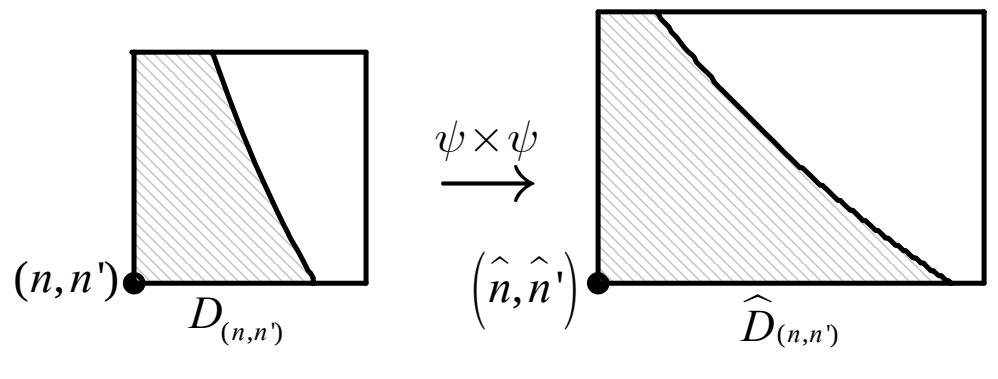

(Fig. 29)

The transformation that maps $D_{(n,n')}$ in $\widehat{D}_{(n,n')}$ is $\begin{cases} \hat{x} = \psi_n(x) \\ \hat{y} = \psi_{n'}(y) \end{cases}$

The Jacobian for this transformation is $J = \begin{vmatrix} \frac{\partial \hat{x}}{\partial x} & \frac{\partial \hat{x}}{\partial y} \\ \frac{\partial \hat{y}}{\partial x} & \frac{\partial \hat{y}}{\partial y} \end{vmatrix} = \begin{vmatrix} \psi_n'(x) & 0 \\ 0 & \psi_{n'}'(y) \end{vmatrix} = \psi_n'(x)\psi_{n'}'(y)$

Thus, $\hat{A}_{(n,n')} = \iint_{\widehat{D}_{(n,n')}} d\hat{x}d\hat{y} = \iint_{D_{(n,n')}} \left|\psi_n'(x)\psi_{n'}'(y)\right| dxdy$

Since $\psi$ is an $\mathbb{R}^+$ prime coding function, then $|J| = \xi_n \xi_{n'}$ and as a result the relationship between the areas of the essential regions in $xy$ and in $\hat{x}\hat{y}$ is

$$\hat{A}_{(n,n')} = \iint_{D_{(n,n')}} \xi_n \xi_{n'} \, dxdy = \xi_n \xi_{n'} A_{(n,n')}$$

## 2.4 $\widehat{A_I}, \widehat{A_S}, \widehat{A_T}$ functions

Let $\alpha$ be an even number. We will assume, for technical reasons that $\alpha \geq 16$.
Let $k \in [4, \alpha/2]$. Consider the subsets of $\mathbb{R}^2$

$$\text{i) } D_I(k) = \left\{(x,y) \in \mathbb{R}^2 : x \geq 2, y \geq x, xy \leq k\right\}$$

$$\text{ii) } D_S(k) = \left\{(x,y) \in \mathbb{R}^2 : x \geq 2, y \geq x, \alpha - k \leq xy \leq \alpha - 4\right\}$$

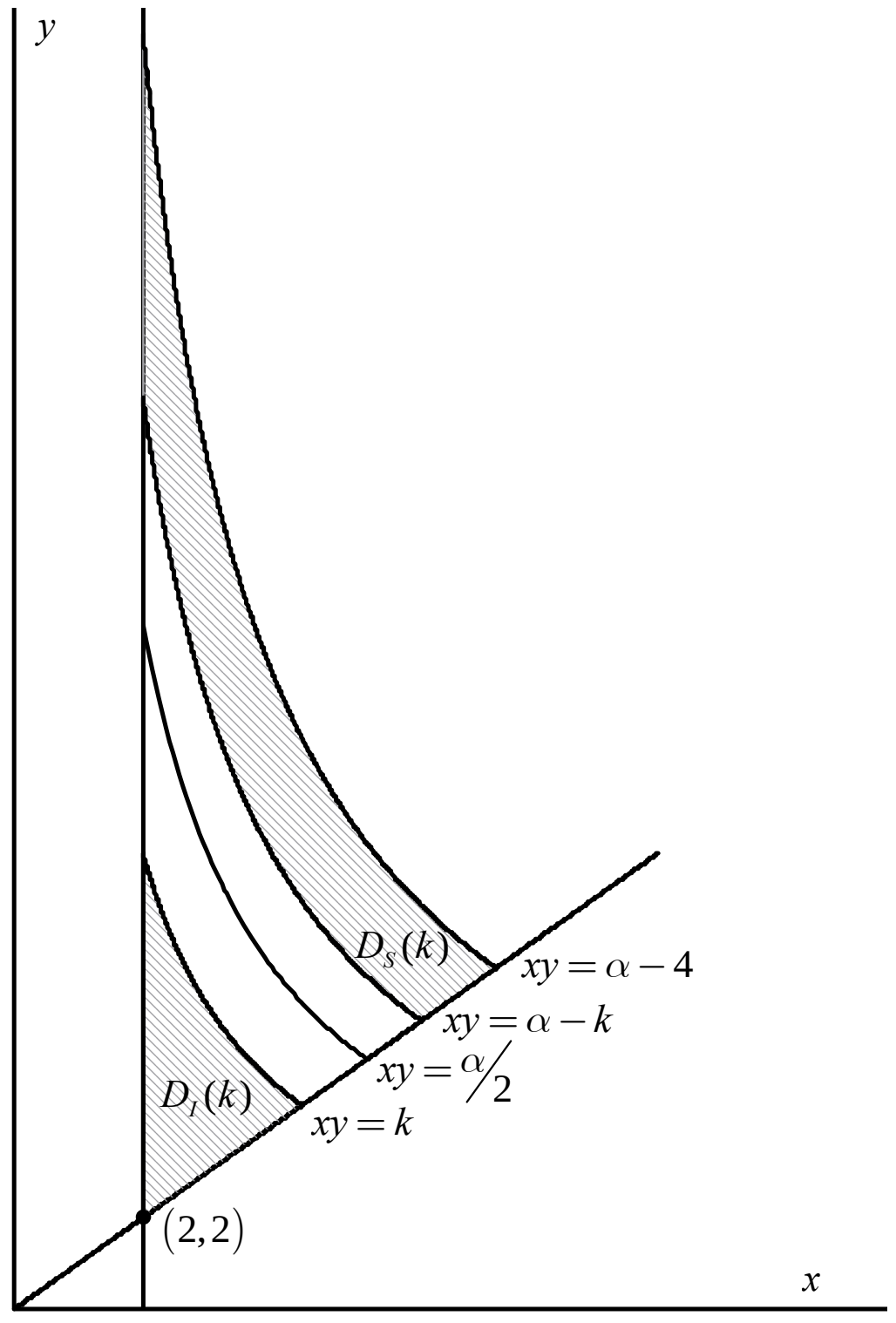

(Fig. 30)

Be an $\mathbb{R}^+$ prime coding function and consider the subsets of $\left[0, M_\psi\right)^2$

$$\widehat{D_I}\left(\hat{k}\right) = (\psi \times \psi)\left(D_I(k)\right), \widehat{D_S}\left(\hat{k}\right) = (\psi \times \psi)\left(D_S(k)\right)$$

We now define the functions

1) $\widehat{A_I} : \left[\hat{4}, \widehat{\alpha \div 2}\right] \to \mathbb{R}^+,\ \hat{k} \to \widehat{A_I}\left(\hat{k}\right)$ (Area of $\widehat{D_I}\left(\hat{k}\right)$ )

2) $\widehat{A_S} : \left[\hat{4}, \widehat{\alpha \div 2}\right] \to \mathbb{R}^+,\ \hat{k} \to \widehat{A_S}\left(\hat{k}\right)$ (Area of $\widehat{D_S}\left(\hat{k}\right)$ )

3) $\widehat{A_T} : \left[\hat{4}, \widehat{\alpha \div 2}\right] \to \mathbb{R}^+,\ \widehat{A_T} = \widehat{A_I} + \widehat{A_S}$

## 2.5 $\left(\widehat{A_I}\right)''$ function

Let $\alpha$ be an even number$(\alpha \geq 16)$. We take $k_0 = 4, 5, ..., \frac{\alpha}{2} - 1$ and we study the second derivative of $\widehat{A_I}$ at each closed interval$\left[\hat{k}_0, \hat{k}_0 \oplus \hat{1}\right]$ . For this, we consider the corresponding function $A_I(k)$ (Area of $D_I(k)$). Then $\forall k \in \left[k_0, k_0 + 1\right]$ we verify $A_I(k) = A_I(k_0) + A_I(k) - A_I(k_0)$. Additionally, $A_I(k) - A_I(k_0)$ is the sum of the areas in the essential regions associated with the $xy = k$ hyperbola, minus the area in the essential regions associated with the $xy = k_0$ hyperbola.

Then, $A_I(k) - A_I(k_0) = \sum_{(n,i(n)) \in E_S(k_0)} \left[A_{n,i(n)}(k) - A_{n,i(n)}(k_0)\right]$.

We know that functions $A_{(n,i(n))}(k)$ have a second derivative in $\left[k_0, k_0 + 1\right]$, therefore

$$A_I''(k) = \sum_{(n,i(n)) \in E_S(k_0)} A''_{(n,i(n))}(k) \quad \left(\forall k \in [k_0, k_0 + 1]\right)$$

We now want to find the expression of $\left(\widehat{A_I}\right)''$ as a function of the variable $\hat{k}$, where $\hat{k} \in \left[\hat{k}_0, \hat{k}_0 \oplus \hat{1}\right]$. In section 2.3 we proved that

$$\widehat{A}_{(n,n')}\left(\hat{k}\right) = \xi_n \xi_{n'} A_{(n,n')}(k) \text{ (when we have an } \mathbb{R}^+ \text{ prime coding function).}$$

If we derive with respect to $\hat{k}$, we obtain$\left(\widehat{A}_{(n,n')}\right)'\left(\hat{k}\right) = \xi_n \xi_{n'} A'_{(n,n')}(k) \frac{dk}{d\hat{k}}$

(Since at $k \in [k_0, k_0 + 1]$, the expression of $\hat{k}$ is $\hat{k} = \xi_{k_0}(k - k_0) + B_{k_0}$ (from 1.6))

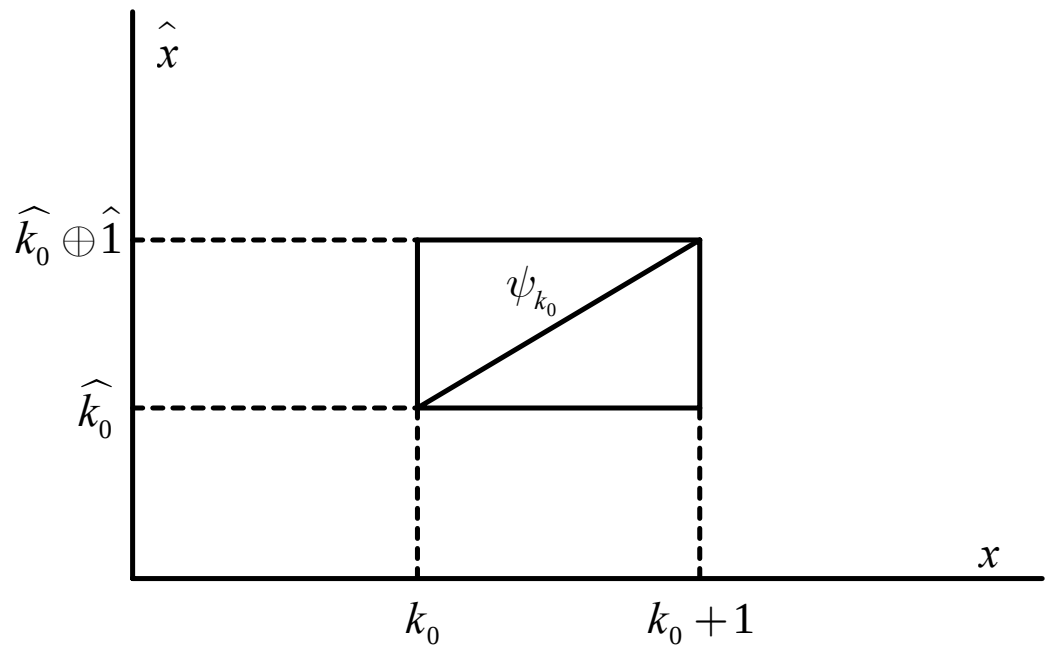

(Fig. 31)

Then $\frac{dk}{d\hat{k}} = \frac{1}{\xi_{k_0}}$ , therefore $\left(\hat{A}_{(n,n')}\right)'\left(\hat{k}\right) = \frac{\xi_n \xi_{n'}}{\xi_{k_0}} \cdot A'_{(n,n')}\left(k\right)$

Deriving once again

$$\left(\hat{A}_{(n,n')}\right)''\left(\hat{k}\right) = \frac{\xi_n \xi_{n'}}{\xi^2_{k_0}} \cdot A''_{(n,n')}\left(k\right)$$

Combining all this with 2.2.1, we obtain the following proposition:

**2.5.1 Proposition** Let $\alpha$ be an even number$\left(\alpha \geq 16\right)$, then

a)$\left(\widehat{A_I}\right)''\left(\hat{k}\right) = \sum_{(n,i(n))\in E_S(k_0)} \left(\hat{A}_{(n,i(n))}\right)''(\hat{k})$ , $\forall \hat{k} \in \left[\hat{k}_0, \hat{k}_0 \oplus \hat{1}\right]$, $\hat{k}_0 \in \left\{\hat{4}, \hat{5}, ..., \hat{\alpha} \div \hat{2} \sim \hat{1}\right\}$

b) For $\left[\hat{k}_0, \hat{k}_0 \oplus \hat{1}\right]$ where $\hat{k}_0 \in \left\{\hat{4}, \hat{5}, ..., \hat{\alpha} \div \hat{2} \sim \hat{1}\right\}$ and bearing in mind the different types of essential regions, we obtain

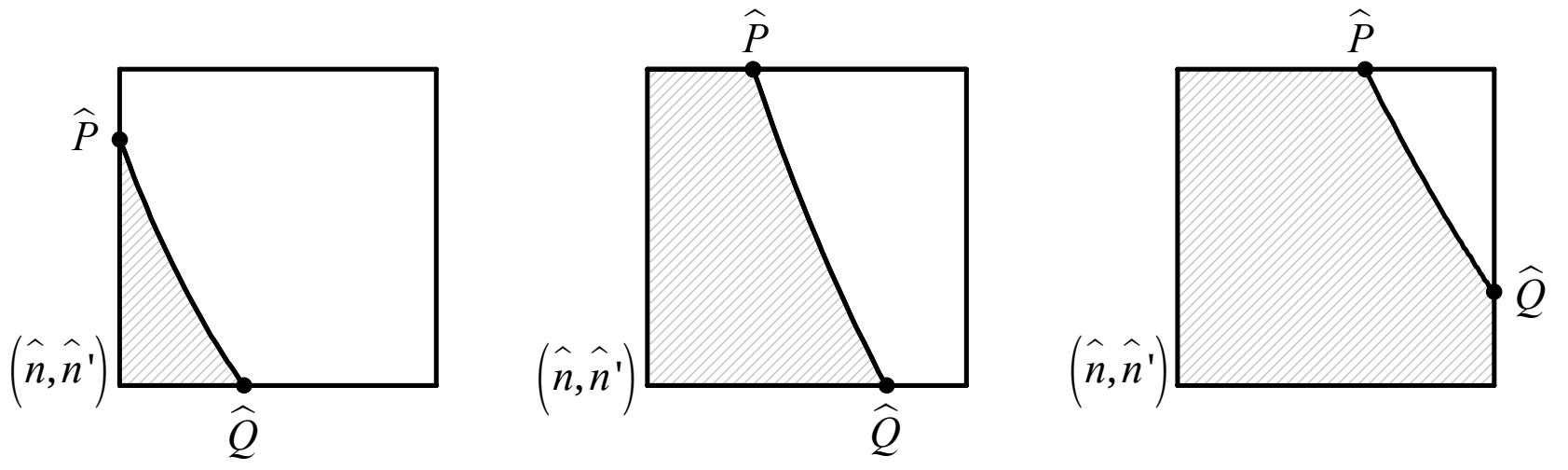

$\left(\hat{A}_{(n,n')}\right)''(\hat{k}) = \frac{\xi_n \xi_{n'}}{\xi^2_{k_0}} \cdot \frac{1}{k}$ $\quad$ $\left(\hat{A}_{(n,n')}\right)''(\hat{k}) = 0$ $\quad$ $\left(\hat{A}_{(n,n')}\right)''(\hat{k}) = -\frac{\xi_n \xi_{n'}}{\xi^2_{k_0}} \cdot \frac{1}{k}$

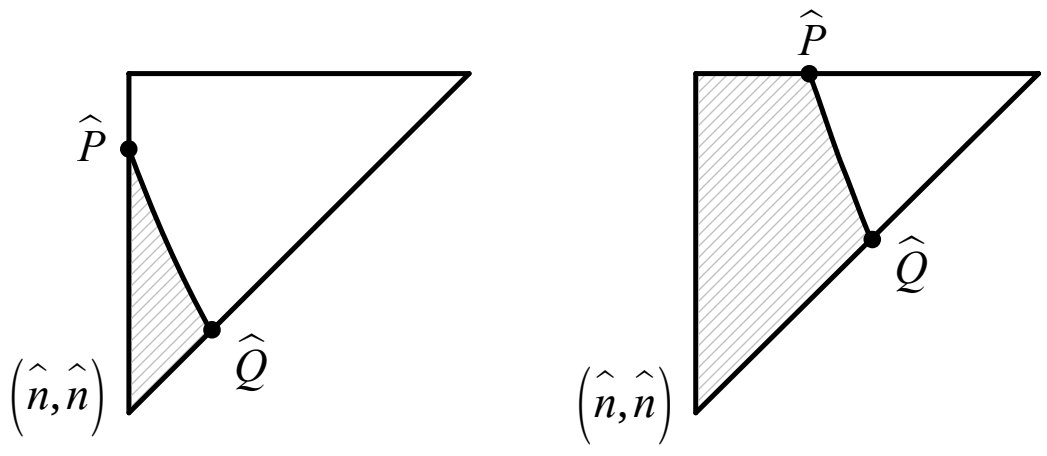

$\left(\hat{A}_{(n,n)}\right)''(\hat{k}) = \frac{\xi^2_n}{\xi^2_{k_0}} \cdot \frac{1}{2k}$ $\quad$ $\left(\hat{A}_{(n,n)}\right)''(\hat{k}) = -\frac{\xi^2_n}{\xi^2_{k_0}} \cdot \frac{1}{2k}$

(Fig. 32)

**2.5.2 Example** We will find $\left(\hat{A}_I\right)''(\hat{k})$ in $\left[\widehat{12}, \widehat{13}\right]$ ($\hat{\alpha} \geq \widehat{26}$)

$$\left(\hat{A}_I\right)''(\hat{k}) = \frac{\xi_2 \xi_6}{\xi^2_{12}} \cdot \frac{1}{k} - \frac{\xi_2 \xi_4}{\xi^2_{12}} \cdot \frac{1}{k} + \frac{\xi_3 \xi_4}{\xi^2_{12}} \frac{1}{k} - \frac{1}{2} \cdot \frac{\xi^2_3}{\xi^2_{12}} \cdot \frac{1}{k} \quad (\forall \hat{k} \in \left[\widehat{12}, \widehat{13}\right])$$

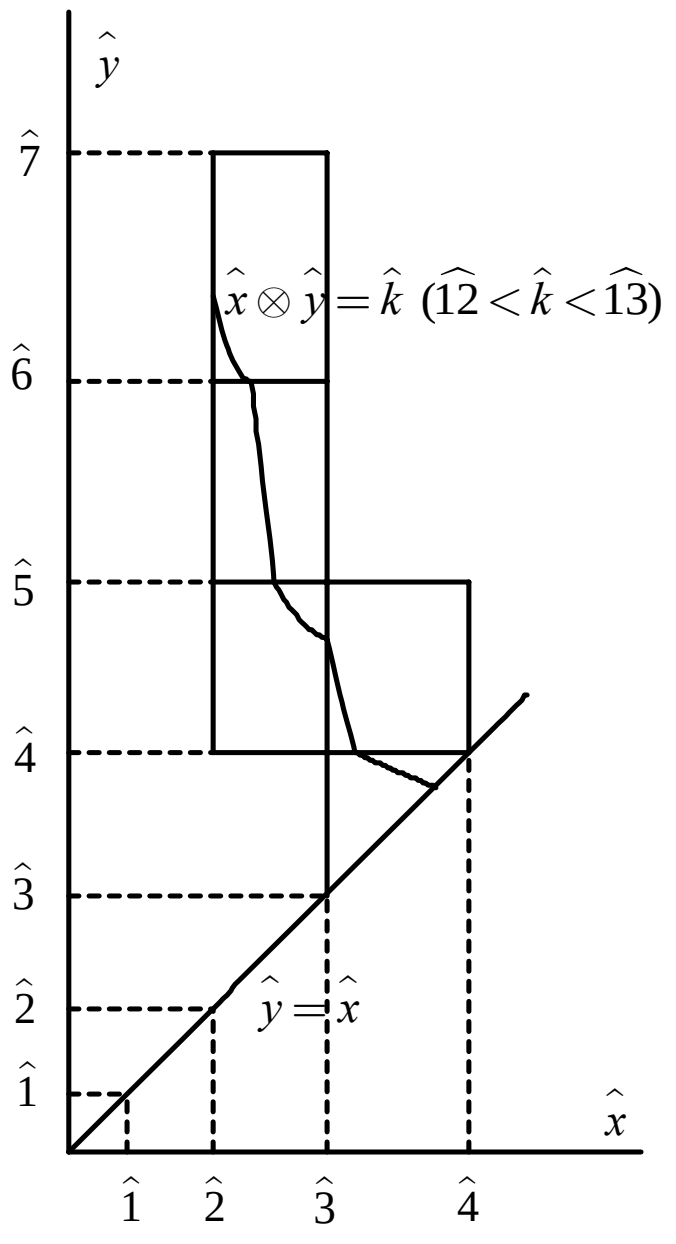

(Fig. 33)

Removing the common factor $\frac{1}{\xi_{12}^2}\cdot\frac{1}{k}$ then

$$\left(\widehat{A}_I\right)''(\hat{k})=\frac{1}{\xi_{12}^2\cdot k}\left(\xi_2\xi_6-\xi_2\xi_4+\xi_3\xi_4-\frac{1}{2}\xi_3^2\right)$$

Consider the polynomial $p\left(x_2,x_3,...,x_6\right)=x_2x_6-x_2x_4+x_3x_4-\frac{1}{2}x_3^2$.

We call this polynomial a lower essential polynomial of $k_0=12$ and we write it as $P_{I,k_0}$.

**2.5.3 Definition** Let $\alpha$ be an even number$(\alpha\ge 16)$. The polynomial obtained naturally by removing the common factor function $\frac{1}{\xi_{k_0}^2}\cdot\frac{1}{k}$ in $\left(\widehat{A}_I\right)''(\hat{k})$ in the interval$\left[\hat{k}_0,\hat{k}_0\oplus\hat{1}\right]$ $\left(k_0=4,5,...,\frac{\alpha}{2}-1\right)$ is called a lower essential polynomial of $k_0$. It is written as $P_{I.k_0}$.

**2.5.4 Remarks** i) Lower essential polynomials are homogeneous polynomials of degree 2. ii) The variables that intervene in $P_{I,k_0}$ are at most $x_n,x_{i(n)}$ and $\left(n,i(n)\right)\in E_s\left(k_0\right)$, some of which may be missing (those which correspond to essential regions in which the second derivative is $0$). iii) We will use $P_{I.k_0}$ as the coefficient of $\frac{1}{\xi_{k_0}^2}\cdot\frac{1}{k}$ in $\left(\widehat{A}_I\right)''(\hat{k})$

**2.5.5 Corollary** Let $\alpha$ be an even number$(\alpha\ge 16)$. Then,

$\forall\hat{k}\in\left[\hat{k}_0,\hat{k}_0\oplus\hat{1}\right]$, ( $\hat{k}_0\in\left\{\hat{4},\hat{5},...,\hat{\alpha}\div\hat{2}\sim\hat{1}\right\}$ ) we verify $\left(\widehat{A}_I\right)''(\hat{k})=\frac{1}{\xi_{k_0}^2 k}\cdot P_{I,k_0}$

## 2.6 $\hat{A}''_S$, $\hat{A}''_T$ functions

Let $\alpha$ be an even number$(\alpha \geq 16)$. We take $k_0 \in \left\{4, 5, ..., \frac{\alpha}{2} - 1\right\}$ and we examine the second derivative of $\hat{A}_S$ at each closed interval$\left[\hat{k}_0, \hat{k}_0 \oplus \hat{1}\right]$. For this we take the function. $A_S(k)$(Area of $D_S(k)$) Then, $\forall k \in [k_0, k_0 + 1]$ we verify

$$A_s(k) = A_s(k_0) + A_s(k) - A_s(k_0).$$

Additionally, $A_S(k) - A_S(k_0)$ is the area included between the curves

$$xy = \alpha - k_0, xy = \alpha - k, x = 2, y = x.$$

As a result, it is the sum of the areas in the essential regions of the $xy = \alpha - k_0$ hyperbola, minus the area in the essential regions of $xy = \alpha - k$. We obtain

$$A_S(k) - A_S(k_0) = \sum_{(n,i(n)) \in E_S(\alpha - k_0 - 1)} \left[A_{(n,i(n))}(\alpha - k_0) - A_{(n,i(n))}(\alpha - k)\right].$$

We derive

$$A''_S(k) = -\sum_{(n,i(n)) \in E_S(\alpha - k_0 - 1)} A''_{(n,i(n))}(\alpha - k)$$

Of course, the same relationships as in the lower areas are maintained with the expression $\hat{A}''_S$ as a function of $\hat{k}$. We are left with

$$\left(\hat{A}_{(n,i(n))}\right)''(\hat{k}) = -\frac{\xi_n \xi_{n'}}{\xi^2_{\alpha - k_0 - 1}} \cdot A''_{(n,i(n))}(\alpha - k)$$

We define upper essential polynomial in the same way we defined lower essential polynomial and we write them as $P_{S,k_0}$. The same remarks are maintained.

**2.6.1 Remarks** i) Upper essential polynomials are homogeneous polynomials of degree 2. ii) The variables that intervene in $P_{S,k_0}$ are, at most, $x_n, x_{i(n)}$ and $(n, i(n)) \in E_S(\alpha - k_0 - 1)$, some of which may be missing (those which correspond to essential regions in which the second derivative is $0$). iii) We will also use $P_{S,k_0}$ as the coefficient of $\frac{1}{\xi^2_{\alpha - k_0 - 1}} \cdot \frac{1}{\alpha - k}$ in $A''_S(k)$.

**2.6.2 Definition** Let $\alpha$ be an even number where $\alpha \geq 16$. Let $k_0 \in \left\{4, 5, ..., \frac{\alpha}{2} - 1\right\}$ and consider an $\mathbb{R}^+$ prime coding function**.** For every $k \in [k_0, k_0 + 1]$ we define

$$A_{k_0}(k) = \frac{1}{\xi^2_{k_0} k}, \quad B_{k_0}(k) = \frac{1}{\xi^2_{\alpha - k_0 - 1}} \cdot \frac{1}{\alpha - k}$$

We obtain the following theorem:

**2.6.3 Theorem** Let $\alpha$ be an even number where $\alpha \geq 16$. Consider an $\mathbb{R}^+$ prime coding function and $k_0 \in \left\{4, 5, ..., \frac{\alpha}{2} - 1\right\}$. Then $\forall \hat{k} \in \left[\hat{k}_0, \hat{k}_0 \oplus \hat{1}\right]$ we verify

$$\left(\widehat{A_T}\right)''\left(\hat{k}\right) = A_{k_0}\left(k\right) P_{I,k_0} + B_{k_0}\left(k\right) P_{S,k_0}$$

Note that $k$ is a function of $\hat{k} : k = k\left(\hat{k}\right)$

## 2.7 Bounds for the coefficients of essential polynomials

Let $\psi$ be an $\mathbb{R}^+$ prime coding with coefficients $\xi_i$, for an $\alpha$ even number where $\alpha \geq 16$ and $k_0 = 4, 5, ..., \frac{\alpha}{2} - 1$, we obtain the following variation for $\alpha - k_0 - 1$

| $k_0$ | 4 | 5 | 6 | $\cdots$ | $\alpha/2 - 1$ |
|---|---|---|---|---|---|
| $\alpha - k_0 - 1$ | $\alpha - 5$ | $\alpha - 6$ | $\alpha - 7$ | $\cdots$ | $\alpha/2$ |

We know that the expression of $\left(\widehat{A_T}\right)''$ in the interval $\left[\hat{k}_0, \hat{k}_0 \oplus \hat{1}\right]$ is

$$\left(\widehat{A_T}\right)''\left(\hat{k}\right) = A_{k_0}\left(k\right) P_{I,k_0} + B_{k_0}\left(k\right) P_{S,k_0} \left(k = k\left(\hat{k}\right)\right)$$

We set the bounds of $A_{k_0}\left(k\right)$ and $B_{k_0}\left(k\right)$ in the interval $\left[k_0, k_0 + 1\right]$.

$$\frac{1}{\left(k_0 + 1\right)\xi_{k_0}^2} \leq A_{k_0}(k) \leq \frac{1}{k_0 \xi_{k_0}^2}, \quad \frac{1}{\left(\alpha - k_0\right)\xi_{\alpha - k_0 - 1}^2} \leq B_{k_0}(k) \leq \frac{1}{\left(\alpha - k_0 - 1\right)\xi_{\alpha - k_0 - 1}^2}$$

Additionally, since $k_0 \leq \frac{\alpha}{2} - 1$, we obtain $2k_0 \leq \alpha - 2, k_0 + 1 \leq \alpha - k_0 - 1$.

Furthermore, since $\frac{\alpha}{2} - \frac{1}{2} > k_0$ we obtain $\alpha - 1 > 2k_0, \alpha - k_0 - 1 > k_0, \xi_{\alpha - k_0 - 1}^2 > \xi_{k_0}^2$.

Thus,

$$\left(\alpha - k_0 - 1\right)\xi_{\alpha - k_0 - 1}^2 > \left(k_0 + 1\right)\xi_{k_0}^2 \Rightarrow A_{k_0}(k) \geq \frac{1}{\left(k_0 + 1\right)\xi_{k_0}^2} > \frac{1}{\left(\alpha - k_0 - 1\right)\xi_{\alpha - k_0 - 1}^2} \geq B_{k_0}(k)$$

In other words, $\forall k \in \left[k_0, k_0 + 1\right]$ we verify $B_{k_0}\left(k\right) < A_{k_0}\left(k\right)$.

We will now see how $A_{k_0}(k)$ and $B_{k_0}(k)$ change when moving along different intervals. We obtain

$$A_{k_0 + 1}(k) \leq \frac{1}{\left(k_0 + 1\right)\xi_{k_0 + 1}^2} < \frac{1}{\left(k_0 + 1\right)\xi_{k_0}^2} \leq A_{k_0}(k)$$

Additionally, $B_{k_0 + 1}(k) \geq \frac{1}{\left(\alpha - k_0 - 1\right)\xi_{\alpha - k_0 - 2}^2} > \frac{1}{\left(\alpha - k_0 - 1\right)\xi_{\alpha - k_0 - 1}^2} \geq B_{k_0}(k)$

Finally,

$$B_{\frac{\alpha}{2} - 1}(k) \leq \frac{1}{\frac{\alpha}{2} \cdot a_{\alpha/2}^2} < \frac{1}{\frac{\alpha}{2} \cdot a_{\alpha/2 - 1}^2} \leq A_{\frac{\alpha}{2} - 1}(k)$$

**2.7.1 Proposition** Let $\alpha$ be an even number where $\alpha \geq 16$ and $\psi$ an $\mathbb{R}^{+}$ prime coding. We call $m_{A_{k_0}} = \inf\left\{A_{k_0}(k) : k \in [k_0, k_0+1]\right\}$, or, for short,

$$M_{A_{k_0}} = \sup\left\{A_{k_0}(k)\right\}, m_{B_{k_0}} = \inf\left\{B_{k_0}(k)\right\},\ M_{B_{k_0}} = \sup\left\{B_{k_0}(k)\right\}$$

Then , we verify

$$m_{B_4} < M_{B_4} < \ldots < m_{B_{\frac{\alpha}{2}-1}} < M_{B_{\frac{\alpha}{2}-1}} < m_{A_{\frac{\alpha}{2}-1}} < M_{A_{\frac{\alpha}{2}-1}} < \ldots < m_{A_4} < M_{A_4}$$

## 2.8 Signs of the essential point coordinates

**2.8.1 Definition** Given an $\mathbb{R}^{+}$ prime coding and $\alpha$ an even number $(\alpha \geq 16)$ we write $P_{k_0} = \left(x_{k_0}, y_{k_0}\right) = \left(P_{I,k_0}, P_{S,k_0}\right), \forall k_0 \in \left\{4, 5, \ldots, \frac{\alpha}{2} - 1\right\}$. We call $P_{k_0}$ the essential points associated with. $\psi$ .

The formula from proposition 2.5.1 is

$$\left(\widehat{A_I}\right)''\left(\hat{k}\right) = \sum_{(n,i(n)) \in E_S(k_0)} \left(\widehat{A}_{(n,i(n))}\right)''(\hat{k}) \quad (\forall \hat{k} \in \left[\hat{k}_0, \hat{k}_0 \oplus \hat{1}\right],\ \hat{k}_0 \in \left\{\hat{4}, \hat{5}, \ldots, \widehat{\alpha} \div \hat{2} \sim \hat{1}\right\})$$

where the $E_S(k_0)$ sub-indexes are

$$(1)\begin{cases} n = 2, 3, \ldots, \left[\sqrt{k_0}\right] - 1 \\ i(n) = \left[\dfrac{k_0}{n+1}\right], \ldots, \left[\dfrac{k_0}{n}\right] \end{cases} \qquad (2)\begin{cases} n = \left[\sqrt{k_0}\right] \\ i(n) = \left[\sqrt{k_0}\right], \ldots, \left[\dfrac{k_0}{\left[\sqrt{k_0}\right]}\right] \end{cases}$$

Thus, for sub-index *n* in (1), in $\left(\widehat{A_I}\right)''$ only intervene $i(n) = \left[\dfrac{k_0}{n+1}\right]$ and $\left[\dfrac{k_0}{n}\right]$, since we have already seen that all the sub-indexes included between them two, $\left(\widehat{A}_{(n,n')}\right)''(\hat{k}) = 0$ , as the regions are of type 3. In the lower essential polynomial we obtain $\xi_n(\xi_{\left[\frac{k_0}{n}\right]} - \xi_{\left[\frac{k_0}{n+1}\right]}) > 0$ (for an $\mathbb{R}^{+}$ prime coding ). For $n = \left[\sqrt{k_0}\right]$ we obtain the cases:

$$\text{(i)} \quad \left[\sqrt{k_0}\right] = \left[\frac{k_0}{\left[\sqrt{k_0}\right]}\right] \qquad \text{(ii)} \ \left[\sqrt{k_0}\right] < \left[\frac{k_0}{\left[\sqrt{k_0}\right]}\right]$$

In case (i) we would obtain the addend $\frac{1}{2}\xi^2_{\left[\sqrt{k_0}\right]}$.

In case (ii) we would obtain $\xi_{\left[\sqrt{k_0}\right]}\xi_{\left[\frac{k_0}{\left[\sqrt{k_0}\right]}\right]} - \frac{1}{2}\xi^2_{\left[\sqrt{k_0}\right]} = \xi_{\left[\sqrt{k_0}\right]}(\xi_{\left[\frac{k_0}{\left[\sqrt{k_0}\right]}\right]} - \frac{1}{2}\xi_{\left[\sqrt{k_0}\right]}) > 0$

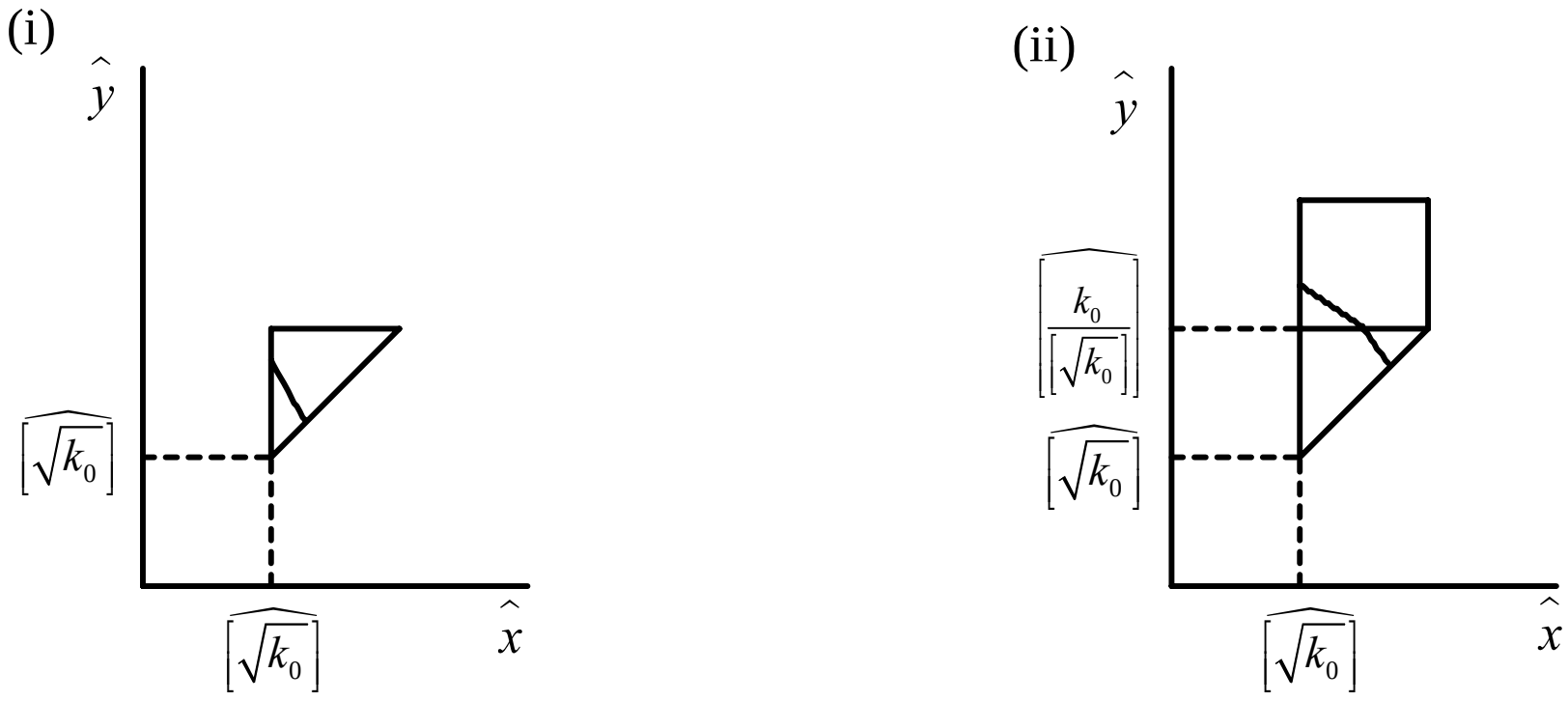

(Fig. 34)

As a result, for an $\mathbb{R}^+$ prime coding we obtain $x_4 > 0, x_5 > 0, ..., x_{\alpha/2-1} > 0$. The reasoning is entirely analogous for the upper essential polynomials that is, $y_4 < 0, y_5 < 0, ..., y_{\alpha/2-1} < 0$. We will now arrange the coordinates for the essential points.

*1.- Lower essential polynomials*

If $k_0$ is not prime, there is at least one natural coordinate point $(\hat{n}, \hat{n}')$ such that $\hat{2} \le \hat{n} \le \hat{n}'$ which the $\hat{x} \otimes \hat{y} = \hat{k}_0$ "hyperbola" goes through.

If $2 < n < n'$ we obtain

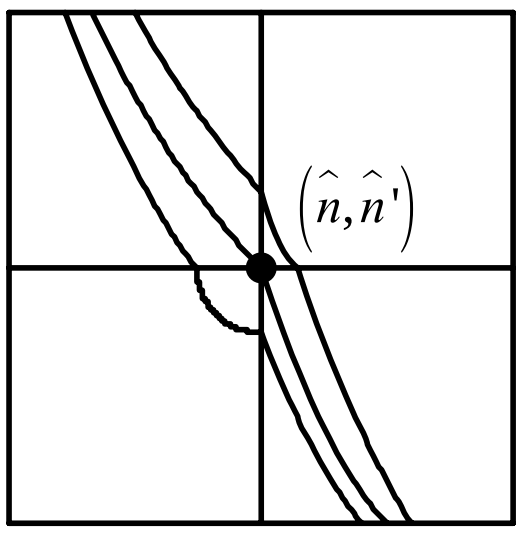

Transformations

| $P_{I,k_0-1}$ | $P_{I,k_0}$ |
|---|---|
| 0 | $-\xi_{n-1}\xi_{n'}$ |
| $-\xi_{n-1}\xi_{n'-1}$ | $\xi_n\xi_{n'}$ |
| $\xi_n\xi_{n'-1}$ | 0 |

(Fig. 35)

If $2 < n = n'$ we obtain

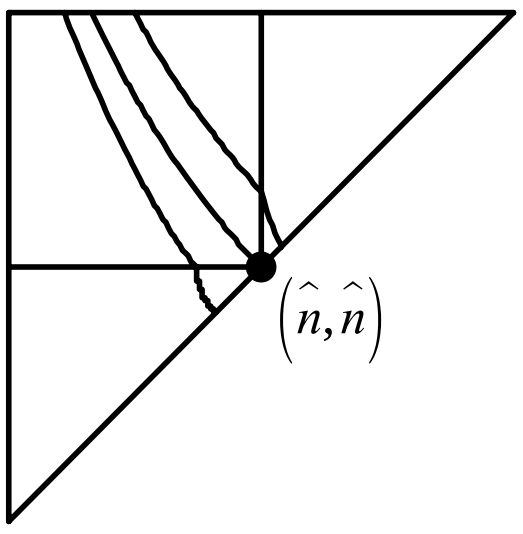

Transformations

| $P_{I,k_0-1}$ | $P_{I,k_0}$ |
|---|---|
| 0 | $-\xi_{n-1}\xi_n$ |
| $-\frac{1}{2}\xi_{n-1}^2$ | $\frac{1}{2}\xi_n^2$ |

(Fig. 36)

If $2 = n < n'$ we obtain

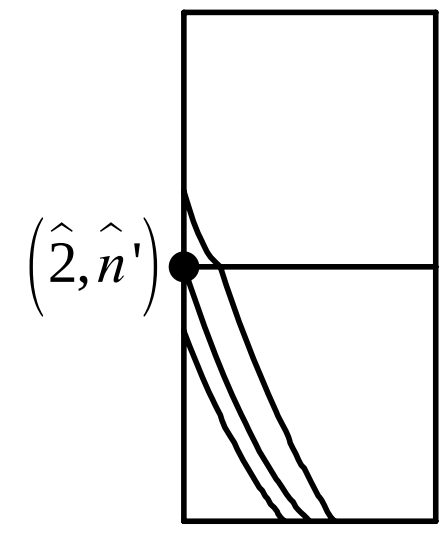

Transformations

| $P_{I,k_0-1}$ | $P_{I,k_0}$ |
|---|---|
| 0 | $\xi_2\xi_{n'}$ |
| $\xi_2\xi_{n'-1}$ | 0 |

(Fig. 37)

Then $P_{I,k_0} - P_{I,k_0-1} > 0$, since where there are transformations we obtain, for an $\mathbb{R}^+$ prime coding, either

$$\xi_n\xi_{n'} - \xi_{n-1}\xi_{n'} + \xi_{n-1}\xi_{n'-1} - \xi_n\xi_{n'-1} = \xi_{n'}\left(\xi_n - \xi_{n-1}\right) - \xi_{n'-1}\left(\xi_n - \xi_{n-1}\right) = \left(\xi_n - \xi_{n-1}\right)\left(\xi_{n'} - \xi_{n'-1}\right) > 0$$

or

$$\frac{1}{2}\xi_n^2 - \xi_{n-1}\xi_n + \frac{1}{2}\xi_{n-1}^2 = \frac{1}{2}\left(\xi_n^2 - 2\xi_{n-1}\xi_n + \xi_{n-1}^2\right) = \frac{1}{2}\left(\xi_n - \xi_{n-1}\right)^2 > 0$$

or

$$\xi_2\xi_{n'} - \xi_2\xi_{n'-1} = \xi_2\left(\xi_{n'} - \xi_{n'-1}\right) > 0$$

If $k_0$ is prime then $P_{I,k_0-1} = P_{I,k_0}$ since the same essential regions exist for the hyperbolas $\hat{x} \otimes \hat{y} = \hat{k}$ in $\left(\hat{k}_0 \sim \hat{1}, \hat{k}_0\right) \cup \left(\hat{k}_0, \hat{k}_0 \oplus \hat{1}\right)$.

*2.- Upper essential polynomials*

If $\hat{\alpha} \sim \hat{k}_0$ is not prime, and reasoning in the same way, we obtain $P_{S,k_0} - P_{S,k_0-1} > 0$. For $\widehat{\alpha - k_0}$ prime we obtain $P_{S,k_0} = P_{S,k_0-1}$ since the same essential regions exist for the hyperbolas $\hat{x} \otimes \hat{y} = \hat{\alpha} - \hat{k}$ in $\left(\hat{k}_0 \sim \hat{1}, \hat{k}_0\right) \cup \left(\hat{k}_0, \hat{k}_0 \oplus \hat{1}\right)$ ($k_0 \in \left\{5, 6, ..., \frac{\alpha}{2} - 1\right\}$). We obtain the theorem:

**2.8.2 Theorem** Let $\alpha$ be an even number$\left(\alpha \geq 16\right)$, and $\psi$ an $\mathbb{R}^+$ prime coding. Let $P_{k_0} = \left(x_{k_0}, y_{k_0}\right)$ be the essential points. Then:

i) $0 < x_4 \leq x_5 \leq ... \leq x_{\alpha/2-1}$. Additionally, $x_{k_0-1} = x_{k_0} \Leftrightarrow k_0$ is prime.

ii) $y_4 \leq y_5 \leq ... \leq y_{\alpha/2-1} < 0$. Additionally, $y_{k_0-1} = y_{k_0} \Leftrightarrow \alpha - k_0$ is prime.

**2.8.3 Corollary** In the hypotheses from the above theorem, we obtain:

*"The even number $\alpha$ is the sum of two primes $k_0$ and $\alpha - k_0$ $\left(k_0 \in \left\{5, 6, ..., \frac{\alpha}{2} - 1\right\}\right)$ iff, the consecutive essential points $P_{k_0-1}$ and $P_{k_0}$ are repeated, that is $P_{k_0-1} = P_{k_0}$*

# 3. STATEMENT EQUIVALENT TO THE GOLDBACH CONJECTURE

## 3.1 Construction of the Goldbach Conjecture function $\mathfrak{G}$

**3.1.0 Quotation** "I have sometimes thought that the profound mystery which envelops our conceptions relative to prime numbers depends upon the limitations of our faculties in regard to time which, like space may be in essence poly-dimensional and that this and other such sort of truths would become self-evident to a being whose mode of perception is according to superficially as opposed to our own limitation to linearly extended time."

(J.J. Sylvester, from "On certain inequalities relating to prime numbers", *Nature* 38 (1888) 259-262, and reproduced in *Collected Mathematical Papers*, Volume 4, page 600 ,Chelsea, New York, 1973)

**3.1.1 Proposition** Let $\alpha$ be an even number $(\alpha \geq 16)$, and $\psi$ be an $\mathbb{R}^+$ prime coding with $\xi_i$ coefficients. Let $P_{k_0} = \left(x_{k_0}, y_{k_0}\right) \left(k_0 = 4,5,\ldots,\frac{\alpha}{2}-1\right)$ be the essential points. Then,

i) $x_{k_0}$ depends at most on $\xi_2, \xi_3, \ldots, \xi_{\left[\frac{k_0}{2}\right]}$

ii) $y_{k_0}$ depends at most on $\xi_2, \xi_3, \ldots, \xi_{\left[\frac{\alpha-k_0-1}{2}\right]}$

<u>*Proof*</u> i) In $E_S(k_0) = \{(n, i(n))\}$ we verify that $n \leq i(n)$. The smallest $n$ is $2$ and the biggest $i(n)$ is $\left[\frac{k_0}{2}\right]$. ii) In $E_S(\alpha - k_0 - 1) = \{(n, i(n))\}$ we verify that $n \leq i(n)$. The smallest $n$ is $2$ and the biggest $i(n)$ is $\left[\frac{\alpha - k_0 - 1}{2}\right]$.

**3.1.2. Consequence** Since the biggest sub-index coefficient that appears at the essential point coordinates is $\xi_{\left[\frac{\alpha-4-1}{2}\right]} = \xi_{\alpha/2-3}$ we conclude that, knowing the coefficients $\xi_2, \xi_3, \ldots, \xi_{\alpha/2-3}$ all the essential points are determined. Note that where $0 < \xi_2 < \xi_3 < \ldots < \xi_{\alpha/2-3}$ the 2.7.3 corollary is met. This leads to the following definition.

**3.1.3 Definition** Let $\alpha$ be an even number $(\alpha \geq 16)$, and $\psi$ an $\mathbb{R}^+$ prime coding such that its $\xi_i$ coefficients verify $0 < \xi_0 < \xi_1 < \xi_2 < \ldots < \xi_{\alpha/2-1}$. We say that $\psi$ a *prime* $\mathbb{R}^+$ *coding adapted to* $\alpha$.

Generally, $\left(\widehat{A_T}\right)''\left(\hat{k}_0 -\right) \neq \left(\widehat{A_T}\right)''\left(\hat{k}_0 +\right)$. The following proposition provides sufficient conditions for the $\left(\widehat{A_T}\right)''$ function to be well-defined and continuous in the $\left[\hat{4}, \hat{\alpha} \div \hat{2}\right]$ closed interval.

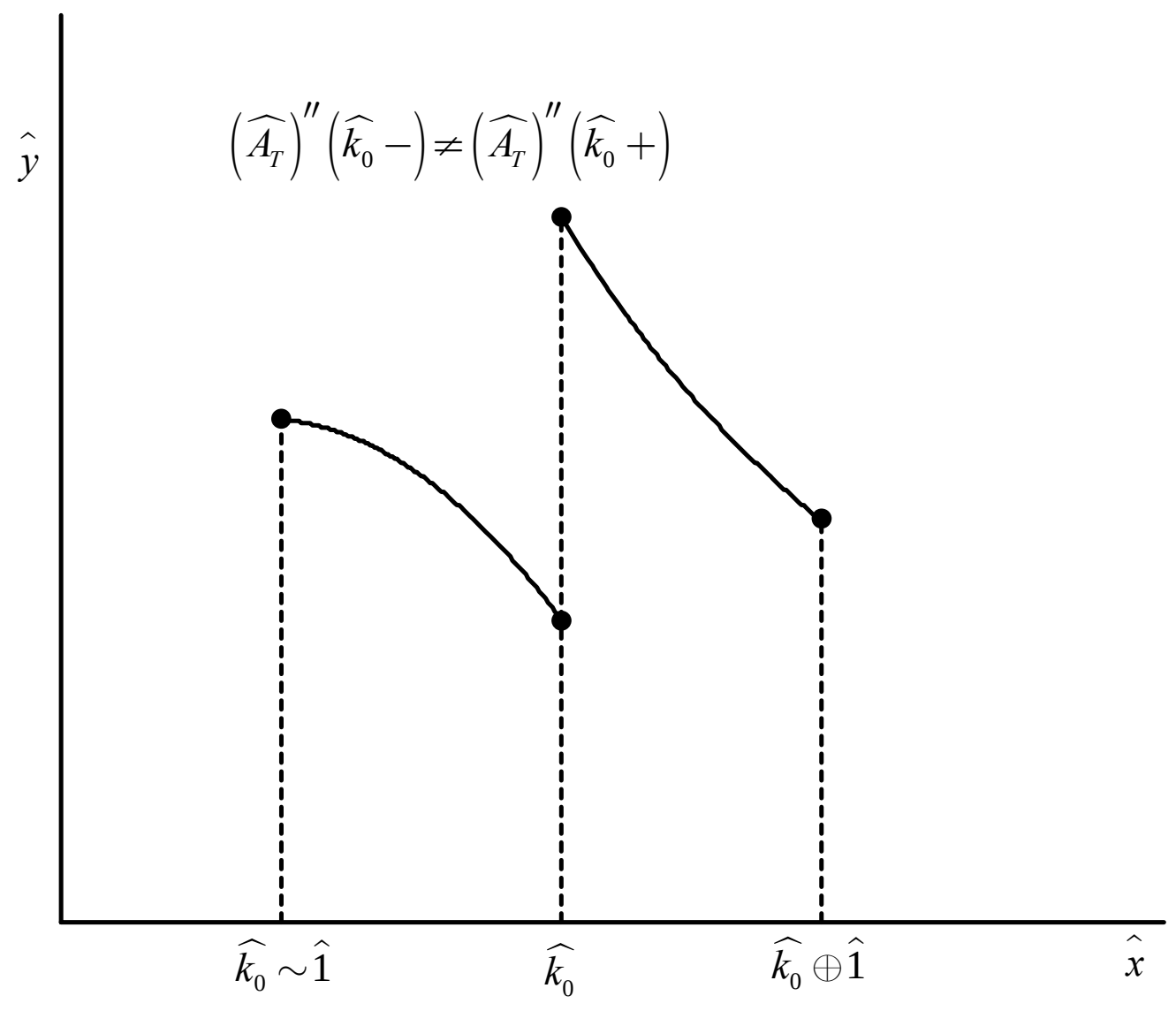

(Fig. 38)

**3.1.4 Proposition** Let $\alpha$ be an even number $(\alpha \geq 16)$ and $\psi$ be a prime $\mathbb{R}^+$ coding adapted to $\alpha$. Assume that

i) $\xi_{k_0}^2 x_{k_0-1} = \xi_{k_0-1}^2 x_{k_0} \wedge \xi_{\alpha-k_0-1}^2 y_{k_0-1} = \xi_{\alpha-k_0}^2 y_{k_0}$ $\forall k_0 \in \left\{5, 6, ..., \frac{\alpha}{2}-1\right\}$ where $k_0$ is non-prime.

ii) $\xi_{\alpha-p_0}^2 = \left|y_{p_0-1}\right|\left(\frac{\left|y_{p_0}\right|}{\xi_{\alpha-p_0-1}^2} + \frac{\alpha-p_0}{p_0}\cdot x_{p_0-1}\cdot\left(\frac{1}{\xi_{p_{0-1}}^2} - \frac{1}{\xi_{p_0}^2}\right)\right)^{-1}$ $\forall p_0 \in \left\{5, 6, ..., \frac{\alpha}{2}-1\right\}$ where $p_0$ is prime. Then

$$\left(\widehat{A_T}\right)''\left(\hat{k}_0-\right) = \left(\widehat{A_T}\right)''\left(\hat{k}_0+\right) \quad \left(\forall k_0 \in \left\{5, 6, ..., \frac{\alpha}{2}-1\right\}\right)$$

<u>*Proof*</u> $\left(\widehat{A_T}\right)''\left(\hat{k}_0-\right) = \frac{x_{k_0-1}}{k_0\xi_{k_0-1}^2} + \frac{y_{k_0-1}}{(\alpha-k_0)\xi_{\alpha-k_0}^2}$ , $\left(\widehat{A_T}\right)''\left(\hat{k}_0+\right) = \frac{x_{k_0}}{k_0\xi_{k_0}^2} + \frac{y_{k_0}}{(\alpha-k_0)\xi_{\alpha-k_0-1}^2}$

$$(\forall k_0 \in \left\{5, 6, ..., \frac{\alpha}{2}-1\right\})$$

Then, $\left(\widehat{A_T}\right)''\left(\hat{k}_0-\right) = \left(\widehat{A_T}\right)''\left(\hat{k}_0+\right) \Leftrightarrow \frac{1}{k_0}\left(\frac{x_{k_0-1}}{\xi_{k_0-1}^2} - \frac{x_{k_0}}{\xi_{k_0}^2}\right) = \frac{1}{\alpha-k_0}\left(\frac{\left|y_{k_0-1}\right|}{\xi_{\alpha-k_0}^2} - \frac{\left|y_{k_0}\right|}{\xi_{\alpha-k_0-1}^2}\right)$

When $k_0$ is non-prime, i) implies the equality above. Note that if $k_0$ is not prime, then $x_{k_0-1} < x_{k_0}$, and consequently, $\xi_{k_0-1}^2 < \xi_{k_0}^2$, in other words, it is consistent with the hypothesis that $\psi$ is a prime $\mathbb{R}^+$ coding adapted to $\alpha$.

If $p_0$ is prime, then $x_{p_0-1} = x_{p_0}$ therefore $\left(\widehat{A_T}\right)''\left(\widehat{p}_0-\right) = \left(\widehat{A_T}\right)''\left(\widehat{p}_0+\right)$ is equivalent to

$$\frac{\alpha-p_0}{p_0}x_{p_0-1}\cdot\left(\frac{1}{\xi_{p_0-1}^2} - \frac{1}{\xi_{p_0}^2}\right) = \frac{\left|y_{p_0-1}\right|}{\xi_{\alpha-p_0}^2} - \frac{\left|y_{p_0}\right|}{\xi_{\alpha-p_0-1}^2}$$

which in turn is equivalent to ii).

Now, let $\alpha$ be an even number $(\alpha \geq 16)$. We will construct a prime $\mathbb{R}^+$ coding adapted to $\alpha$ in such a way that $\left(\widehat{A_T}\right)''\left(\hat{k}_0 -\right) = \left(\widehat{A_T}\right)''\left(\hat{k}_0 +\right)\ \forall k_0 \in \left\{5,6,\ldots,\frac{\alpha}{2}-1\right\}$. We would then have constructed the continuous function

$$\mathfrak{G}:\left[\hat{4},\widehat{\alpha}\div\hat{2}\right]\to\mathbb{R}^+,\quad \mathfrak{G}\left(\hat{k}\right)=\left(\widehat{A_T}\right)''\left(\hat{k}\right)$$

For this we select, at random, $0<\xi_2<\xi_3<\xi_4<\xi_5$. According to 3.1.1, $x_4,x_5,\ldots,x_{11}$ are readily determined. We select $\xi_6^2=\frac{x_6}{x_5}\xi_5^2$ , then $\xi_6>\xi_5$, and $x_{12}$ and $x_{13}$ are readily determined. We select $\xi_7>\xi_6$ at random, and $x_{14}$ and $x_{15}$ are readily determined. We now take $\xi_8^2=\frac{x_8}{x_7}\xi_7^2, \xi_9^2=\frac{x_9}{x_8}\xi_8^2, \xi_{10}^2=\frac{x_{10}}{x_9}\xi_9^2$, then $\xi_7<\xi_8<\xi_9<\xi_{10}$, and $x_{16},\ldots,x_{21}$ are readily determined. We select $\xi_{11}>\xi_{10}$ at random, and $x_{22}$ and $x_{23}$ are readily determined. Note that for a prime $i$ we are selecting $\xi_i$ at random with the sole condition $\xi_i>\xi_{i-1}$.

Let $s_0$ be the largest prime such that $s_0\leq\frac{\alpha}{2}-1$. Then, following the same principle, we take $\xi_{s_0}>\xi_{s_0-1}$ at random, and $x_{2s_0}$ and $x_{2s_0+1}$ are readily determined. Finally, we select

$$\text{i) }\xi_{s_0+1}^2=\frac{x_{s_0+1}}{x_{s_0}}\xi_{s_0}^2,\ \xi_{s_0+2}^2=\frac{x_{s_0+2}}{x_{s_0+1}}\xi_{s_0+1}^2,\ \ldots,\ \xi_{\alpha/2-1}^2=\frac{x_{\alpha/2-1}}{x_{\alpha/2-2}}\xi_{\alpha/2-2}^2\quad \left(\text{if } s_0\leq\frac{\alpha}{2}-2\right).$$

$$\text{ii) }\xi_{\alpha/2-1}>\xi_{\alpha/2-2}\ \text{ at random }\left(\text{if } s_0=\frac{\alpha}{2}-1\right).$$

Following 3.1.2, all the essential points $P_{k_0}$ associated with the number $\alpha$ have been determined. We select $\xi_{\alpha/2}^2$ at random and are only have to determine which are to be the remaining coefficients.

i) If $s_0=\frac{\alpha}{2}-1$, we select

$$\xi_{\alpha/2+1}^2=\xi_{\alpha-(\alpha/2-1)}^2=\xi_{\alpha-s_0}^2=\left|y_{s_0-1}\right|\left(\frac{\left|y_{s_0}\right|}{\xi_{\alpha-s_0-1}^2}+\frac{\alpha-s_0}{s_0}\cdot x_{s_0-1}\cdot\left(\frac{1}{\xi_{s_0-1}^2}-\frac{1}{\xi_{s_0}^2}\right)\right)^{-1}$$

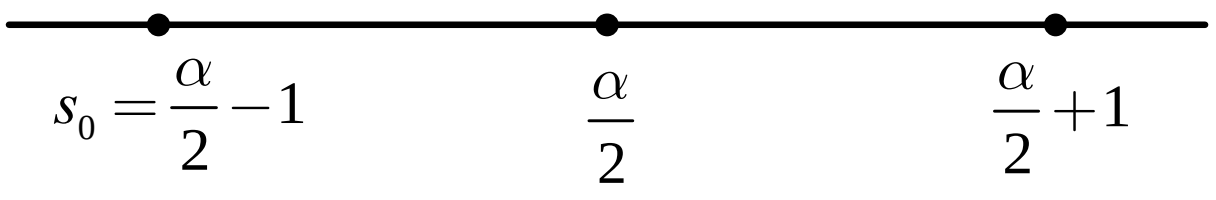

(Fig. 39)

ii) If $s_0<\frac{\alpha}{2}-1$ we select

$$\xi_{\alpha/2+1}^2=\xi_{\alpha-(\alpha/2-1)}^2=\left|y_{\alpha/2-2}\right|\left|y_{\alpha/2-1}\right|^{-1}\xi_{\alpha/2}^2,\quad \xi_{\alpha/2+2}^2=\xi_{\alpha-(\alpha/2-2)}^2=\left|y_{\alpha/2-3}\right|\left|y_{\alpha/2-2}\right|^{-1}\xi_{\alpha/2+1}^2,$$

…

$$\xi^2_{\alpha-s_0-2} = \left|y_{s_0+1}\right|\left|y_{s_0+2}\right|^{-1}\xi^2_{\alpha-s_0-3} \qquad ,\xi^2_{\alpha-s_0-1} = \left|y_{s_0}\right|\left|y_{s_0+1}\right|^{-1}\xi^2_{\alpha-s_0-2}$$

We also verify $\xi^2_{\alpha-s_0-1} = \left|y_{s_0}\right|\left|y_{\alpha/2-1}\right|^{-1}\xi^2_{\alpha/2}$ .

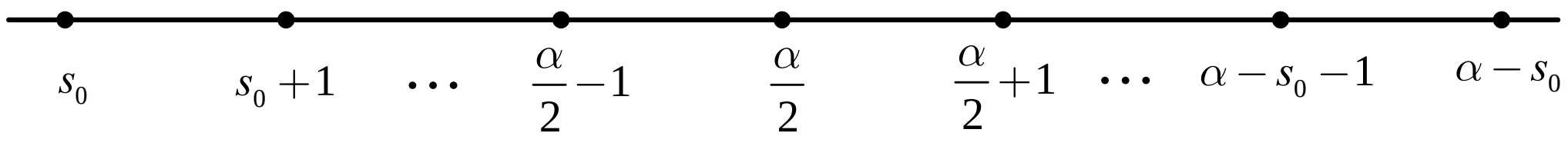

(Fig. 40)

We now take :

$$\xi^2_{\alpha-s_0} = \left|y_{s_0-1}\right|\left(\frac{\left|y_{s_0}\right|}{\xi^2_{\alpha-s_0-1}} + \frac{\alpha-s_0}{s_0}\cdot x_{s_0-1}\cdot\left(\frac{1}{\xi^2_{s_0-1}} - \frac{1}{\xi^2_{s_0}}\right)\right)^{-1}$$

Having selected these first coefficients, we construct the remaining coefficients in the following way:
For each prime $r_0$ where $5 \le r_0 < s_0$ we select

$$\xi^2_{\alpha-r_0} = \left|y_{r_0-1}\right|\left(\frac{\left|y_{r_0}\right|}{\xi^2_{\alpha-r_0-1}} + \frac{\alpha-r_0}{r_0}\cdot x_{r_0-1}\cdot\left(\frac{1}{\xi^2_{r_0-1}} - \frac{1}{\xi^2_{r_0}}\right)\right)^{-1}$$

Between two consecutive primes $p_0$ and $q_0$, such that $5 \le p_0 < q_0 \le s_0$, we select

$$\xi^2_{\alpha-q_0+1} = \left|y_{q_0-2}\right|\left|y_{q_0-1}\right|^{-1}\xi^2_{\alpha-q_0} \ , \ \xi^2_{\alpha-q_0+2} = \left|y_{q_0-3}\right|\left|y_{q_0-2}\right|^{-1}\xi^2_{\alpha-q_0+1}$$

$$\cdots$$

$$\xi^2_{\alpha-p_0-2} = \left|y_{p_0+1}\right|\left|y_{p_0+2}\right|^{-1}\xi^2_{\alpha-p_0-3} \qquad , \ \xi^2_{\alpha-p_0-1} = \left|y_{p_0}\right|\left|y_{p_0+1}\right|^{-1}\xi^2_{\alpha-p_0-2}$$

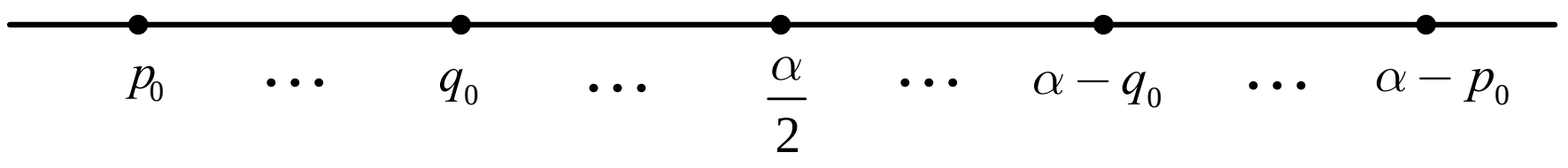

(Fig. 41)

We also verify

$$\xi^2_{\alpha-p_0-1} = \left|y_{p_0}\right|\left|y_{q_0-1}\right|^{-1}\xi^2_{\alpha-q_0}$$

We have now chosen the coefficients $\xi_2,\xi_3,...,\xi_{\alpha/2-1},\xi_{\alpha/2},\xi_{\alpha/2+1},...,\xi_{\alpha-5}$. The remaining coefficients of the prime $\mathbb{R}^+$ coding adapted to $\alpha$ are irrelevant. Due to the actual construction of these coefficients, the hypotheses in proposition 3.1.4. are verified, and we have therefore constructed the following continuous function

$$\mathfrak{G}:\left[\hat{4},\hat{\alpha}\div\hat{2}\right]\to\mathbb{R}^+, \ \ \mathfrak{G}\left(\hat{k}\right) = \left(\widehat{A_T}\right)''\left(\hat{k}\right)$$

**3.1.5 Definition** We will call the "Goldbach Conjecture function associated to $\alpha$" any function $\mathfrak{G}$ constructed in this manner.

**3.1.6 Proposition** Let $\mathfrak{G}$ be a Goldbach Conjecture function with coefficients $\xi_i$ associated to $\alpha$. Let $\mathfrak{P} = \left\{ r_0 : r_0 \text{ is prime, where } 5 \le r_0 \le \frac{\alpha}{2} - 1 \right\}$ and let $s_0$ be the maximum of $\mathfrak{P}$. We call $F_{r_0} = \frac{\alpha - r_0}{r_0} \cdot x_{r_0 - 1} \cdot \left( \frac{1}{\xi^2_{r_0-1}} - \frac{1}{\xi^2_{r_0}} \right)$. Then,

$$\xi^2_{\alpha-5} = |y_4| \left( \left| y_{\alpha/2-1} \right| \xi^{-2}_{\alpha/2} + \sum_{r_0 \in \mathfrak{P}} F_{r_0} \right)^{-1}$$

*Proof* According to the construction of any Goldbach Conjecture function $\mathfrak{G}$, we verify $\xi^2_{\alpha-s_0} = |y_{s_0-1}| \left( \left| y_{\alpha/2-1} \right| \xi^{-2}_{\alpha/2} + F_{s_0} \right)^{-1}$ regardless of the fact that $s_0 = \frac{\alpha}{2} - 1$ or $s_0 < \frac{\alpha}{2} - 1$

We now define the function $\gamma : \mathfrak{P} - \{5\} \to \mathfrak{P} - \{s_0\}$, $\gamma(p)$: as the prime number before $p$. Let $q_0 \in \mathfrak{P} - \{5\}$ and assume that

$$\xi^2_{\alpha-q_0} = |y_{q_0-1}| \left( \left| y_{\alpha/2-1} \right| \xi^{-2}_{\alpha/2} + F_{s_0} + F_{\gamma(s_0)} + F_{\gamma^2(s_0)} + \ldots + F_{\gamma^h(s_0)} \right)^{-1} \quad \text{where } \gamma^h(s_0) = q_0 .$$

Now, let $\gamma(q_0) = p_0$. Thus, due to the construction of the $\mathfrak{G}$ function we verify

$$\xi^2_{\alpha-p_0} = |y_{p_0-1}| \left( |y_{p_0}| \xi^{-2}_{\alpha-p_0-1} + F_{p_0} \right)^{-1} = |y_{p_0-1}| \left( |y_{q_0-1}| \xi^{-2}_{\alpha-q_0} + F_{p_0} \right)^{-1} =$$

$$= |y_{p_0-1}| \left( \left| y_{\alpha/2-1} \right| \xi^{-2}_{\alpha/2} + F_{s_0} + F_{\gamma(s_0)} + \ldots + F_{\gamma^h(s_0)} + F_{\gamma^{h+1}(s_0)} \right)^{-1}$$

As a consequence, and taking $p_0 = 5$, we obtain

$$\xi^2_{\alpha-5} = |y_4| \left( \left| y_{\alpha/2-1} \right| \xi^{-2}_{\alpha/2} + F_5 + F_7 + F_{11} + \ldots + F_{s_0} \right)^{-1} = |y_4| \left( \left| y_{\alpha/2-1} \right| \xi^{-2}_{\alpha/2} + \sum_{r_0 \in \mathfrak{P}} F_{r_0} \right)^{-1}$$

**3.1.7 Simplification of the variables in the $\mathfrak{G}$ function.** Let $\alpha$ be an even number $(\alpha \ge 16)$. Let $\mathfrak{G}$ be any Goldbach Conjecture function associated to $\alpha$. Then, the $\mathfrak{G}$ coefficients can be expressed in the following way:

$$\xi_2^2 > 0 \,, \xi_3^2 = \lambda_3^2 \xi_2^2 \,, \xi_4^2 = \lambda_4^2 \xi_3^2 \,, \xi_5^2 = \lambda_5^2 \xi_4^2 \,, \xi_{p_0} = \lambda^2_{p_0} \xi_{p_0-1} (\forall p_0 \in \mathfrak{P} - \{5\})$$
$$(\xi_2^2 \in (0,+\infty), \lambda_i^2 \in (1,+\infty) \forall i \in \mathfrak{I} = \{3,4\} \cup \mathfrak{P})$$

According to the construction of $\mathfrak{G}$, all the coefficients depend exclusively on the variables $\xi_2^2$, $\lambda_i^2$ and $\xi^2_{\alpha/2} \in (0,+\infty)$. We denote $\overline{\lambda} = (\lambda_i^2)$ $\forall i \in \mathfrak{I}$ thus, any Goldbach conjecture function can be written

$$\mathfrak{G} = \mathfrak{G}(\xi_2^2, \xi^2_{\alpha/2}, \overline{\lambda}) \ , (\xi_2^2 \in (0,+\infty) \, , \ \xi^2_{\alpha/2} \in (0,+\infty) \, , \ \lambda_i^2 \in (1,+\infty) \ \ \forall i \in \mathfrak{I} )$$

**3.1.8 Proposition** Let $\mathfrak{G}=\mathfrak{G}(\xi_2^2,\xi_{\alpha/2}^2,\overline{\lambda})$ be a Goldbach Conjecture function for the even number $\alpha$ $(\alpha \geq 16)$ and the coefficients $\xi_i$. Let us denote for every $p_0 \in \mathfrak{P}-\{5\}$, $P(p_0):=\{p : p \text{ prime } \wedge\ 5 \leq p \leq \gamma(p_0)\}$ Then, $\forall p_0 \in \mathfrak{P}-\{5\}$ we verify

$$\text{i)}\ \frac{x_{p_0}}{\xi_{p_0-1}^2}=\frac{1}{2\lambda_3^2\lambda_4^2}\prod_{j\in P(p_0)}\frac{1}{\lambda_j^2}.\qquad \text{ii)}\ F_5=\frac{\alpha-5}{5}\cdot\frac{1}{2\lambda_3^2\lambda_4^2}\left(1-\frac{1}{\lambda_5^2}\right)$$

$$\text{iii)}\ F_{p_0}=\frac{\alpha-p_0}{p_0}\cdot\frac{1}{2\lambda_3^2\lambda_4^2}\prod_{j\in P(p_0)}\frac{1}{\lambda_j^2}\left(1-\frac{1}{\lambda_{p_0}^2}\right)$$

*Proof* i) The equality is true when $p_0=7$. In fact, according to the construction of $\mathfrak{G}$, we have $\frac{x_7}{\xi_6^2}=\frac{x_6}{\xi_6^2}=\frac{x_5}{\xi_5^2}=\frac{x_4}{\lambda_5^2\xi_4^2}=\frac{\xi_2^2}{2\lambda_5^2\lambda_4^2\lambda_3^2a_2^2}=\frac{1}{2\lambda_3^2\lambda_4^2}\prod_{j\in P(7)}\frac{1}{\lambda_j^2}$. Assume that the equality is true for a prime $p_0 \in \mathfrak{P}-\{5,s_0\}$, we prove that it is also true for the next prime $q_0$. With the actual construction of $\mathfrak{G}$, we obtain:

$$\xi_{q_0-1}^2=\frac{x_{q_0-1}}{x_{p_0}}\xi_{p_0}^2=\frac{x_{q_0}}{x_{p_0}}\xi_{p_0}^2 \Rightarrow \frac{x_{q_0}}{\xi_{q_0-1}^2}=\frac{x_{p_0}}{\xi_{p_0}^2}=\frac{x_{p_0}}{\lambda_{p_0}^2\xi_{p_0-1}^2}=\frac{1}{\lambda_{p_0}^2}\cdot\frac{1}{2\lambda_3^2\lambda_4^2}\prod_{j\in P(p_0)}\frac{1}{\lambda_j^2}=$$

$$=\frac{1}{2\lambda_3^2\lambda_4^2}\prod_{j\in P(q_0)}\frac{1}{\lambda_j^2}.$$

ii) $$F_5=\frac{\alpha-5}{5}x_4\left(\frac{1}{\xi_4^2}-\frac{1}{\xi_5^2}\right)=\frac{\alpha-5}{5}\cdot\frac{x_4}{\xi_4^2}\left(1-\frac{1}{\lambda_5^2}\right)=$$

$$=\frac{\alpha-5}{5}\cdot\frac{1}{2\lambda_3^2\lambda_4^2}\left(1-\frac{1}{\lambda_5^2}\right).$$

iii) $$\forall p_0 \in \mathfrak{P}-\{5\} \Rightarrow F_{p_0}=\frac{\alpha-p_0}{p_0}\cdot x_{p_0-1}\left(\frac{1}{\xi_{p_0-1}^2}-\frac{1}{\xi_{p_0}^2}\right)=\frac{\alpha-p_0}{p_0}\cdot\frac{x_{p_0}}{\xi_{p_0-1}}\left(1-\frac{1}{\lambda_{p_0}^2}\right)=$$

$$=\frac{\alpha-p_0}{p_0}\cdot\frac{x_{p_0}}{\xi_{p_0-1}^2}\left(1-\frac{1}{\lambda_{p_0}^2}\right)=\frac{\alpha-p_0}{p_0}\cdot\frac{1}{2\lambda_3^2\lambda_4^2}\prod_{j\in P(p_0)}\frac{1}{\lambda_j^2}\left(1-\frac{1}{\lambda_{p_0}^2}\right).$$

**3.1.9 Example.** We construct the elements that intervene in any Goldbach Conjecture function $\mathfrak{G}$ where $\alpha=18$. In this case, $\frac{\alpha}{2}=9$, $\frac{\alpha}{2}-1=8$, $\frac{\alpha}{2}-3=6$. Then, the coefficients $\xi_2^2,\xi_3^2,\xi_4^2,\xi_5^2$ can be thus expressed:

$$\xi_2^2>0,\ \xi_3^2=\lambda_3^2\xi_2^2,\ \xi_4^2=\lambda_3^2\lambda_4^2\xi_2^2,\ \xi_5^2=\lambda_3^2\lambda_4^2\lambda_5^2\xi_2^2.$$

Then, $x_4,x_5,...,x_{11}$ are readily determined. Using 2.8.1 we obtain the expressions of $x_i(i\ natural\ number\ ,\ 4\leq i\leq 11)$.

$$x_4=\frac{1}{2}\xi_2^2=x_5(5prime)$$

$$x_6=\xi_2\xi_3-\frac{1}{2}\xi_2^2=(\lambda_3-\frac{1}{2})\xi_2^2=x_7(7prime)$$

$$x_8=\xi_2\xi_4-\frac{1}{2}\xi_2^2=(\lambda_4\lambda_3-\frac{1}{2})\xi_2^2$$

$$x_9=|y_8|=\xi_2\xi_4-\xi_2\xi_3+\frac{1}{2}\xi_3^2=(\lambda_4\lambda_3-\lambda_3+\frac{1}{2}\lambda_3^2)\xi_2^2$$

$$x_{10}=|y_7|=\xi_2\xi_5-\xi_2\xi_3+\frac{1}{2}\xi_3^2=(\lambda_5\lambda_4\lambda_3-\lambda_3+\frac{1}{2}\lambda_3^2)\xi_2^2=x_{11}=|y_6|(11prime)$$

$$\xi_6^2=\frac{x_6}{x_5}\xi_5^2=2(\lambda_3-\frac{1}{2})\lambda_5^2\lambda_4^2\lambda_3^2\xi_2^2 .$$

Now $x_{12}$ and $x_{13}$ are readily determined:

$$x_{12}=|y_5|=\xi_2\xi_6-\xi_2\xi_4+\xi_3\xi_4-\frac{1}{2}\xi_3^2=(\sqrt{2(\lambda_3-\frac{1}{2})}\cdot\lambda_5\lambda_4\lambda_3-\lambda_4\lambda_3+\lambda_4\lambda_3^2-\frac{1}{2}\lambda_3^2)\xi_2^2=$$

$$=x_{13}=|y_4|(13prime).$$

$$\xi_7^2=\lambda_7^2\xi_6^2,\ \xi_8^2=\frac{x_8}{x_7}\xi_7^2$$

Also,

$$F_5=\frac{\alpha-5}{5}\cdot\frac{1}{2\lambda_3^2\lambda_4^2}\left(1-\frac{1}{\lambda_5^2}\right), F_7=\frac{\alpha-7}{7}\cdot\frac{1}{2\lambda_3^2\lambda_4^2\lambda_5^2}\left(1-\frac{1}{\lambda_7^2}\right).$$

Choosing at random $\xi_9^2>0$, the remaining coefficients are readily determined:

$$\xi_{10}^2=\xi_{\alpha-8}^2=\frac{|y_7|}{|y_8|}\xi_9^2$$

$$\xi_{11}^2=\xi_{\alpha-7}^2=|y_6|\left(|y_7|\xi_{\alpha-8}^{-2}+F_7\right)^{-1}=|y_6|\left(|y_8|\xi_9^{-2}+F_7\right)^{-1}$$

$$\xi_{12}^2=\xi_{\alpha-6}^2=\frac{|y_5|}{|y_6|}\xi_{\alpha-7}^2=|y_5|\left(|y_8|\xi_9^{-2}+F_7\right)^{-1}$$

$$\xi_{13}^2=\xi_{\alpha-5}^2=|y_4|\left(|y_5|\xi_{\alpha-6}^{-2}+F_5\right)^{-1}=|y_4|\left(|y_8|\xi_9^{-2}+F_5+F_7\right)^{-1}$$

For $\xi_2^2\in(0,+\infty),\ \lambda_j^2\in(1,+\infty),\ \xi_9^2\in(0,+\infty)$ we obtain all the Goldbach Conjecture functions $\mathfrak{G}$ associated to $\alpha=18$: $\mathfrak{G}=\mathfrak{G}(\xi_2^2,\xi_{\alpha/2}^2,\overline{\lambda})$.

**3.1.10 Note** It is convenient to denote:

$$\mu_4=\lambda_3^2\lambda_4^2,\ \mu_5=\lambda_3^2\lambda_4^2\lambda_5^2,\ \mu_{p_0}=\lambda_3^2\lambda_4^2\prod_{j\in P(p_0)\cup\{p_0\}}\lambda_j^2\ (\ \forall p_0\in\mathfrak{P}-\{5\}).$$

We also define $\gamma(5)=4$.

**3.1.11 Proposition.** Let $\alpha$ be an even number ($\alpha\geq16$), and $\mathfrak{G}=\mathfrak{G}(\xi_2^2,\xi_{\alpha/2}^2,\overline{\lambda})$ be any Goldbach conjecture function associated to $\alpha$, $n(\alpha):=card(\mathfrak{I})$ ($\mathfrak{I}=\{3,4\}\cup\mathfrak{P}$).
Then, these functions exist:

$$f_i,\ g_j,\ h_j\ :\ (1,+\infty)^{n(\alpha)}\to\mathbb{R}\,(i\in\mathbb{N},\ j\in\mathbb{N},2\leq i\leq\frac{\alpha}{2}-3,4\leq j\leq\frac{\alpha}{2}-1$$

such that

$$\text{i) }\xi_i^2=f_i(\overline{\lambda})\xi_2^2\qquad\text{ii) }x_j=g_j(\overline{\lambda})\xi_2^2\qquad\text{iii) }|y_j|=h_j(\overline{\lambda})\xi_2^2$$

*Proof* Considering that $\left|y_j\right| = x_{\alpha-j-1}$ ( $\forall j\in\mathbb{N}: 4\le j\le \frac{\alpha}{2}-1$ ) is sufficient to prove that these functions exist:

$$f_i,\ g_j\ :\ (1,+\infty)^{n(\alpha)} \to \mathbb{R}\,(i\in\mathbb{N},\ \ j\in\mathbb{N}, 2\le i\le \frac{\alpha}{2}-3, 4\le j\le \alpha-5)$$

such that

$$i')\ \xi_i^2 = f_i(\overline{\lambda})\xi_2^2, \qquad ii')\ \ x_j = g_j(\overline{\lambda})\xi_2^2$$

then we would choose $h_j = g_{\alpha-j-1}$ $(j\in\mathbb{N}\ 4\le j\le \frac{\alpha}{2}-1)$.

Following 3.1.9, i) and ii) are true for the natural numbers *i, j* where $2\le i\le 5$, $4\le j\le 11$, that is, i) and ii) are true for every $\xi_i^2, x_j$ naturally associated to the prime $p_0$ =5. Now, regardless of 3.1.9 , we prove that i) and ii) are true for the natural numbers *i, j* where $6\le i\le 7$ , $12\le j\le 13$ .

In fact,

$$\xi_6^2 = \frac{x_6}{x_5}\xi_5^2 = \frac{g_6(\overline{\lambda})}{g_5(\overline{\lambda})} f_5(\overline{\lambda})\xi_2^2 = f_6(\overline{\lambda})\xi_2^2 \text{ ( if we define } f_6 = \frac{g_6}{g_5} f_5 \text{ ).}$$

The addends that appear in $x_{12}$ and $x_{13}$ have the form $\pm\,\xi_l\xi_k$ or $\pm\,\frac{1}{2}\xi_h^2$ ( *l, h, k,* natural numbers where $2\le l\le 6$, $2\le k\le 6$, $2\le h\le 6$), that is, we have addends of the form $\pm\sqrt{f_l(\overline{\lambda}) f_k(\overline{\lambda})\xi_2^2}$ or $\pm\frac{1}{2} f_h(\overline{\lambda})\xi_2^2$. Thus, $x_{12}$ and $x_{13}$ can be written

$$x_{12} = g_{12}(\overline{\lambda})\xi_2^2,\ \ x_{13} = g_{13}(\overline{\lambda})\xi_2^2$$

Now, $\xi_7^2 = \lambda_7^2\xi_6^2 = \lambda_7^2 f_6(\overline{\lambda})\xi_2^2 = f_7(\overline{\lambda})\xi_2^2$ ( if we define $f_7 = \lambda_7^2 f_6$ ) then, $x_{14}$ and $x_{15}$ are readily determined and their addends have the form $\pm\,\xi_l\xi_k$ or $\pm\,\frac{1}{2}\xi_h^2$ ( *l, h, k,* natural numbers where $2\le l\le 7$, $2\le k\le 7$, $2\le h\le 7$) . Following the reasoning stated above, $x_{14}$ and $x_{15}$ can be expressed $x_{14} = g_{14}(\overline{\lambda})\xi_2^2, x_{15} = g_{15}(\overline{\lambda})\xi_2^2$. We now consider the prime $p_0$ (where $7 < p_0 \le s_0$ ). Following the previous outline we easily prove that if $i')$ and $ii')$ are true for every *i, j* natural numbers (where $2\le i\le \gamma(p_0)$, $4\le j\le 2\gamma(p_0)+1$) then $i')$ and $ii')$ are also true for every *i,j* natural numbers where $2\le i\le p_0$, $4\le j\le 2p_o+1$.

**3.1.12 Corollary** Let $\alpha$ be an even number $(\alpha\ge 16)$ and $\mathfrak{G} = \mathfrak{G}(\xi_2^2, \xi_{\alpha/2}^2, \overline{\lambda})$ any Goldbach Conjecture function associated to $\alpha$, then $\left|x_{k_0-1}\right|\cdot\left|x_{k_0}\right|^{-1}$ and $\left|y_{k_0-1}\right|\cdot\left|y_{k_0}\right|^{-1}$ do not depend on $\xi_2^2$ $(\forall k_0$ natural number, $5\le k_0\le \frac{\alpha}{2}-1)$.

**3.1.13 Definition** Let $\alpha$ be an even number $(\alpha\ge 16)$ and $\mathfrak{G} = \mathfrak{G}(\xi_2^2, \xi_{\alpha/2}^2, \overline{\lambda})$ be any Goldbach conjecture function associated to $\alpha$. If $\lambda_i = u\in(1,+\infty)\ \forall i\in\{3,4\}\cup\mathfrak{P}$, we say that $\mathfrak{G}$ is a scalar Goldbach conjecture function associated to $\alpha$. We denote such a function by $\mathfrak{G}=\mathfrak{G}(\xi_2^2,\xi_{\alpha/2}^2,u)$.

**3.1.14 Proposition** Let $\alpha$ be an even number ($\alpha\geq 16$) and $\mathfrak{G}=\mathfrak{G}(\xi_2^2,\xi_{\alpha/2}^2,u)$ any scalar Goldbach conjecture function associated to $\alpha$. Then, $\forall k_0\in\mathbb{N}(4\leq k_0\leq\alpha-5)$ we verify $\lim_{u\to 1^+} x_{k_0}=\frac{1}{2}\xi_2^2$.

*Proof* Choosing $\xi_3^2=\lambda_3^2\xi_2^2=u^2\xi_2^2,\xi_4^2=\lambda_3^2\lambda_4^2\xi_2^2=u^4\xi_2^2,\xi_5^2=\lambda_3^2\lambda_4^2\lambda_5^2\xi_2^2=u^6\xi_2^2$ we readily determine $x_4,x_5,...,x_{11}$.

Following 3.1.9,

$$x_4=x_5=\frac{1}{2}\xi_2^2\ \ x_6=x_7=(u-\frac{1}{2})\xi_2^2,\quad x_8=(u^2-\frac{1}{2})\xi_2^2$$

$$x_9=(\frac{3}{2}u^2-u)\xi_2^2\ \ x_{10}=x_{11}=(u^3-u+\frac{1}{2}u^2)\xi_2^2$$

and therefore we verify $\lim_{u\to 1^+} x_{k_0}=\frac{1}{2}\xi_2^2\quad(\forall k_0\in\mathbb{N},4\leq k_0\leq 11)$.

Now, $\xi_6^2=\frac{x_6}{x_5}\xi_5^2$ thus we verify

$$\lim_{u\to 1^+}\xi_6^2=\lim_{u\to 1^+}\frac{x_6}{x_5}u^6\xi_2^2=\xi_2^2.$$

We have readily determined $x_{12}$ and $x_{13}$. Following 2.8.1, for any Goldbach conjecture function $\mathfrak{G}=\mathfrak{G}(\xi_2^2,\xi_{\alpha/2}^2,u)$ and for every natural number $k_0$ ($4\leq k_0\leq\alpha-5$) the expression of $x_{k_0}$ is

$$x_{k_0}=\xi_2(\xi_{\left[\frac{k_0}{2}\right]}-\xi_{\left[\frac{k_0}{3}\right]})+\xi_3(\xi_{\left[\frac{k_0}{3}\right]}-\xi_{\left[\frac{k_0}{4}\right]})+...+\xi_{[\sqrt{k_0}]-1}(\xi_{\left[\frac{k_0}{[\sqrt{k_0}]-1}\right]}-\xi_{\left[\frac{k_0}{[\sqrt{k_0}]}\right]})+\frac{1}{2}\xi_{[\sqrt{k_0}]}^2$$

$$(\text{if } [\sqrt{k_0}]=\left[\frac{k_0}{[\sqrt{k_0}]}\right])$$

$$x_{k_0}=\xi_2(\xi_{\left[\frac{k_0}{2}\right]}-\xi_{\left[\frac{k_0}{3}\right]})+\xi_3(\xi_{\left[\frac{k_0}{3}\right]}-\xi_{\left[\frac{k_0}{4}\right]})+...+\xi_{[\sqrt{k_0}]}(\xi_{\left[\frac{k_0}{[\sqrt{k_0}]}\right]}-\frac{1}{2}\xi_{[\sqrt{k_0}]})\quad(\text{ if } [\sqrt{k_0}]<\left[\frac{k_0}{[\sqrt{k_0}]}\right])$$

Considering that $\lim_{u\to 1^+}\xi_i^2=\xi_2^2$ ($i=2,3,4,5,6$) we conclude that

$$\lim_{u\to 1^+}x_{12}=\lim_{u\to 1^+}x_{13}=\frac{1}{2}\xi_2^2$$

Now, $\xi_7^2=\lambda_7^2\xi_6^2=u^2\xi_6^2$ thus, $\lim_{u\to 1^+}\xi_7^2=\xi_2^2$. According to the construction of $\mathfrak{G}$ and by a simple induction process we obtain

$$\lim_{u\to 1^+}x_{k_0}=\frac{1}{2}\xi_2^2\ \ (\forall k_0\in\mathbb{N},4\leq k_0\leq\alpha-5).$$

**3.1.15 Corollary** Since $\left|y_j\right| = x_{\alpha-j-1}$ $(\forall j \in \mathbb{N}, 4 \le j \le \frac{\alpha}{2} - 1)$, the following is true

$$\lim_{u\to1^+} x_{k_0} = \lim_{u\to1^+}\left|y_{k_0}\right| = \frac{1}{2}\xi_2^2 \ (\forall k_0 \in \mathbb{N}, 4 \le k_0 \le \frac{\alpha}{2} - 1).$$

## 3.2 Time and Goldbach Conjecture

**3.2.1 Proposition** The following set is infinite.

$$\mathfrak{N} = \left\{\alpha \in \mathbb{N} : (\alpha \ \textit{even}) \wedge (\alpha \ge 16) \wedge (\frac{\alpha}{2} \ \textit{non-prime}) \ \wedge (\alpha - 3 \ \textit{non-prime})\right\}$$

*Proof* Consider $\mathfrak{N}_1 = \{12k : k \in \mathbb{N} \ (k \ge 2)\}$, obviously $12k$ is even, $\alpha \ge 16$, ${}^{12k}\!/_{2} = 6k$ is not a prime and $12k - 3 = 3(4k - 1)$ is not a prime either. Thus $\mathfrak{N}_1 \subset \mathfrak{N}$ and $\mathfrak{N}_1$ is infinite. As a consequence, $\mathfrak{N}$ is an infinite set.

**3.2.2 $\psi^{-1}$ as a movement. Conclusion** Let $\psi : \mathbb{R}^+ \to [0, M_\psi)$ be an $\mathbb{R}^+$ coding function, $\alpha$ an even number $(\alpha \ge 16)$, consider the restriction to $[4, \alpha]$, that is,

$$\psi : [4, \alpha] \to [\widehat{4}, \widehat{\alpha}]$$

If we consider the variable $\widehat{k}$ as time and the variable $k$ as space, then $\psi^{-1}$ represents the movement of a particle in the closed interval $[4, \alpha]$. Suppose that we have chosen $\psi$ in such a way that the function $\mathfrak{G} = \mathfrak{G}(\xi_2^2, \xi_{\alpha/2}^2, u)$ is a scalar Goldbach Conjecture function associated to $\alpha$. Then, $\mathfrak{G}$ is the acceleration of the $\widehat{A_T}$ area. Following to 3.1.15, if we choose $\xi_2^2 = \xi_{\alpha/2}^2 = 1$ we verify:

$$\lim_{u\to1^+} x_{k_0} = \frac{1}{2} \text{ and } \lim_{u\to1^+} y_{k_0} = -\frac{1}{2} \ (\forall k_0 \in \mathbb{N} : 4 \le k_0 \le \frac{\alpha}{2} - 1)$$

Also, considering the construction of $\mathfrak{G}$ : $\lim_{u\to1^+} \xi_i^2 = 1$ $(\forall i \in \mathbb{N} : \ 2 \le i \le \frac{\alpha}{2} - 1)$.

This means that in the limit position $\psi = \psi^{-1} = I : [4, \alpha] \to [4, \alpha]$ ( $I$ identity function ) and in this limit position the essential points have been transformed into

$$P_{k_0} = (x_{k_0}, y_{k_0}) = \left(\frac{1}{2}, -\frac{1}{2}\right) \quad (\forall k_0 \in \mathbb{N} : 4 \le k_0 \le \frac{\alpha}{2} - 1)$$

In other words, the characterization 2.7.3 about the fact $\alpha \in \mathfrak{N}$ being the sum of two prime numbers has been lost. This leads to the following conclusion:

"*There exists at least one characterization of the Goldbach Conjecture in an infinite set of even numbers that depends on time*".

REFERENCES

This paper is not based on, nor continues previously published works. Consequently, it is not possible to make a list of references. We only wish to reference, in a general manner, all those people who have provided the mathematical concepts and tools used in this paper.